\documentclass[11pt]{article}
\usepackage{amsfonts}
\usepackage{amssymb}
\usepackage{amsmath}
\usepackage{graphicx}
\usepackage{hyperref}
\usepackage{rotating}
\usepackage{multirow}
\usepackage{subcaption}
\usepackage{changepage}
\usepackage{xcolor}
\newtheorem{theorem}{Theorem}
\newtheorem{lemma}{Lemma}

\newtheorem{assumption}{Assumption}

\begin{document}

\title{Asymptotically Optimal Control of a Centralized Dynamic Matching
Market with General Utilities}
\author{Jose H. Blanchet\thanks{Management Science and Engineering Department, Stanford University, Stanford, CA 94305, jose.blanchet@stanford.edu}, Martin I. Reiman\thanks{
Department of Industrial Engineering and Operations Research, Columbia
University, New York, NY 10027, martyreiman@gmail.com}, Virag Shah\thanks{Management Science and Engineering Department, Stanford University, Stanford, CA 94305, virag@stanford.edu}, Lawrence M. Wein\thanks{
Graduate School of Business, Stanford University, Stanford, CA 94305,
lwein@stanford.edu, corresponding author}, Linjia Wu\thanks{Management Science and Engineering Department, Stanford University, Stanford, CA 94305, linjiawu@stanford.edu}}
\maketitle

\begin{abstract}
\setlength{\baselineskip}{12pt plus2pt minus1pt}
We consider a matching market where buyers and sellers arrive
according to independent Poisson processes at the same rate and independently abandon the market if not matched after an exponential amount of time with the same mean. In this centralized market, the utility for the system manager from
matching any buyer and any seller is a general random variable. We consider a sequence of systems indexed by $n$ where the arrivals in the $n^{\mathrm{th}}$ system are sped up by a factor of $n$. We analyze two families of one-parameter policies: the population threshold policy immediately matches an arriving agent to its best available mate only if the number of mates in the system is above a threshold, and the utility threshold policy matches an arriving agent to its best available mate only if the corresponding utility is above a threshold. Using an asymptotic fluid analysis of the two-dimensional Markov process of buyers and sellers, we show that when the matching utility distribution is light-tailed, 
the population threshold policy with threshold $\frac{n}{\ln n}$ is asymptotically optimal among all policies that make matches only at agent arrival epochs. In the heavy-tailed case
, we characterize the optimal threshold level for both policies. 
We also study the utility threshold policy in an unbalanced matching market with heavy-tailed matching utilities, and find that the buyers and sellers have the same asymptotically optimal utility threshold. To illustrate our theoretical results, we use extreme value theory to derive optimal thresholds when the matching utility distribution is exponential, uniform, Pareto, and correlated Pareto. In general, we find that as the right tail of the matching utility distribution gets heavier, the threshold level of each policy (and hence market thickness) increases, as does the magnitude by which the utility threshold policy outperforms the population threshold policy. 

\vskip 0.4truecm

\end{abstract}

\section{Introduction}

\label{sec-intro}

We consider a symmetric centralized dynamic matching
market (the asymmetric case is also discussed for heavy-tailed utilities). 
Two types of agents, which we call buyers and sellers, arrive to the
market according to independent Poisson processes with rate $\lambda$, and each agent abandons
(i.e., exits) the market after an independent exponential amount of time with rate $\eta$ if he has not
yet been matched. The utility of a match between any buyer and any seller is
a general random variable. In this
centralized model, the agents make no explicit decisions, and at the time of
an agent arrival, the system manager observes all matching utilities between
the arrival and all potential mates (e.g., sellers if the arrival is a
buyer) who are currently in the market. Using information about the number
of buyers and sellers and their matching utilities, the system manager
decides when to make matches and which agents to match.

Centralized dynamic matching markets occur in settings such as organ
transplants, public housing, labor markets and various online platforms. In
practice, matching utilities include information about tissue type matching
and the geographical distance between the donor and the recipient for organ
transplants; the location and desirability of the residence and the distance
between the residence and the applicant's current residence in public housing;
and the match between the needs of the employer and the experience and
skills of the job applicant in the labor market. This information can lead
to wide variations in the matching utilities between different buyers and
sellers, and our goal is to understand how best to exploit this variation
when managing the market. However, in our idealized model, the details
about this information are suppressed (e.g., we do not use covariates describing the agents to
help make decisions) and aggregated into the matching utility distribution
between buyers and sellers.

A key issue in centralized dynamic matching markets is to find the optimal
market thickness; i.e., rather than match a new agent upon its arrival, it
may be preferable to place the arriving agent in the market and allow more
agents to arrive in the hope of making a higher-utility match in the future.
In our model, we aim to maximize the long-run expected average utility rate
(i.e., utility of matches per unit time) of all matches. Although we do not explicitly 
include agent waiting costs, a strategy that forces agents to wait too long
for the market to thicken can backfire because agents may abandon the market
before they are matched.

Due to the challenging nature of this problem, we resort to asymptotic
methods. We consider a sequence of systems where the arrival rates in the $%
n^{\mathrm{th}}$ system are multiplied by $n>0$. In the absence of any
matching, the number of agents of each type would be precisely the number of
customers in a $M/M/\infty$ queue, which would be $O(n)$ (a generic function $f(n)$ is $O(n)$ if $\limsup_{n\to\infty}\frac{f(n)}{n}\le c$ for some finite constant $c>0$). We use two types
of asymptotic methods: one is a fluid analysis of the two-dimensional
Markov process for the number of buyers and sellers in the market when the
arrival rates are large. The other is extreme value theory (Gumbel 1958,
Galambos 1978) and regularly varying functions (Resnick 1987), which are used because the utility of a match under the
policies we consider is the maximum of a (typically) large number of random variables. In our study, a fluid analysis of the queueing process is sufficient to derive our results, and leads to a decoupling of the extremal behavior of the utilities and the dynamics of the queueing system. This decoupling in turn allows us to consider correlated utilities, which is a feature that is lacking in other dynamic matching models. 

In this asymptotic regime, we compute an upper bound on the utility rate of
any policy that makes matches only at agent arrival epochs, and compare it to the utility rate of two families of threshold policies: the population threshold policy and the utility threshold policy. Under the population threshold policy, the system manager immediately matches an arriving agent to
the available mate with the highest matching utility (at which point, the
arriving agent and its matched mate exit the system and their matching
utility is collected by the system manager) only if the number of available
mates in the market exceeds a specified threshold; otherwise, the arriving
agent is not immediately matched and is instead placed in the market. Under the utility threshold policy, the arriving agent is immediately matched to its best available mate only if the corresponding matching utility exceeds a specified threshold. 

Although possibly not optimal among all policies, these single-parameter policies are easy to implement and describe, and allow for quite explicit results. In fact, the population threshold policy can be implemented without ever calculating the utility of individual matches (although the probability distribution of matches is required to compute the optimal threshold): all that is required is a ranked ordering of the possible matches. As discussed below, the utility threshold policy outperforms the population threshold policy in our examples, but the latter policy is asymptotically optimal in certain cases. Another natural class of policies to consider is a batching policy, where the system manager -- after a certain amount of time or after a certain number of buyers and/or sellers collect in the market -- matches a set of agents. This approach requires an optimization algorithm to perform the matching and hence is more computationally demanding than our two threshold policies. Moreover, if there are many agents who abandon quickly after arrival, as in some call centers (e.g., Fig.~20 in Gans \textsl{et al.} 2003), a batching policy may not be very robust in practice. Nonetheless, in~\S\ref{sec-batch} we consider a batch-and-match policy that  periodically (with an asymptotically optimal time window) optimally matches all agents on the thinner side of the market with an equal number of agents randomly selected from the thicker side of the market. 

\subsection{Preview of Results}

\label{ssec-preview} In extreme value theory, the limiting
distribution of the maximum of many random variables can be one of three
types, loosely based on whether the underlying distribution of these random
variables has an exponential right tail, has a heavier (e.g., power law)
right tail, or is bounded from above, and our results are qualitatively
different in each case. Although our main results are couched in terms of regularly varying functions, we preview our results with three canonical examples (Table~\ref{table1}) -- one from each of the three domains of attraction -- in the symmetric case, which are analyzed in~\S %
\ref{sec-examples} in the Appendix. When matching utilities have an
exponential distribution, the population threshold policy with a threshold of $\frac{n}{%
\ln n}$ is asymptotically optimal with a utility rate that is $O(n\ln n)$
and twice as large as the utility rate of the greedy policy -- i.e., the population threshold policy with a
threshold of zero -- in the limit.
When the matching utilities have a Pareto ($c,\beta)$ distribution with shape parameter $%
\beta>1$ (and hence a finite mean), the  population threshold  $\frac{\lambda}{%
\eta(1+\beta)}n$ is asymptotically optimal. Although the utility rate of this threshold policy does not
converge to the loose upper bound in this case, the utility rate and the
upper bound are both $O(n^{1+1/\beta})$, whereas the utility rate under the
greedy policy is only $O(n^{1+1/(2\beta)})$. When the matching utilities
have a uniform distribution, the greedy policy is asymptotically optimal (i.e., 0 is an asymptotically optimal population threshold) and
the optimal utility rate is $O(n)$. 

In the Pareto case, the  utility threshold  $0.763\sqrt{\pi n}$ is asymptotically optimal when $c=1$ and $\beta=2$, and the corresponding utility rate is $O(n^{1+1/\beta})$. 
This asymptotic utility rate is computed explicitly and it is shown to be larger than the utility rate of the asymptotically optimal population threshold policy. In the exponential and uniform cases, where we have already identified an asymptotically optimal policy, we use heuristics to 
compute, in the pre-limit, utility threshold policies that are consistent with the asymptotically optimal descriptions, but outperform in simulation results the population threshold policy. 
We also consider a positively correlated Pareto case in~\S\ref{ssec-examples-correlated}, and show that the asymptotically optimal population threshold is independent of the correlation, the asymptotically optimal utility threshold decreases as the correlation increases, and the utility rates of both threshold policies decrease  as the correlation increases. In~\S\ref{sec-unbalanced}, we consider an unbalanced market, where buyers have a different arrival rate and abandonment rate than sellers, and  analyze the utility threshold policy in the heavy-tailed case. Surprisingly, although we allow the buyers and sellers to have a different utility  threshold, we find that they have the same asymptotically optimal utility threshold. Finally, we show in~\S\ref{sec-batch} that in the Pareto case, the utility threshold policy outperforms the batch-and-match policy. 

Taken together, the optimal amount of
patience -- and hence market thickness - increases with the right tail of
the matching utility distribution, as does the optimal utility rate and the performance gap between the utility threshold policy and the population threshold policy. Our limited analysis of an unbalanced market suggests that the optimal market thickness also increases with the amount of imbalance. In our particular model of correlation, increased positive correlation among matching utilities decreases the benefit of increased patience (i.e., the system manager is less likely to observe a future utility that is much better than the best existing utility), while the cost of increased patience (i.e., the number of abandonments) is independent of the correlation. Among the three one-parameter policies considered here, the utility threshold policy displays the best performance.

\subsection{Related Work}

\label{ssec-lit-review} Matching markets is a large and active area of
research, and we restrict our review to centralized dynamic markets.
Although our model lacks the contextual richness of some of the models for
specific types of markets, the most distinctive feature of our model is the
general matching utilities, which allows us to understand how the right tail of
the matching utility distribution impacts the optimal thickness of the
market (Table~\ref{table1}). In contrast, much of the recent work in dynamic
(centralized or decentralized) matching markets, either via two-type agents
(e.g., easy-to-match or hard-to-match agents, or matches that are preferred or non-preferred, Baccara \textsl{et al.} 2015, Ashlagi \textsl{et
al.} 2019a, Ashlagi \textsl{et al.} 2019b) or a compatibility network (Ashlagi \textsl{et al.} 2013, Anderson \textsl{et al.} 2017, Akbarpour \textsl{et al.} 2019, Varma \textsl{et al.} 2019) essentially
lead to dichotomous outcomes for a match. Exceptions include $\ddot {\mathrm{U}}
$nver (2010), who considers blood type compatibility for a dynamic kidney
exchange model, Emek \textsl{et al.} 2016 and Ashlagi \textsl{et al.} 2017a, who consider minimizing mismatch costs when agents arrive on a finite metric space in a non-bipartite and bipartite setting, respectively, and Ashlagi \textsl{et al.} 2018, who allow general matching utilities in a discrete time model with a constant time until abandonment. They perform a primal-dual analysis to derive competitive ratios for algorithms when there is no prior information about the match values or arrival times. 

The analysis of multiclass matching queues is an active area. Hu and Zhou (2016)
consider a discrete-time, multiclass, discounted variant of our problem that
includes waiting costs. They show that the optimal policy is of threshold
form under vertical and unidirectionally horizontal differentiated types. Ding \textsl{et al.}
(2016) allow the matching utilities to depend on the class of buyer and
seller, and performs a fluid analysis of a greedy policy, and B$\breve {%
\mathrm{u}}$si$\acute {\mathrm{c}}$ and Meyn (2016) minimize linear holding
costs in a system without class-dependent matching utilities or abandonment,
but also find that matches are not made until there are a sufficient number
of agents in the market; see Moyal and Perry (2017), where these systems are
referred to as matching queues, for other references to these types of
models. 

Gurvich and Ward (2014) and Nazari and Stolyar (2019) study a control problem in a more general
setting than the studies above, where arriving customers wait to be matched to agents of other
classes. Gurvich and Ward (2014) minimize cumulative holding costs over a finite horizon, and show that a myopic discrete-review matching algorithm is asymptotically optimal. Nazari and Stolyar (2019) maximize the long-run average revenue rate subject to maintaining stable queues, and construct a greedy primal-dual approach that is asymptotically optimal. It is difficult to compare this powerful result to our results, given that we assume abandonment rather than stability, and we have a single-class model with general rewards rather than a multiclass model with class-dependent rewards. 

Two other studies consider fluid and diffusion limits of simplified
versions of our model where either a match occurs with a certain probability
for each buyer-seller pair (B$\ddot {\mathrm{u}}$ke and Chen 2017) or
everyone matches when there is an available mate (Liu \textsl{et al.} 2015), which corresponds to our
greedy policy, but with a deterministic utility (i.e., a matching utility
distribution that is a point mass at one value). In both cases, the system
state reduces to a one-dimensional quantity (the number of sellers minus the
number of buyers), whereas our model requires a two-dimensional state space
for a non-greedy policy. 

Perhaps the most closely related paper is Mertikopoulos \textsl{et al.} (2020), which also considers a symmetric centralized dynamic matching market. Compared to our study, they assume independent exponential mismatch costs rather than general matching utilities, and consider waiting times rather than abandonment, and are interested in minimizing the sum of mismatch and waiting costs over a finite horizon. They consider a class of policies that make the $k^{\rm th}$ match (which has the lowest mismatch cost among possible matches) when the short side of the market grows to a certain one-parameter function of $k$. They analyze the performance of the policy (using the celebrated $\pi^2/6$ result for the expected minimum weight matching due to Mezard and Parisi (1987) and rigorously proved by Aldous (2001)) under various values of the parameter, and also identify a policy that balances the mismatch and waiting costs. It is difficult to draw qualitative comparisons between our results for exponential utilities (which incorporate abandonments) and their results (which incorporate waiting costs); indeed, our approach depends on the right tail of the exponential distribution via extreme value theory, whereas their approach depends on the left tail of the exponential distribution via minimum weighted matching.

We briefly mention other work that is only peripherally related. Originally
motivated by public housing (Kaplan 1988), Caldentey \textsl{et al.} (2009)
and Adan and Weiss (2012) consider infinite bipartite matching of servers
and customers under the first-come first-served policy. There is also
a stream of work in online bipartite matching in an adversarial setting
(Karp \textsl{et al.} 1990), where agents do not wait in the market if they
are not matched immediately. Finally, there is a body of literature (e.g., Duffie \textsl{et al.} 2018 and references therein) that uses the law of large numbers to analyze the performance of static and dynamic matching models used in economics, finance and genetics, but these models are descriptive rather than prescriptive. 

\subsection{Organization}

The paper is organized as follows.
We formulate the model in~\S\ref{sec-model} and state our main theoretical results in~\S\ref{sec-results}, which are proved in~\S\ref{sec-proofs}. After analyzing a greedy policy in~\S\ref{sec-greedy}, we apply our main results to specific matching utility distributions in~\S\ref{sec-examples} in the Appendix and assess the accuracy of these results in a simulation study in~\S\ref{sec-simulation}. The unbalanced case is studied in~\S\ref{sec-unbalanced}, the batch-and-match policy is analyzed in~\S\ref{sec-batch}, and concluding remarks are offered in~\S\ref{sec-conclusion}. 

\subsection{Notation}

For the convenience of the reader we collect together the notational conventions used in this paper.
Although we have already introduced the notation $O(n)$ above, we repeat it here:
a generic function $f(n)$ is $O(n)$ if $\limsup_{n\to\infty}\frac{f(n)}{n}\le c$ for some finite constant $c>0$.
In a similar vein, we introduce  $o(n)$, $\Omega(n)$, and $\Theta(n)$.
A generic function $f(n)$ is $o(n)$ if $\lim_{n\to\infty}\frac{f(n)}{n}=0$,
is $\Omega(n)$ if there exist $c>0$ and an integer $n_o\ge 1$ such that $f(n)\ge cn$ for all integers $n\ge n_o$, and 
is $\Theta(n)$ if $f(n)$ is both $O(n)$ and $\Omega(n)$.
We use $x_n \sim y_n$ as shorthand for $\frac{x_n}{y_n}\to 1$ as $n\to\infty$. 

We let $\mathbb{R}$ denote the real line, and, for any finite integer $k \ge 1$, let $\mathbb{R}^{k}$ denote the $k$-dimensional Euclidean  space. The Euclidean norm of $x\in \mathbb{R}^{k}$ is denoted by $|x|$.
We let $\mathbb{R}_+$ denote the set of nonnegative reals, and $\mathbb{Z}_+$ denote the set of nonnegative integers.
The stochastic processes we consider take values in $\mathbb{R}^{k}$, and are assumed to be elements of $\mathbb{D}^{k} [0, \infty)$, the space of right continuous functions mapping $[0,\infty)$ into $\mathbb{R}^{k}$ that have left limits, endowed with the Skorokhod topology.

For $x \in \mathbb{R}_+, \lceil x \rceil$ is the smallest integer that is not smaller than $x$.
The standard stochastic order between two distribution functions $F_1$ and $F_2$ is denoted by $F_1 \le_{st} F_2$.
We use $\stackrel{d}{=}$ to denote equality in distribution. More specifically, we write $X\stackrel{d}{=}$ Poisson$(x)$ to denote that the random variable $X$ has a Poisson distribution with mean $x$.

\section{The Model\label{sec-model}}
\textsl{Dynamics.} Buyers and sellers arrive to the market according to independent Poisson processes with rate $\lambda$. The agents are impatient, in that each buyer and each seller independently abandons the market after an independent and identically distributed (i.i.d.) exponential amount of time with rate $\eta$ if they are not matched within this time.  If an agent is matched prior to his abandonment then the agent leaves at the time of matching.  

Let $B(t)$ and $S(t)$ be the number of buyers and sellers in the system at 
time $t$; these agents have arrived but have not yet abandoned or been matched. The utility of a match between any buyer and any seller is a random variable $V\ge 0$ with cumulative distribution function (CDF) $F(v)$. When a buyer (seller, respectively) arrives to this centralized system to find it in state $(B(t),S(t))$ then $S(t)$ ($B(t)$, respectively) instances of $V$ are observed by the system manager, which represent the matching utilities of the arriving agent with all currently available potential mates. Thus, at any point in time the system manager knows the utility that would be generated by matching any buyer to any seller. 

\textsl{Policies.} Our goal is to maximize the long run expected average rate of utility from matches, which we refer to as the utility rate. While the system manager could conceivably make matches at any point in time, we restrict our attention to \textsl{arrival-only policies}, where a match may occur only at the arrival epoch of one of the agents being matched. In particular, we consider the following two classes of arrival-only policies. 

\begin{enumerate}
\item Population threshold policies: A buyer who arrives at time $t$ is matched immediately to a seller if the number of sellers in the system satisfies $S(t) \ge z$; in this case, the arriving buyer is matched to the seller who has the highest matching utility with the buyer, with ties broken arbitrarily. If $S(t) < z$, then the arriving buyer waits in the market, and leaves upon being matched to a later-arriving seller or upon abandonment. Similarly, a seller who arrives at time $t$ is immediately matched to the highest-matching buyer if $B(t) \ge z$, and waits otherwise. We refer to the parameter $z$ as the population threshold.  

\item Utility threshold  policies: A buyer who arrives at time $t$ is matched immediately to the seller with matching value $\max_{1\le i\le S(t)}V_i$ if $\max_{1\le i\le S(t)} V_i>v$ for some fixed $v\ge 0$, with ties broken arbitrarily. If $\max_{1\le i\le S(t)} V_i\le v$, then the 
arriving buyer waits in the market, and leaves upon being matched to a later-arriving seller or upon abandonment. Similarly, a seller who arrives at time $t$ is immediately matched to the buyer with matching value $\max_{1\le i\le B(t)}V_i$ if $\max_{1\le i\le B(t)} V_i>v$, and waits otherwise. We refer to the parameter $v$ as the utility threshold. 

\end{enumerate}

Given the symmetry of the underlying stochastic model, it seems natural to restrict ourselves to single-parameter policies, where buyers and sellers have the same threshold ($z$ or $v$). In the analysis of the unbalanced case in~\S7, we allow different utility thresholds for buyers and sellers ($v_b$ and $v_s$) and find that the asymptotically optimal values satisfy $v_b=v_s$ under Pareto matching utilities. This result suggests that a single-parameter threshold policy is not only easier to use in practice and easier to analyze than a two-parameter threshold policy, but also does not sacrifice performance.

\textsl{Utilities.} In our model, the utilities of potential matches of a new arrival with  agents on the other side of the market may be correlated. However, we make the following assumption. 

\begin{assumption}\label{as:UtilityDependence}
There exists a sequence of distributions $F_1, F_2, \ldots$ such that $F_{k-1} \le_{st} F_k$ for each $k$ and, for an arriving agent who finds $k$ agents on the other side of the market, $\max\{V_1,\ldots,V_k\}$ is independent of the past, and has distribution $F_k$. 
\end{assumption}


For example, if the utilities of different matches are i.i.d. with distribution $F$, then 
$F_{k}\left( x\right) = \left( F\left( x\right) \right)^k$ in Assumption~\ref{as:UtilityDependence}. But Assumption 1 allows us to deal with correlated utilities, which is natural when there is contextual information (e.g., covariates) that can be used to inform the utilities based on the types of buyers and sellers to be matched. More specifically, Assumption 1 holds if the utilities are conditionally independent given a context observed at the time of arrival. In this case, the equation $F_{k}\left( x\right) = \left( F\left( x\right) \right)^k$ holds with an additional expectation, and the stochastic ordering in $k$ still holds.  


Let the random variable $M(k) \triangleq \max\{V_1,\ldots,V_k\}$ have distribution $F_k$. We impose the following assumption on $M(k)$.

\begin{assumption} \label{as:Regularly_Varying} 
For each $x \in \mathbb R_+$, define 
\[
m\left( x\right) =E\left( M\left( \left\lceil x\right\rceil \right) \right) , 
\]%
and suppose that $m\left( \cdot \right) $ is regularly varying with index $%
\alpha \in \lbrack 0,1)$. That is, for every $x>0$, 
\begin{equation}
\lim_{t\rightarrow \infty }\frac{m\left( tx\right) }{m\left( t\right) }%
=x^{\alpha }.  \label{RV_max}
\end{equation}%
A regularly varying function with index $\alpha =0$ is also known as slowly
varying.
\end{assumption}

For the case of i.i.d. utilities, Assumption \ref{as:Regularly_Varying} covers every utility distribution such that $E\left( V^{1+\delta }\right) <\infty $ for some $\delta >0$. 
All distributions that belong to the maximum domain of
attraction of a generalized extreme value distribution -- which unifies the Type I (Gumbel), Type II (Frechet) and Type III (Weibull) laws within a single parametric family -- satisfy \eqref{RV_max} (including, e.g., uniform, beta, gamma, lognormal, Pareto).  
There are also other distributions that do not belong to any domain of
attraction in extreme value theory for which \eqref{RV_max} holds; e.g., the geometric, negative binomial, and Poisson distributions satisfy \eqref{RV_max} with $\alpha =0$. For ease of reference, we collect some basic facts about extreme value theory and regularly varying functions in~\S\ref{sec-evt}. 

The case $\alpha =0$ corresponds to distributions for which all
moments exist (i.e., the tail of $V$ decays faster than any polynomial),
whereas $\alpha >0$ corresponds to the case in which the tails of $%
V$ decrease roughly like a polynomial with degree $1/\alpha $. The
condition that $\alpha <1$ is imposed to guarantee that $E\left( V^{1+\delta
}\right) <\infty $ for some $\delta >0$. We will refer to $\alpha
=0 $ as the light-tailed case and $%
\alpha \in \left( 0,1\right) $ as the 
heavy-tailed case.

\textsl{Scaling.} To make further progress, we consider a sequence of systems indexed by $n=1,2,\ldots$, and some quantities in the $n^{\rm th}$ system include the subscript $n$. The arrival rate in the $n^{\rm th}$ system is $n\lambda$, and the abandonment rate in the $n^{\rm th}$ system is $\eta$. Alternatively and equivalently, we could leave the arrival rate unscaled and slow down the abandonment rate by a factor of $n$, as in Liu \textsl{et al.} (2015). The matching utilities are unscaled. In the $n^{\rm th}$ system, we denote the system state by $(B_n(t),S_n(t))$, the population threshold by $z_n$, the utility threshold by $v_n$, and the utility rate by $U_n$.

\section{Main Results}\label{sec-results}

Results for the population threshold policy and the utility threshold policy are given in Theorem~\ref{thm:Population_Based} in~\S\ref{ssec-results-population} and in Theorem~\ref{thm:utility_based} in~\S\ref{ssec-results-utility}, respectively. Theorem~\ref{thm:Population_Based} shows that the optimal population threshold policy is asymptotically optimal among the class of arrival-only policies when $\alpha=0$, and provides the asymptotically optimal population threshold when $\alpha\in(0,1)$.  Theorem~\ref{thm:utility_based} provides the asymptotically optimal utility threshold when $\alpha\in (0,1)$. The proofs of Theorems~\ref{thm:Population_Based} and~\ref{thm:utility_based} appear in~\S\ref{sec-proofs}. 

\subsection{Population Threshold Policy}\label{ssec-results-population}

%
%



We begin by providing a dynamic description of the system using Poisson processes. Denote the indicator function of event $x$ by $I_{\{x\}}$ and let
$N_B^+\left( \cdot \right)
, N_B^-\left( \cdot \right) ,N_S^+\left(\cdot \right)
,N_S^-\left( \cdot \right) $
be independent
Poisson processes with unit rate, which are used to construct buyer arrivals, buyer abandonments, seller arrivals and seller abandonments, respectively. Under the population threshold policy with threshold $z_n$, the state $(B_n,S_n)$ of the $n^{\rm th}$ system at time $t$ satisfies%
\begin{align}
B_{n}\left(  t\right)    & =B_{n}\left(  0\right)  +\int_{0}^{t}I_{\{
S_{n}\left(  r_{-}\right)  <z_{n}\}}  dN_{B}^{+}\left(  \lambda nr\right)
-N_{B}^{-}\left(  \eta\int_{0}^{t}B_{n}\left(  r\right)  dr\right)
\nonumber\\
& -\int_{0}^{t}I_{\{B_{n}\left(  r_{-}\right)  \geq z_{n}\}}
dN_{S}^{+}\left(  \lambda nr\right)  ,\label{eq1b}\\
S_{n}\left(  t\right)    & =S_{n}\left(  0\right)  +\int_{0}^{t}I_{\{
B_{n}\left(  r_{-}\right)  <z_{n}\}}  dN_{S}^{+}\left(  \lambda nr\right)
-N_{S}^{-}\left(  \eta\int_{0}^{t}S_{n}\left(  r\right)  dr\right)
\nonumber\\
& -\int_{0}^{t}I_{\{S_{n}\left(  r_{-}\right)  \geq z_{n}\}}
dN_{B}^{+}\left(  \lambda nr\right)  .\label{eq2b}%
\end{align}


The process $\{B_n(t), S_n(t), t\ge 0\}$ is a non-negative (entry wise) irreducible two-dimensional birth-and-death process on a subset of $\mathbb{Z}_+\times \mathbb{Z}_+$ and each coordinate is bounded by that of an infinite-server queue, for each $n>0$. Thus, the process $\{B_n(\cdot), S_n(\cdot)\}$ is a positive-recurrent continuous-time Markov chain and therefore it possesses a stationary distribution, which we denote by $(B_n(\infty), S_n(\infty))$. By symmetry  the utility rate $U^p_n(z_n)$ of the population threshold policy with threshold $z_n$ can be expressed as 
\begin{eqnarray}
U_{n}^p(z_n) &=&\lambda nE\left[m\left( B_{n}\left( \infty \right) \right) I_{\{
B_{n}\left( \infty \right) \geq z_{n}\}} \right]  + \lambda nE\left[m\left( S_{n}\left( \infty \right) \right) I_{\{
S_{n}\left( \infty \right) \geq z_{n}\}} \right], \nonumber \\
 &=& 2\lambda nE\left[m\left( B_{n}\left( \infty \right) \right) I_{\{
B_{n}\left( \infty \right) \geq z_{n}\}} \right] . \label{eq3b}
\end{eqnarray}

The theorem below shows that for $\alpha = 0$, the population threshold policy is asymptotically optimal among the family of arrival-only policies. Also, for each $\alpha\in [0,1)$, it characterizes the scaling of the optimal population threshold, i.e., the threshold $z_n$ that maximizes the utility rate asymptotically as $n\to \infty$. 
\begin{theorem}\label{thm:Population_Based}

Suppose that Assumption \ref{as:Regularly_Varying} holds.

i) If $\alpha =0$ then there exists an $o(n)$ sequence of population thresholds $z_n^*$ such that $\lim_{n\to \infty} \frac{m(n)}{m(z_n^*)} = 1$. For any such sequence of thresholds, the population threshold policy is
asymptotically optimal in the following sense. Let $U_n^p(z_n^*)$ and $U_n$ be the utility rates under the above policy and any other arrival-only policy, respectively. Then $\lim\inf_{n\to \infty} \frac{U_n^p(z_n^*)}{U_n} \ge 1$. The associated utility rate satisfies
\begin{equation}
\lim_{n\to \infty} \frac{U_n^p(z_n^*)}{nm(n)}  = \lambda.
\label{eq4b}
\end{equation}

ii) If $\alpha \in \left( 0,1\right) $ then  the population threshold policy with $z_{n}^*=z_{\ast }n$ where $z_{\ast
} = \frac{\lambda \alpha}{\eta (1+\alpha)}$ is asymptotically
optimal among the class of population threshold policies. The associated utility rate satisfies
\begin{equation}
\lim_{n\to \infty} \frac{U_n^p(z_n^*)}{n m(n)} = \lambda z_{\ast }^\alpha \left(1 - \frac{\eta z_{\ast }}{\lambda} \right).
\label{eq5b}
\end{equation}

\end{theorem}

For $\alpha=0$, it remains to compute an $o(n)$ sequence of thresholds $z_n^*$ such that $\lim_{n\to \infty} \frac{m(n)}{m(z_n^*)} = 1$. This is usually not difficult to do. For example, when utilities are i.i.d.\ with an exponential distribution, then $z_n^* = \frac{n}{\ln n}$ satisfies this property. More generally, as shown in Theorem 1 in Bojanic and Seneta (1971), for a large class of distributions, setting $z_n^* = \frac{n}{m(n)^\delta}$ for any positive real $\delta$ is sufficient. By setting $z_{n}$ in this way (i.e., $o(n)$ but not too small), we  simultaneously ensure the following: (1) the fraction of agents that abandon the system tends to $0$, and (2) the market thickness, i.e. $B_n(\infty)$, is almost linear in $n$. In other words, almost all agents experience maximal utility. This can be seen most clearly in equation~(\ref{eq4b}), where the utility rate under the optimal population threshold policy satisfies $U_n^p(z_n^*)\sim n\lambda m(n)$, which is the arrival rate of buyers times the expected value of the maximum of $n$ matching utilities. 

However, for heavy-tailed distributions in part ii) of Theorem~\ref{thm:Population_Based}, $m(z_n)$ for any $o(n)$ sequence $z_n$ is vanishingly small compared to $m(n)$. Thus, it is not possible to ensure that most users see maximal utility, implying that our simple upper bound is unachievable. Moreover, to maximize the utility rate, it is not obvious whether the system manager should set $z_n=o(n)$ to guarantee that most agents are matched instantly, or should set $z_n = O(n)$ to ensure that market thickness is maximal even if a nontrivial fraction of users abandon the system.  Part ii) of Theorem~\ref{thm:Population_Based} implies that the latter option is the right choice under heavy-tailed distributions. 

We conclude this subsection with a brief sketch of the proof of Theorem~\ref{thm:Population_Based}, which relies on a fluid analysis of equations~(\ref{eq1b})-(\ref{eq2b}). We define $\bar B_n(t)=n^{-1}B_n(t)$ and $\bar S_n(t)=n^{-1}S_n(t)$. Because the formal limit of $(B_n(t),S_n(t))$ involves indicator functions that are not continuous (see~(\ref{eq64a})-(\ref{eq64b}) in~\S\ref{ssec-proofs-population} in the Appendix), we need to study the limiting dynamical system as the solution to the following Skorokhod problem:
\begin{align}
\bar{B}\left(  t\right)    & =\bar{B}\left(  0\right)  +\lambda t-\eta\int%
_{0}^{t}\bar{B}\left(  r\right)  dr-L_{z}^{\bar{B}}\left(  t\right)
-L_{z}^{\bar{S}}\left(  t\right)  , \label{eq62aa-mt}\\
\bar{S}\left(  t\right)    & =\bar{S}\left(  0\right)  +\lambda t-\eta\int%
_{0}^{t}\bar{S}\left(  r\right)  dr-L_{z}^{\bar{B}}\left(  t\right)
-L_{z}^{\bar{S}}\left(  t\right)  , \label{eq62bb-mt}
\end{align}
where $L_{z}^{\bar{B}}\left(  \cdot\right)  ,$ $L_{z}^{\bar{S}}\left(
\cdot\right)  $ are nondecreasing processes such that $L_{z}^{\bar{B}}\left(
0\right)  =L_{z}^{\bar{S}}\left(  0\right)  =0$ and%
\begin{equation}
\int_{0}^{t}\left(  \bar{B}\left(  r\right)  -z\right)  dL_{z}^{\bar{B}%
}\left(  r\right)  =\int_{0}^{t}\left(  \bar{S}\left(  r\right)  -z\right)
dL_{z}^{\bar{S}}\left(  r\right) =0 , \label{eq62cc-mt}
\end{equation}
and $\bar{B}\left(  t\right)  ,\bar{S}\left(  t\right)  \leq z$. 
To obtain explicit expressions for the Skorokhod problem, we use the change of variables $\bar{B}_{z}\left(  t\right)  =z-\bar{B}\left(  t\right)$, $\bar{S}_{z}\left(  t\right)  =z-\bar{S}\left(  t\right) $ and 
$\bar{\lambda}_{z}=\lambda/\eta-z\geq0$. This allows us to reduce~(\ref{eq62cc-mt}) to the one-dimensional condition
\[
\int_{0}^{t}\min\left(  \bar{B}_{z}\left(  r\right)  ,\bar{S}_{z}\left(
r\right)  \right)  dL\left(  r\right)  =0,\text{ \ }L\left(  0\right)  =0,
\]
which enables us to obtain an explicit solution to~(\ref{eq62aa-mt})-(\ref{eq62cc-mt}). With this solution in hand, we show uniqueness and then apply a standard Picard iteration to argue existence. 

We use martingale arguments to show that $(\bar S_n(\cdot),\bar B_n(\cdot))\to (\bar S(\cdot),\bar B(\cdot))$ uniformly on compact sets in probability. The dynamical system describing $(\bar B,
\bar S)$ has the unique attractor $(z,z)$ if $\lambda/\eta\ge z$, given the initial condition $\bar B(0)\le z$, $\bar S(0)\le z$. We then show that the limit interchange $(t\to \infty$ and $n\to\infty$) holds, and prove that $(\bar B_n(\infty),\bar S_n(\infty))\to (z,z)$ almost surely as $n\to\infty$. 

The next step in the proof is to compute the utility rate. Taking expectations on both sides of equation~(\ref{eq1b}) yields 
\begin{equation}
\label{eq62a-mt}
\eta E\left( \bar{B}_{n}\left( \infty \right) \right) =\lambda \{P\left( 
\bar{S}_{n}\left( \infty \right) <z\right) -P\left( \bar{B}_{n}\left( \infty
\right) \geq z\right) \},
\end{equation}
from which we can obtain, using symmetry arguments, that 
\begin{equation} \label{eq:Fluid_Stationary-mt}
\lim_{n\rightarrow \infty }P\left( \bar{B}_{n}\left( \infty \right) \geq
z\right) =\frac{1}{2}\left( 1-\frac{\eta z}{\lambda }\right) . 
\end{equation}%

The following key lemma, which is proved in~\S\ref{sec-proofs}, allows us to compute the utility rate in (\ref{eq3b}). Recall that $m(n)$ is defined in Assumption~\ref{as:Regularly_Varying}.
\begin{lemma}\label{lemma:Expectation_Max}
 Let $\left\{ N_{n}\right\} _{n\geq 1}$ be a sequence of
positive random variables taking values on the positive integers and let $%
\bar{N}_{n}=E\left( N_{n}\right) <\infty $. Assume that $\bar{N}%
_{n}\rightarrow \infty $, and that $P\left( \left\vert N_{n}-\bar{N}%
_{n}\right\vert >\varepsilon \bar{N}_{n}\right) \rightarrow 0$. Then $E\left[
m\left( N_{n}\right) \right] \sim m\left( \bar{N}_{n}\right) $ as $%
n\rightarrow \infty $.
\end{lemma}

Using Lemma~\ref{lemma:Expectation_Max} and equation~(\ref{eq:Fluid_Stationary-mt}) and setting $z_n=nz$ allows us to compute the utility rate
\begin{equation}
U_{n}^p(z_n) = \lambda nm\left( zn\right) \left( 1-\frac{\eta z}{\lambda }\right) \left(
1+o\left( 1\right) \right) \label{eq70a-mt}
\end{equation}%
as $n\rightarrow \infty $, and combining~(\ref{eq70a-mt}) with equation~(\ref{RV_max}) yields  
\begin{equation} \label{eq:Utility_Limit-mt}
\frac{U_n^p(z_n)}{n m(n)} = \lambda z^\alpha \left( 1-\frac{\eta z}{\lambda }\right) \left(
1+o\left( 1\right) \right).
\end{equation}

In the $\alpha\in(0,1)$ case, we optimize the right side of~(\ref{eq:Utility_Limit-mt}) with respect to $z$ to obtain the asymptotically optimal population threshold $z_{n}^*=z_{\ast }n$, where $z_{\ast
} = \frac{\lambda \alpha}{\eta (1+\alpha)}$.

In the $\alpha=0$ case, we use the fluid limit analysis similar to above to show that $E[B_n(\infty)]$ is $o(n)$ for any sequence of thresholds $z_n$ that is $o(n)$. Furthermore, the arguments used to obtain~(\ref{eq:Fluid_Stationary-mt}) also imply that  
\begin{equation}\label{eq:ProbAlpha0-mt}
\lim_{n\rightarrow \infty }P\left( {B}_{n}\left( \infty \right) \geq
z_n \right) = \lim_{n\rightarrow \infty }P\left( {S}_{n}\left( \infty \right) \geq
z_n \right) = \frac{1}{2}.
\end{equation}
A PASTA (Poisson Arrivals See Time Averages) argument implies that $U_n^p(z_n) \ge \lambda n m(z_n) \left( 1 + o(1) \right).$ Consequently, for any sequence $z_n=o(n)$ such that $$\lim_{n\to \infty} \frac{m(z_n)}{m(n)} = 1,$$ we would have that  $U_n^p(z_n) \ge \lambda n m(n) \left( 1 + o(1) \right)$. Lemma~\ref{lemma:VerySlowlyVarying} in~\S\ref{ssec-proofs-lemmas} in the Appendix guarantees that such a sequence exists. 

Finally, asymptotic optimality in part i) of Theorem~\ref{thm:Population_Based} follows from the above results by constructing the following simple upper bound (see~\S\ref{sec-proofs} for a proof of Lemma~\ref{lemma:bound}) on the performance of any arrival-only policy, which uses Lemma~\ref{lemma:Expectation_Max} and assumes that all agents are matched (and hence the arrival rate in Lemma~\ref{lemma:bound} is $\lambda n$) and that -- when computing $B_n(\infty)$ in equation~(\ref{eq3b}) -- agents leave only upon abandonment (implying that $B_n(\infty)\stackrel{d}{=} {\rm Poisson}(\lambda n/\eta))$.

\begin{lemma}\label{lemma:bound}
Let $U_n$ be the utility rate for any arrival-only policy. Then an upper bound $U_n^+$ is given by 
$$U_{n}\leq U_n^+ = \lambda nm\left( \frac{\lambda n}{\eta}\right).$$
\end{lemma}

%
%



\subsection{Utility Threshold Policy}
\label{ssec-results-utility}


Because the population threshold policy is asymptotically optimal within the class of arrival-only policies when $\alpha=0$, we focus on the case $\alpha\in (0,1)$ in Theorem~\ref{thm:utility_based}. In order to describe the dynamics of the utility threshold policy, we
introduce two independent arrays of nonnegative i.i.d. random variables,
$\left\{  V_{i,j}^B:i\geq1,j\geq1\right\}  $ and $\left\{  V_{i,j}^{S
}:i\geq1,j\geq1\right\}  $ having CDF $F\left(
\cdot\right)  $. We let  $\left\{  A_{j}^B:j\geq1\right\}  $ be the sequence of arrival times associated with the process 
$N_B^+(n\lambda \cdot)$ and $\left\{  A_{j}^S:j\geq1\right\}  $ be the sequence of arrival times associated with the process
$N_S^+(n\lambda \cdot)$. The dynamics can be described path-by-path as follows:
\begin{align}
B_n\left(  t\right)   &  =B_n\left(  0\right)  +\sum
_{j=1}^{N_B^+\left(  n\lambda t\right)  }I_{\{\max_{i=1}^{S_n\left(
A^B_{j-}\right)  }V_{i,j}^B\leq v\}}  -\sum_{j=1}^{N_S^+\left(n\lambda 
t\right)  }I_{\{\max_{i=1}^{B_n(A^S_{j-})}V_{i,j}%
^{S}>v\}} \label{MP0}\\
&  -N_B^-\left( \eta \int_{0}^{t}B_n\left(  r_{-}\right)  dr\right)
,\nonumber\\
S_n\left(  t\right)   &  =S_n\left(  0\right)  +\sum
_{j=1}^{N_S^+\left( n\lambda t\right)  }I_{\{\max_{i=1}^{B_n%
(A_{j-}^S)}V_{i,j}^{S}\leq v\}}  -\sum_{j=1}^{N_B^+\left(
n\lambda t\right)  }I_{\{\max_{i=1}^{S_n(A_{j-}^B)}V_{i,j}^B%
>v\}} \nonumber\\
&  -N_S^-\left( \eta \int_{0}^{t}S_n\left(  r\right)
dr\right)  .\nonumber
\end{align}

By symmetry and ergodicity, we can express the utility rate $U_n^u(v_n)$ for the utility threshold policy with threshold $v_n$ as 
\begin{equation}
U_{n}^u(v_n) = 2 \lambda nE\left[E[M(B_{n}(\infty)) I_{\{
M\left( B_{n}\left( \infty \right) \right) \geq v_{n}\}} |B_{n}\left(
\infty \right)]\right]. 
\label{eq6b}
\end{equation}

Because the analysis of the utility threshold policy considers the entire distribution of the maximum rather than only its expected value, we need to strengthen Assumption \ref{as:Regularly_Varying} by imposing the following additional assumption.

\begin{assumption}\label{as:Max_Pareto} In addition to Assumption \ref{as:Regularly_Varying}, suppose that $\alpha \in \left( 0,1\right) $ and
\[
\frac{M(n)}{m(n)} \Rightarrow X ~~{\rm as} ~ n\rightarrow \infty,
\]%
where $P\left( X>t\right) =1-e^{-\kappa /t^{1/\alpha
}}$  and $\kappa $ is a
normalizing constant such that $E\left( X\right) =1$.
\end{assumption}

That is, $X=(\kappa^{-1}T)^{-\alpha}$ is an exponential random variable with mean one. Assumption \ref{as:Max_Pareto} is satisfied if the utilities belong to the domain of attraction of the Frechet law, which in turn is equivalent, in the i.i.d. case, to requiring the distribution of utilities to be regularly varying with index $1/\alpha$ (see Section 1.2, Proposition 1.11 of Resnick 1987).


\begin{theorem}
\label{thm:utility_based}

Suppose that Assumption \ref{as:Max_Pareto} holds. For
$x\in\lbrack0,\lambda/\eta]$, define 
\[
v\left(  x\right)  = \left(\frac{\kappa x}{\ln\left(\frac{2\lambda}{\eta x+\lambda}\right)}\right)^{\alpha}.
\]
Then there exists a unique solution $x_{\ast}\in\left(  0,\lambda/\eta\right)  $
satisfying%
\[
x_{\ast}^{1-\alpha}v\left(  x_{\ast}\right)  \frac{\eta}{2\lambda\alpha\kappa^{\alpha}}%
=\int_{0}^{\kappa x_{\ast}/v\left(  x_{\ast}\right)  ^{1/\alpha}%
}t^{-\alpha}e^{-t} dt.
\]
Moreover, a threshold policy with utility threshold $v_{n}^*=v\left(  x_{\ast}\right)  m\left(
n\right)  $ is asymptotically optimal among the class of utility threshold policies and the associated utility rate satisfies
\[
\lim_{n\rightarrow\infty}\frac{U_{n}^u(v_n^*)}{nm(n)}=2\lambda x_{\ast}^{\alpha
}E\left[  XI_{\{X\geq\frac{v(x_{\ast})}{x_{\ast}^{\alpha}}\}}  \right] .
\]

\end{theorem}

As in the population threshold policy, the above result shows that for heavy-tailed distributions it is beneficial to ensure that market thickness is maximal at the cost of abandonment of a nontrivial fraction of users in the system. Although we do not prove any results for the utility threshold policy in the $\alpha=0$ case (since asymptotic optimality is already achieved for the population threshold policy), we show in~\S\ref{sec-examples} how heuristics inspired by Theorems~\ref{thm:Population_Based} and~\ref{thm:utility_based} can lead to effective utility thresholds in the $\alpha=0$ case.

The proof of Theorem~\ref{thm:utility_based} uses the same general approach as in the proof of Theorem \ref{thm:Population_Based} part ii), and we briefly outline it here. We assume that the thresholds satisfy 
$$\frac{v_{n}}{m(n)} \to v \text{ for some } v \ge 0,$$
and use Assumption~\ref{as:Max_Pareto} to show that the putative fluid limit of $B_n(t)=n^{-1}B_n(t)$ and $S_n(t)=n^{-1}S_n(t)$ is
\begin{equation}
    \label{fl1}
\bar{B}\left( t\right) =\bar{B}\left( 0\right)  +\lambda \int_{0}^{t} e^{-\kappa 
\bar{S}\left( r\right) /v^{1/\alpha }} -\eta \int_{0}^{t}%
\bar{B}\left( r\right) dr - \lambda \int_{0}^{t} \left( 1 -e^{-\kappa \bar{B}%
\left( r\right) /v^{1/\alpha }}\right)  dr, 
\end{equation}
\begin{equation}
\bar{S}\left( t\right) = \bar{S}\left( 0\right) +\lambda \int_{0}^{t}e^{-\kappa 
\bar{B}\left( r\right) /v^{1/\alpha }} -\eta \int_{0}^{t} 
\bar{S}\left( r\right) dr  - \lambda \int_{0}^{t} \left(1-e^{-\kappa \bar{S}%
\left( r\right) /v^{1/\alpha }} \right) dr.
\label{fl2}
\end{equation}
A martingale decomposition similar to that given in the proof of part ii) of Theorem~\ref{thm:Population_Based} shows that $\bar{B}_{n}\left( \cdot \right) \rightarrow \bar{B}\left( \cdot \right)$ and $\bar{S}_{n}\left( \cdot \right) \rightarrow \bar{S}%
\left( \cdot \right)$  uniformly on compact sets in probability. Because~(\ref{fl1})-(\ref{fl2}) does not pose the degeneracies
involving the Skorokhod map encountered in the case of Theorem \ref{thm:Population_Based} part ii), we can use Theorem 7.2 of Chapter 3 in Ethier and Kurtz (2005) to show that the family $\left\{\left((\bar{B}_n(t): t\ge 0),(\bar{S}_n(t):t\ge 0)\right)\right\}_{n\ge 1}$ is tight in the Skorokhod topology. 

The unique solution to the fluid limit satisfies
$$
0 = - \eta\bar{x} - \lambda+ 2 \lambda e^{-\kappa\bar{x} /
v^{1/\alpha}},
$$
which can be expressed as 
\begin{equation}
v\left(  \bar{x}\right) = \left(\frac{\kappa \bar x}{\ln\left(\frac{2\lambda}{\eta \bar x+\lambda}\right)}\right)^{\alpha}.
\label{Eq_v_zbar-mt}%
\end{equation}
or
$$
\bar{x}\left(  v\right)  =-\frac{\lambda}{\eta}+\frac{v^{1/\alpha}%
}{\kappa}W\left(  \frac{2\lambda\kappa}{\eta v^{1/\alpha}}\exp\left(
\frac{\lambda\kappa}{\eta v^{1/\alpha}}\right)  \right)  , 
$$
where $W(x)$ is the Lambert W function. 

We use Assumptions~\ref{as:Regularly_Varying} and \ref{as:Max_Pareto} to optimize the utility rate with respect to $\bar x(v)$, yielding the optimization problem
\begin{equation}
\sup_{\bar{x}\in\left(  0,\lambda/\eta\right)  }2\lambda\bar{x}^{\alpha
}\kappa^\alpha\int_{0}^{\kappa\bar{x}/v\left(  \bar{x}\right)  ^{1/\alpha}%
}t^{-\alpha}e^{-t} dt.
\label{Opt_z-mt}%
\end{equation}
The solution to~(\ref{Opt_z-mt}) reduces to $\bar x^*$ uniquely satisfying
$$
 \frac{v\left(  \bar{x}\right) \eta}{2\lambda\kappa^{\alpha}}=\alpha
\bar{x}^{\alpha-1}\int_{0}^{\kappa\bar{x}/v\left(  \bar{x}\right)  ^{1/\alpha
}}t^{-\alpha}e^{-t}dt,
$$
and substituting $\bar x^*$ into~(\ref{Eq_v_zbar-mt}) gives the optimal utility threshold.


\bigskip

\section{A Greedy Policy}

\label{sec-greedy} In~\S\ref{sec-examples} in the Appendix, we apply the results in Theorems~\ref{thm:Population_Based} and~\ref{thm:utility_based} to several different matching utility distributions, and then assess the accuracy of these analyses via simulation in~\S\ref{sec-simulation}. To provide a natural benchmark for comparison, we first analyze the greedy policy, which corresponds to the population threshold policy with threshold $z_n=0$. That is, under the greedy policy, each arriving agent is
matched to the available mate with the highest matching utility, and waits
in the market if there are no available mates.

Under the greedy policy, the state of the $n^{\mathrm{th}}$ system can
be described by $B_{n}(t)-S_{n}(t)$ because there are never both buyers and
sellers in the system at the same time. By Theorem~4.5 in Liu \textsl{et al.}
(2015), the steady-state distribution of $\frac{B_{n}(t)-S_{n}(t)}{\sqrt{n}}$
converges to $N(0,\lambda /\eta )$ as $n\rightarrow \infty $.

The probability that a buyer or seller abandons is the long-run expected
number of abandonments per unit time divided by the total arrival rate of
agents (i.e., buyers plus sellers), which can be approximated by 
\begin{eqnarray}
\frac{\sqrt{n}\eta E[|N(0,\lambda /\eta )|]}{2\lambda n} &=&\frac{\sqrt{n}%
\eta \sqrt{\frac{2}{\pi }}\frac{\lambda }{\eta }}{2\lambda n},  \nonumber \\
&=&\frac{1}{\sqrt{2\pi n}},  \label{eq13a} \\
&\rightarrow &0.  \label{eq14}
\end{eqnarray}%
By~(\ref{eq14}), the matching rate (i.e., the average number of matches per unit time) for the greedy policy converges to $n\lambda 
$ as $n\rightarrow \infty $.

When a match occurs (i.e., when there is at least one available mate upon an
agent's arrival), the expected number of available mates when an agent
arrives can be approximated by 
\begin{equation}
\sqrt{n}E[N(0,\lambda /\eta )|N(0,\lambda /\eta )>0] = \frac{\lambda }{\eta }\sqrt{\frac{2n}{\pi }}.  \label{eq15}
\end{equation}
By~(\ref{eq14})-(\ref{eq15}) and Lemma \ref{lemma:Expectation_Max}, the utility rate of the greedy policy, which
is denoted by $U_{n}^{g}$, satisfies 
\begin{equation}
U_{n}^{g}\sim n\lambda m\left(\frac{\lambda }{\eta }\sqrt{\frac{2n}{\pi }}\right).  \label{eq16}
\end{equation}%

\section{Simulation Results}

\label{sec-simulation} To assess the accuracy of our asymptotic results, we
consider special cases of the three canonical examples in~\S \ref%
{sec-examples} in the Appendix: exp(1), Pareto(1,2) and U[0,1]. For all cases, we let $%
\lambda =\eta =1$ and $n=1000$, so that the mean number of buyers and
sellers in a match-free system is 1000. We initialize the system with 1000
buyers and 1000 sellers, simulate the system for 1500 time units, discarding
the first 150 time units, and then repeat this procedure 100 times.

To find the optimal population threshold levels, we compute the utility rate for the population threshold policy for each integer
threshold value in the range [0,1000], using the same set of random numbers
for each threshold level. We repeat the same procedure for the utility threshold policy, and discretize the utility threshold values by 0.1 for the exp(1) and Pareto(1,2) cases, and by 0.01 for the U[0,1] case. 

\noindent{\sl Exponential(1) case}. In the exponential case, we predict that the optimal population threshold level is $z_n^*=\frac{1000}{%
\ln 1000}=144.8$, and the utility rate under
this threshold policy approaches the upper bound and is twice as large as the utility rate of the greedy
policy (see~\S \ref{ssec-examples-exp}). The optimal threshold level found
via simulation is 148, and the suboptimality of the utility rate under
the threshold 144.8 vs. the threshold 148 is 0.004\% (Table~\ref{table2}%
). Our heuristic utility threshold is  $v_n^*=5.56$ from~(\ref{eq65}) in the Appendix, which coincides with the optimal threshold found via simulation (with a discretization of 0.1) of 5.6. 

However, the predicted utility rates are less accurate than our determination of the best threshold levels. By~(\ref{eq51a}) in the Appendix, our best estimate for the utility rate under the optimal population threshold policy is 5553, which is 14.9\% higher than the simulated value in Table~\ref{table2}. By~(\ref%
{eq13a}) and~(\ref{eq23a}) in the Appendix, our best estimate of the utility rate under the
greedy policy is 
\begin{eqnarray}
U_n^g & \approx & \frac{\lambda n}{\nu}\Biggl(1-\frac{1}{\sqrt{2\pi n}}%
\Biggr)\Biggl(\gamma+\ln\Biggl(\frac{\lambda}{\eta}\sqrt{\frac{2n}{\pi}}%
\Biggr)\Biggr),  \nonumber \\
& = & 3757,  \nonumber
\end{eqnarray}
which is 8.5\% higher than the simulated value in Table~\ref{table2}. Our
best estimate of the upper bound is given in~(\ref{eq22a}) in the Appendix, which yields
7485. The optimal-to-greedy ratio of the simulated utility rates is $\frac{4833}{3462}=1.40$
rather than 2. Further simulations reveal that convergence is very slow:
this simulated ratio is 1.48 when $n=10^4$ and 1.54 when $n=10^5$. Most of the inaccuracy in estimating the optimal-to-greedy
ratio is due to the fact that the simulated utility rate of the optimal
threshold policy is not very close to the upper bound.

Finally, the utility rate of the optimal utility threshold policy is 5732 (Table~\ref{table3}). While still far from the upper bound, it is 18.6\% higher than the utility rate achieved by the optimal population threshold policy. 

\noindent{\sl Pareto(1,2) case}. In the Pareto case, we predict that the optimal population threshold level is $\frac{1000}{3}=333.3$. The optimal
population threshold found via simulation is 347, and the utility 
suboptimality of the theoretical threshold is 0.03\% (Table~\ref{table2}). The solution to~(\ref{eq58b}) in the Appendix is $z^*=0.512$. Hence, the optimal utility threshold level in~(\ref{eq58c}) in the Appendix is $v_n^*=42.8$, which is very close to the value of 42.0 found via simulation. 

Our estimate of the utility rate under the optimal population threshold policy is 21,573 by~(\ref{eq57}) in the Appendix, which is 2.4\% less than the simulated value of 22,102 in Table~\ref{table2}. The utility rate under the optimal utility threshold policy in~(\ref{eq58d}) in the Appendix is 43,756, which is nearly identical to the optimal simulated value of 43,750. Our best estimate for the utility rate of the greedy
policy is $\Bigl(1-\frac{1}{\sqrt{2\pi n}}\Bigr)$ times the right side of~(\ref{eq55a-mt}) in the Appendix, or 8791, which is 6.4\% larger than the simulated value in
Table~\ref{table2}. Our estimate of the upper bound in~(\ref{eq54a}) in the Appendix is
56,050. By~(\ref{eq58a}) in the Appendix, the predicted performance ratio between the optimal population threshold policy and the greedy policy is $\frac{2}{3}\left(\frac{1000\pi}{18}\right)^{1/4}=2.42$, compared to the optimal-to-greedy simulated ratio of $\frac{22,102}{%
8259}=2.68$ (Table~\ref{table2}). By~(\ref{eq58}) in the Appendix, the ratio of the upper bound to the utility rate of the optimal population threshold policy is predicted to be $\frac{3\sqrt{3}}{2}=2.60$, compared to the simulated value of $\frac{56,050}{22,102}=2.54$. 

The simulated utility rate of the optimal utility threshold policy is nearly twice as large as the simulated utility rate of the optimal population threshold policy (Table~\ref{table3}), although it is still 21.9\% smaller than the predicted upper bound of 56,050. 

\noindent{\sl Uniform(0,1) case}. In the uniform case, we predict that the greedy policy is asymptotically
optimal. The optimal population threshold level found via simulation is 22, and the resulting 
utility suboptimality of the greedy policy is 4.0\% (Table~\ref{table2}). Note that other population thresholds aside from zero are also asymptotically optimal in this case, including $\ln(n)=\ln(1000)=6.91$, which has a suboptimality of 2.0\%. Our best estimate of the utility rate
under the greedy policy is $\Bigl(1-\frac{1}{\sqrt{2\pi n}}\Bigr)$ times the
right side of~(\ref{eq34a}) in the Appendix, or 948.2, which is 4.4\% larger than the
simulated value in Table~\ref{table2}. The upper bound in~(\ref{eq33a}) in the Appendix
equals 987.4, which is 4.3\% larger than the utility rate corresponding to the optimal population threshold level of 22. The predicted optimal utility threshold from equation~(\ref{eq61b}) in the Appendix is $v_n^*=0.974$, compared to the the value of $0.96$ found via simulation, for a utility suboptimality of 0.24\% (Table~\ref{table3}). 

In summary, our analysis identifies the optimal threshold level
within about 2\% (considering the possible range of [0,1000]) and its
suboptimality is no more than 2\% for the population threshold policy in the uniform case, and is negligible in the other five cases. We also note that the predicted fraction of agents
who abandon the market under the optimal population threshold policy,
which is $\frac{2z_n}{2\lambda_n}=\frac{1}{\ln n}=0.145$ (i.e., the total
abandonment rate divided by the total arrival rate) in the exponential case, 
$\frac{z_n^*}{n}=\frac{1}{3}$ in the Pareto case by~(\ref%
{eq56}) in the Appendix, and $1-\frac{1}{\sqrt{2\pi n}}=0.013$ in the uniform case by~(\ref{eq13a}), are
reasonably close to the simulated values in the fourth column of Table~\ref%
{table3}. As predicted by our analysis, the utility rate of the greedy
policy -- normalized by the mean of the matching distribution -- increases
with the right tail of the matching distribution (this quantity is 1817 for the
uniform, 3462 for the exponential, and 4129 for the Pareto), as does the
ratio of the utility rates between the best threshold policy and the greedy
policy (1.04 for the uniform, 1.40 for the exponential, and 2.68 for the
Pareto under the population threshold policy, and 1.06, 1.66 and 5.30 under the utility threshold policy). In addition, despite the asymptotic optimality result, there is a
large gap between the utility rate of the best population threshold policy and the
upper bound in the exponential case. The improvement of the utility threshold policy over the population threshold policy also increases with the right tail of the matching distribution, with the ratio of the utility rates equaling 1.02, 1.19 and 1.98 for the uniform, exponential and Pareto cases, respectively. This improvement is achieved by being more patient and allowing more agents to abandon the market, particularly in the Pareto case (last column in Table~\ref{table3}).

\section{Unbalanced Markets}

\label{sec-unbalanced} 

In this section, we consider unbalanced markets, where buyers and sellers arrive at rates $n\lambda_b$ and $n\lambda_s$ in the $n^{\rm th}$ system, and abandon at rates $\eta_b$ and $\eta_s$, respectively. We restrict ourselves to the analysis of the utility threshold policy in the case $\alpha\in (0,1)$, which is very similar to the corresponding analysis in the symmetric case. We also note that an analysis of the population threshold policy in the unbalanced case is complicated by the extra degree of freedom that is introduced (and needs to be determined) in equations~(\ref{eq69a})-(\ref{eq:Fluid_Stationary}) in the Appendix, and is beyond the scope of this paper. Under the utility threshold policy in the $n^{\rm th}$ system, an arriving buyer is matched to the seller that yields the maximum utility if this utility exceeds the threshold $v_{n,s}$; similarly, an arriving seller is matched to its highest-matching buyer if the utility exceeds the threshold $v_{n,b}$. The dynamics are given by the equations
\begin{align}
B_{n}\left( t\right) & =B_{n}\left( 0\right) +\sum_{j=1}^{N_{B}^{+}\left(
n\lambda _{b}t\right) }I_{\left\{\max_{i=1}^{S_{n}\left( A_{j-}^B\right)
}V_{i,j}^B\leq v_{n,s}\right\}} -\sum_{j=1}^{N_{S}^{+}\left( n\lambda
_{s}t\right) }I_{\left\{\max_{i=1}^{B_{n}(A_{j-}^S)}V_{i,j}^{S
}>v_{n,b}\right\}}   \label{eq40} \\
& -N_{B}^{-}\left( \eta _{b}\int_{0}^{t}B_{n}\left( r_{-}\right) dr\right) ,
\nonumber \\
S_{n}\left( t\right) & =S_{n}\left( 0\right) +\sum_{j=1}^{N_{S}^{+}\left(
n\lambda _{s}t\right) }I_{\left\{ \max_{i=1}^{B_{n}(A_{j-}^S)}V_{i,j}^{S}\leq v_{n,b}\right\}} -\sum_{j=1}^{N_{B}^{+}\left(
n\lambda _{b}t\right) }I_{\left\{ \max_{i=1}^{S_{n}(A_{j-}^B)}V_{i,j}^B>v_{n,s}\right\}}   \label{eq41} \\
& -N_{S}^{-}\left( \eta _{s}\int_{0}^{t}S_{n}\left( r\right) dr\right) ,
\nonumber
\end{align}%
where $\left\{ A_{j}^B:j\geq 1\right\} $ is the sequence of arrival times associated with $%
N_{B}^{+}(n\lambda _{b}\cdot )$, and $\left\{ A_{j}^S:j\geq 1\right\}
$ is the sequence of arrival times associated with $N_{S}^{+}(n\lambda _{s}\cdot )$. As
in the symmetric case, the $V_{i,j}^B$s and $V_{i,j}^{S}$s are
independent arrays of i.i.d. random variables with distribution $F\left(
\cdot \right) $.

Following the development in the symmetric case (e.g., (\ref{eq6b}%
)), the utility rate takes the form
\begin{eqnarray}
U_{n}^u(v_{n,b},v_{n,s}) &=&n\lambda _{b}E\left[E[M\left( S_{n}\left( \infty \right) \right) I_{\{
M\left( S_{n}\left( \infty \right) \right) \geq v_{n,s}\}} |S_{n}\left(
\infty \right) ]\right] \nonumber \\
&&+n\lambda _{s}E\left[E[M\left( B_{n}\left( \infty \right) \right)
I_{\{M\left( B_{n}\left( \infty \right) \right) \geq v_{n,b}\}}
|B_{n}\left( \infty \right) ]\right] . \label{unbalanced-utility}
\end{eqnarray}%

Our main result is presented in Theorem~\ref{them}. The proof of Theorem~\ref{them} appears in~\S\ref{ssec-proofs-unbalanced} in the Appendix and largely mimics the proof of Theorem~\ref{thm:utility_based}.  

\begin{theorem}\label{them}
	Suppose that Assumption \ref{as:Max_Pareto} holds. Let $v_*$ be the optimal solution for the optimization problem
	\begin{align}
		\max_{v\geq 0} & \lambda_s b^\alpha E[XI_{\{b^\alpha X>v\}}]+\lambda_b s^\alpha E[XI_{\{s^\alpha X>v\}}] \label{opt:theorem}\\
	\text{s.t. }&\eta_b b+\lambda_s=\lambda_b \exp(-\kappa s/v^{1/\alpha})+\lambda_s\exp(-\kappa b/v^{1/\alpha})\label{opt:condition_1}, \\
	&\eta_b b+\lambda_s = \eta_s s +\lambda_b\label{opt:condition_2}.
\end{align}
Then a threshold policy of the form $v_{n,b}^*=v_{n,s}^*=m(n)v_*$ is asymptotically optimal among the class of utility threshold polices. The associated utility rate satisfies
\begin{align*}
	\lim_{n\rightarrow\infty}\frac{U_n^u(v_{n,b}^*,v_{n,s}^*)}{nm(n)}=\lambda_s {b_*}^\alpha E[XI_{\{b_*^\alpha X>v_*\}}]+\lambda_b s_*^\alpha E[XI_{\{s_*^\alpha X>v_*\}}],
\end{align*}
where $b_*,s_*$ are solutions that satisfy constraints \eqref{opt:condition_1}-\eqref{opt:condition_2}.
\end{theorem}

Although the results in Theorem~\ref{them} are beyond our intuitive grasp, we attempt to provide some possible intuition for why $v_b^*=v_s^*$ in the unbalanced case. 
Let us consider a fluid model in which the two utility thresholds are both equal to $v^*$. We can classify the matched sellers into two categories: actively matched (i.e., they arrive to the market and are immediately matched with buyers) and passively matched (i.e., they wait in the market and then are matched with arriving buyers). Now suppose that we change the utility thresholds to $v_b$ and $v_s$, where $v_b<v^*<v_s$, in such a way that the number of additional actively-matched sellers (call it $ds$) equals the reduction in the number of passively-matched sellers. Because the total number of matched sellers does not change, the number of abandoned agents remains the same and we can focus on the matching utilities of these marginal sellers. Let $u_a$ be the utility per match for the $ds$ additional actively-matched sellers and $u_p$ be the utility rate per match for the $ds$ sellers that are no longer passively matched. The utility per match of these marginal sellers is between the new threshold and the old threshold, and therefore $v_b<u_a<v*$ and $v*<u_p<v_s$. Hence, the net change in utility is $(u_a-u_p)ds$, which is negative.

We conclude this section with a numerical example that is a variant of the one in~\S\ref{ssec-examples-pareto}: let $\lambda_b=2$, $\lambda_s=1$, $\eta_b=1$, $\eta_s=1$, $n=1000$, and assume a Pareto(1,2) distribution, so that $\alpha=1/2$, $\kappa=1/\pi$ and $m(n)=\sqrt{1000\pi}$. Then $b(s)=s+1$ in~(\ref{eq51}) in the Appendix, and~(\ref{eq55})  in the Appendix reduces to 
\[
2e^{-s\tau}+e^{-(s+1)\tau}=s+2.
\]
The solution to~(\ref{eq53}) in the Appendix is $s_*=0.365$ and $\tau(0.365)=0.361$, which yields $v_{n,b}^*=v_{n,s}^*=\sqrt{\frac{1000}{0.361}} = 52.7$. Interestingly, this threshold level of 52.7 is higher than in the symmetric case, where $\lambda_b=1$ and $v_n^*=42.8$. Moreover, leaving all parameter values fixed except for $\lambda_b$, we numerically compute $v_{n,b}^*$ in~(\ref{eq54}) in the Appendix and find that it is increasing and concave in $\lambda_b\ge 1$. 

With $\lambda_b=2$, we simulate this system in the same manner as in~\S\ref{sec-simulation}. At a discretization of 0.1, a two-dimensional search of $(v_{n,b},v_{n,s})$ space via simulation for the optimal thresholds yields (52.7,52.3), with a corresponding simulated utility rate of 71,046 and with abandonment fractions of 0.681 for buyers and 0.363 for sellers. The simulated utility rate at $(v_{n,b}^*,v_{n,s}^*)=(52.7,52.7)$ is 71,010, which is suboptimal by 0.05\%. The predicted utility rate, $U_{n}^u(v_{n,b}^*,v_{n,s}^*)$ in~(\ref{eq55a}) in the Appendix, is 70,992, which is 0.03\% less than the simulated value of 71,010. 

Fixing one utility threshold level at 52.7 and varying the other threshold level (Fig. 1 in the Appendix) reveals that the simulated utility rate is slightly more sensitive to $v_s$ than $v_b$, perhaps because arriving buyers see fewer potential matches than arriving sellers. This figure also shows that it is more suboptimal to underestimate the threshold level than to overestimate it. 

\section{Batch-and-Match Policy}
\label{sec-batch} 

In this section, we restrict ourselves to Pareto matching utilities with finite mean, where $F(v)=1-(cv)^{-\beta}$, for $\beta>1$, $c>0$ and $cv\ge 1$, so that $\alpha=1/\beta$ in Assumption~\ref{as:Regularly_Varying}. We consider a one-parameter batch-and-match policy: At times $t = \{\Delta, 2\Delta, 3\Delta,\ldots,\}$, we match $\min\{B(t),S(t)\}$ buyers and sellers by randomly choosing $\min\{B(t),S(t)\}$ agents from the thicker side of the market (e.g., buyers if $B(t)\ge S(t)$) and then maximize the total utility from these matches; this class of policies allows us to consider a balanced random assignment problem, which is easier to analyze than an unbalanced random assignment policy. The goal is to choose the time window $\Delta$ that maximizes the long-run average utility rate. The main qualitative conclusion from this section is that -- for the special case of $\lambda=\eta=c=1$ and $\beta=2$ -- the utility threshold policy easily  outperforms this batch-and-match policy, both in the asymptotic analysis and in the simulation results. 

To analyze the performance of this policy, we consider a random assignment problem, where there are $k$ buyers and $k$ sellers with iid Pareto matching utilities $V_{i,j}$ between buyer $i$ and seller $j$. The matching problem is
\begin{eqnarray*}
	\max_\pi \sum_{i=1}^k V_{i,\pi(i)},
\end{eqnarray*}
where $\pi$ is a permutation function. Let ${\cal M}(k) = \max_\pi \sum_{i=1}^k V_{i,\pi(i)}.$


The main result of this section is given in Theorem~\ref{th:batch}, which is proved in~\S\ref{ssec-proofs-batching} in the Appendix. The corresponding results for the unbalanced case are presented without proof at the end of~\S\ref{ssec-proofs-batching}. 

\begin{theorem} \label{th:batch}
Consider the symmetric model with arrival rate $\lambda$, abandonment rate $\eta$, and Pareto($c,\beta)$ matching utilities with finite mean. Let $U_n^b(\Delta)$ be the utility rate of the batch-and-match policy with time window $\Delta$. 
Then the utility rate satisfies 
\begin{equation}
    \label{batch5}
\lim_{n\to\infty}\frac{U_n^b(\Delta)}{n^{\alpha+1}} \le  \frac{c\Gamma(1-\alpha)\Bigl[\frac{\lambda}{\eta}(1-e^{-\eta \Delta^*})\Bigr]^{\alpha+1}}{\Delta^*},
\end{equation}
where the asymptotically optimal time window $\Delta^*$ is the unique solution to 
\begin{equation}
    \label{batch2}
    e^{\eta \Delta}=(1+\alpha)\eta \Delta+1.
\end{equation}
\end{theorem}
Note that $\Delta^*$ increases in $\alpha$ in~(\ref{batch2}), and so -- as in the population threshold policy -- heavier tails lead to thicker markets. 

We conclude this section with the numerical Pareto(1,2) example from~\S\ref{sec-simulation}, where $\lambda=\eta=1$, $n=1000$, $F(v)=1-v^{-2}$ and $\alpha=1/2$. Equation~(\ref{batch2}) reduces to $e^{\Delta}=\frac{3}{2}\Delta+1$, which has solution $\Delta^*=0.76$, implying from~(\ref{batch6}) in the Appendix that we make approximately $1000(1-e^{-0.76})=532$ matches in each cycle. The upper bound for the utility rate in~(\ref{batch5}) is $\frac{\sqrt{\pi}(1-e^{-\Delta^*})^{3/2}1000^{3/2}}{\Delta^*}=28,644$, which is much smaller than the predicted utility rate of 43,750 for the utility threshold policy from~(\ref{eq58d}). A one-dimensional search using simulation generates an optimal time window of 0.75, confirming the accuracy of our asymptotic analysis. The simulated utility rate under this time window is 25,168 and the lower bound for the utility rate in Lemma~\ref{lemma:batch} is 131, suggesting that the upper bound is useful and the lower bound is very loose.

\section{Concluding Remarks}
\label{sec-conclusion} A fundamental tradeoff in centralized dynamic matching markets relates to market thickness: whether matches should be delayed -- at the risk of antagonizing waiting agents -- in the hope of obtaining better matches in the future. Very little is known about this issue when matching utilities are general. By combining queueing asymptotics (as an aside, we note that perhaps the most surprising part of our study is that rather than requiring a diffusion analysis, a fluid analysis is sufficient to analyze this problem) with extreme value theory, we obtain explicit results that shed light on this issue. For symmetric markets, as the right tail of the matching utility distribution gets heavier, it is optimal to become more patient and let the market thickness (and abandonment rate) increase. While empirical work on matching markets use more complicated covariate models than what we consider (e.g., Hitsch \textsl{et al.} 2010, Boyd \textsl{et al.} 2013, Agarwal 2015), it seems clear from these analyses that matching
utilities typically are not in the domain of attraction of the Weibull law. Therefore, large centralized
matching markets -- whether balanced or unbalanced (see below) -- are likely to benefit from allowing
the market to thicken. 

Enabled by the decoupling of the fluid queueing dynamics and the extremal behavior of the matching utilities, our study appears to be the first to allow for correlated matching utilities, which is likely to be a common phenomenon in practice: an agent who is deemed objectively attractive in a labor, housing or school choice model is likely to have matching utilities with potential mates that are positively correlated rather than i.i.d. In~\S\ref{ssec-examples-correlated}, we find that positive correlation reduces the market thickness in the utility threshold policy but not the population threshold policy, and reduces the utility rate under both policies. 

We note four limitations in our study. First, most of our analysis is restricted to arrival-only policies. In particular, it might be possible to do better by batching sets of agents and then matching them, as in Mertikopoulos \textsl{et al.} (2020). Moreover, generalizing their results to our setting is likely to be quite challenging, in that the $\pi^2/6$ result requires an exponential matching distribution and an objective of minimizing the matching cost (they minimize mismatch plus waiting costs rather than maximizing utility in the presence of abandonment). While they generalize their results in~\S6 of their paper by positing a functional form for how the expected minimum mismatch costs decrease as a function of the number of agents in the market, this functional form does not appear to follow from any more primitive distributional assumptions. In~\S\ref{sec-batch} we consider Pareto utilities and analyze a simple batch-and-match policy, which periodically (with an asymptotically optimal time window) optimally matches all agents on the thinner side of the market with an equal number of agents randomly selected from the thicker side of the market. Perhaps surprisingly, we show that in the Pareto case, the utility threshold policy easily outperforms the batch-and-match policy. Nonetheless, this does not preclude the possibility that more sophisticated batching policies (e.g., optimally -- rather than randomly -- select the agents to match from the thicker side of the market, or include a utility threshold for allowable matches as a second parameter) might outperform the utility threshold policy. 

Second, most of our analysis considers a symmetric market, with buyers and sellers having the same arrival and abandonment rates. While some markets, such as cadaveric organ transplants and public housing, tend to have chronic supply shortages, other
markets have economic forces at play that tend to roughly balance supply and demand.
In a static matching market, even a slight imbalance can 
give rise to a unique stable matching (Ashlagi \textsl{et al.} 2017b). We also note that a greedy policy is optimal in a somewhat different unbalanced market setting, where easy-to-match agents can match with all other agents in the market with a specified probability, but hard-to-match agents can match only with easy-to-match agents with a different specified probability (Ashlagi \textsl{et al.} 2019b). In our analysis of the utility threshold policy in the heavy-tailed case of an unbalanced market, we obtain the somewhat surprising result that the solution is symmetric: i.e., the utility threshold is the same for buyers and sellers. Moreover, we find (in our Pareto example) that the amount of patience increases with the amount of imbalance; i.e., the larger the imbalance, the more agents that are going to be turned away, and the more selective the matching becomes. However, we leave a complete analysis of the unbalanced problem for future work. 

Although our model can be viewed as allowing a continuum of classes via the distribution of the matching utility, the third restriction is that our analysis does not naturally lend itself to a setting where there are a discrete number of classes with class-dependent matching utilities. In particular, in some settings (e.g., organ donation) some classes of buyers/sellers are compatible with only  pre-specified classes of sellers/buyers. 
The decoupling of the queueing fluid dynamics and the extremal behavior of the matching utilities should carry over to the setting with a finite number of classes with some incompatibility among classes. However, it would make sense to consider multiple thresholds in this setting, and a multi-dimensional model with multiple thresholds would be a nontrivial extension. 

The final restriction is exponential abandonment. Relaxing this assumption would require a different approach, such as hazard rate scaling (Reed and Tezcan 2012), and would likely be much more difficult.

Finally, we note that there may be equity issues if a significant number of agents are allowed to abandon the market (Table~\ref{table3}). The consideration of a risk-sensitive objective function would likely require a diffusion approximation, which would be, e.g., a two-dimensional Ornstein-Uhlenbeck process with an unusual Skorokhod condition under a population threshold policy.  

\vskip0.2truecm \centerline{\bf Acknowledgment} \vskip0.2truecm We thank
Can Wang and Halwest Mohammad for running some of the simulations described in~\S \ref{sec-simulation}-\ref{sec-unbalanced} and Kavita Ramanan for advice about the Skorokhod mapping in~\S\ref{sec-proofs}. J. Blanchet also acknowledges NSF support through the grants 1915967, 1820942, 1838576.

\clearpage

\begin{table}
	\centering
	\begin{adjustwidth}{-0.5cm}{}
		\begin{tabular}{|c||c|c|c|}
			\hline
\multirow{2}{*}{}  & \multicolumn{3}{|c|}{Matching Utility Distribution} \\ \cline{2-4}
Policy & Exponential($\nu$) & Pareto ($c$, shape $\beta>1$) & Uniform($a,b)$ \\ \hline\hline
Upper Bound & $U_n^+\sim \frac{\lambda}{\nu} n\ln n$ & $U_n^+ = O(n^{1+1/\beta})$ & $U_n^+\sim \lambda bn$ \\ \hline
Greedy Policy & $U_n^g\sim \frac{\lambda}{2\nu} n\ln n$ & $U_n^g = O(n^{1+1/(2\beta)})$& asymptotically optimal \\ \hline
Population  & $z_n^*=\frac{n}{\ln n}$ is  & $z_n^*=\frac{\lambda}{\eta(1+\beta)}n$ & $z_n^*=0$ is  \\ 
Threshold   & asymptotically  & $U_n^p(z_n^*) = O(n^{1+1/\beta})$ but no & asymptotically optimal \\ 
Policy with &  optimal & convergence to upper  & \\
threshold $z_n$ & &  bound unless $\beta\to\infty$ & \\ \hline
Utility & heuristic $v_n^*=$  & $v_n^* = 1.353\sqrt{n}$ & heuristic $v_n^*=a+$ \\ 
Threshold   & $\frac{\ln n -\ln \ln n-\ln\ln\left(\frac{2\lambda\ln n}{\lambda \ln n + \eta}\right)}{\nu}$ & when $c=1$, $\beta=2$; & $(b-a)\left(\frac{1}{\sqrt{2n\pi}}+\frac{1}{2}\right)^{\frac{\eta}{\lambda}\sqrt{\frac{\pi}{2n}}}$\\ 
Policy with  & & & \\ 
threshold $v_n$ & & $U_n^u(v_n^*) = O(n^{1+1/\beta})$ & \\ \hline
\end{tabular}
\caption{Summary of results for the three canonical cases in~\S\ref{sec-examples}. $U_n^+$, $U_n^g$, $U_n^p(z_n)$ and $U_n^u(v_n)$ are the upper bound on the utility rate for any arrival-only policy, the utility rate for the greedy policy, the utility rate for the population threshold policy with threshold $z_n$, and the utility rate for the utility threshold policy with threshold $v_n$, all for the $n^{\rm th}$ system. 
} 
\label{table1}
	\end{adjustwidth}
\end{table}

\begin{table}[tbp]
\centering
\begin{adjustwidth}{-0.5cm}{}
		\begin{tabular}{|c||c|c||c|c|c|}
			\hline
\multirow{3}{*}{} & \multicolumn{2}{|c||}{Optimal Population Threshold} & \multicolumn{3}{|c|}{Simulated Utility Rate [95\% CI]} \\ \cline{2-6}
Utility & &  & Theoretical & Simulation & Greedy \\
Distribution & Theoretical & Simulation & Threshold & Threshold & Policy\\ \hline\hline
Exponential(1) & 144.8 & 148 & 4833 & 4833 & 3462 \\ 
& & & [4824,4840] & [4827,4841] & [3425,3503] \\ \hline 
Pareto(1,2) & 333.3 & 347 & 22,095 & 22,102 & 8259 \\
& & & [21,997,22,241] & [21,972,22,234] & [8107,8428] \\ \hline
Uniform(0,1) & 0 & 22 & 908.4 & 946.3 & 908.4 \\
& & & [906.0,911.3] & [945.1,947.7] & [906.0,911.3] \\ \hline
\end{tabular}
\caption{Theoretical and simulation results for the population threshold policy.  } 
\label{table2}
	\end{adjustwidth}
\end{table}

\begin{table}[tbp]
\centering
\begin{adjustwidth}{-1.6cm}{}
		\begin{tabular}{|c||c|c|c||c|c|c|}
			\hline
\multirow{3}{*}{} & \multicolumn{3}{|c||}{Population Threshold Policy} & \multicolumn{3}{|c|}{Utility Threshold Policy} \\ \cline{2-7}
Utility & Optimal & Utility & Fraction & Optimal  & Utility & Fraction \\
Distribution & Threshold & Rate & Abandoned & Threshold & Rate & Abandoned \\ \hline\hline
Exponential(1) & 148 & 4833 & 0.140 & 5.6 & 5732 & 0.150 \\ 
& & [4827,4841] & & & [5724,5740] & \\ \hline 
Pareto(1,2) & 347 & 22,102 & 0.334 & 42.0 & 43,750 & 0.503 \\
& & [21,972,22,234] & & & [43,541,43,960] & \\ \hline
Uniform(0,1) & 22 & 946.3 & 0.027 & 0.96 & 963.0 & 0.021 \\
& & [945.1,947.7] & & & [961.7,964.2] & \\ \hline
\end{tabular}
\caption{Simulation results for both threshold policies. Columns~2 and~3 are taken from Table~\ref{table2}.  } 
\label{table3}
	\end{adjustwidth}
\end{table}
\clearpage

\centerline{\Large{\bf Appendix}}

The proofs appear in~\S\ref{sec-proofs}, the main results are applied to canonical examples in~\S\ref{sec-examples}, and some useful facts about extreme value theory and regularly varying functions are collected in~\S\ref{sec-evt}.

\section{Technical Proofs}
\label{sec-proofs}

Before we prove Theorems~\ref{thm:Population_Based},~\ref{thm:utility_based},~\ref{them} and~\ref{th:batch} in~\S\ref{ssec-proofs-population}, ~\S\ref{ssec-proofs-utility}, ~\S\ref{ssec-proofs-unbalanced} and~\S\ref{ssec-proofs-batching} respectively, we prove Lemmas~\ref{lemma:Expectation_Max} and~\ref{lemma:bound}, and state and prove Lemmas~\ref{lemma:VerySlowlyVarying} and~\ref{lemma:Proc_Limit}, in~\S\ref{ssec-proofs-lemmas}. Lemma~\ref{lemma:Proc_Limit} is used in the proof of Theorem~\ref{thm:utility_based}. 

\subsection{Preliminary Lemmas}
\label{ssec-proofs-lemmas}

\noindent{\textbf{Proof of Lemma \ref{lemma:Expectation_Max}}}

By uniform local convergence, we must have that for any $0<a<b<\infty$ (see Resnick 1987, Section 0.4, Proposition 0.5),
\[
\lim_{x\rightarrow\infty}\sup_{a\leq t\leq b}\left\vert \frac{m\left(
xt\right)  }{m\left(  x\right)  }-t^{-\alpha}\right\vert =0.
\]
Therefore,
\begin{eqnarray}
\overline{\lim}_{n\rightarrow\infty}E\left[
\frac{\left\vert m\left(  N_{n}\right)  -m\left(  \bar{N}_{n}\right)
\right\vert I_{\{ \left\vert N_{n}-\bar{N}_{n}\right\vert \leq
\varepsilon\bar{N}_{n}\}}}{m\left(  \bar{N}_{n}\right)}  \right] 
&\leq& \overline{\lim}_{n\rightarrow\infty}\sup_{1-\varepsilon\leq
t\leq1+\varepsilon}\left\vert \frac{m\left(  \bar{N}_{n}t\right)  }{m\left(
\bar{N}_{n}\right)  }-1\right\vert, \nonumber \\ 
&\leq& O\left(  \varepsilon\right)  . \nonumber
\end{eqnarray}
On the other hand, note that
\begin{align*}
& E\left[\left\vert \frac{m\left(  N_{n}\right)  }{m\left(  \bar{N}%
_{n}\right)  }-1\right\vert I_{\{\left\vert N_{n}-\bar{N}_{n}\right\vert
>\varepsilon\bar{N}_{n}\}}  \right]  \\
& \leq E\left[\left\vert \frac{m\left(  N_{n}\right)  }{m\left(\bar{N}%
_{n}\right)  }-1\right\vert I_{\{N_{n}-\bar{N}_{n}>\varepsilon\bar{N}%
_{n}\}}  \right]  +E\left[  \left\vert \frac{m\left(  N_{n}\right)
}{m\left(  \bar{N}_{n}\right)  }-1\right\vert I_{\{\bar{N}_{n}%
-N_{n}>\varepsilon\bar{N}_{n}\}}  \right]  .
\end{align*}
On the event $\bar{N}_{n}\left(  1-\varepsilon\right)  >N_{n}$, since
$m\left(  \cdot\right)  $ is nondecreasing, we have that
\[
\left\vert \frac{m\left(  N_{n}\right)  }{m\left(  \bar{N}_{n}\right)
}-1\right\vert =1-\frac{m\left(  N_{n}\right)  }{m\left(  \bar{N}_{n}\right)
}\leq1 ,
\]
and therefore
\[
E\left[\left\vert \frac{m\left(  N_{n}\right)  }{m\left(  \bar{N}%
_{n}\right)  }-1\right\vert I_{\{\bar{N}_{n}-N_{n}>\varepsilon\bar{N}%
_{n}\}}  \right]  \leq P\left(  \left\vert N_{n}-\bar{N}_{n}\right\vert
>\varepsilon\bar{N}_{n}\right)  =o\left(  1\right) ~~{\rm as} ~~ n\rightarrow\infty . 
\]
On the other hand, again, because $m\left(
\cdot\right)  $ is nondecreasing,
\[
E\left[\left\vert \frac{m\left(  N_{n}\right)  }{m\left(  \bar{N}%
_{n}\right)  }-1\right\vert I_{\{N_{n}-\bar{N}_{n}>\varepsilon\bar{N}%
_{n}\}}  \right]  \leq E\left[  \frac{m\left(  N_{n}\right)  }{m\left(
\bar{N}_{n}\right)  }I_{\{N_{n}>\bar{N}_{n}\left(  1+\varepsilon\right)
\}}  \right]  .
\]
Applying Potter's bound (Bingham, Goldie and Teugels (1987), Theorem 1.5.6 part (iii)) for each
$\delta>0$, there exists $t>0$ such that $m\left(  y\right)  /m\left(
x\right)  \leq2\left(  y/x\right)  ^{\alpha+\delta}$ if $y\geq x\geq t$.
Hence, for any $\delta>0$ there exists $n_{0}>0$ such that $\bar{N}_{n}\geq t$
for all $n\geq n_{0}$, and therefore
\begin{eqnarray*}
E\left[\frac{m\left(  N_{n}\right)  }{m\left(  \bar{N}_{n}\right)
}I_{\{N_{n}>\bar{N}_{n}\left(  1+\varepsilon\right)  \}}  \right]
&\leq& 2E\left[\left(  \frac{N_{n}}{\bar{N}_{n}}\right)  ^{\alpha+\delta
}I_{\{N_{n}>\bar{N}_{n}\left(  1+\varepsilon\right)  \}}  \right] ,  \\
& \leq& 2\left(  E\left[ \left(  \frac{N_{n}}{\bar{N}_{n}}\right)  ^{\left(
\alpha+\delta\right)  r}\right]  \right)  ^{1/r}P\left(  \left\vert N_{n}%
-\bar{N}_{n}\right\vert >\varepsilon\bar{N}_{n}\right)  ^{1/s},
\end{eqnarray*}
for any $r,s>1$ such that $1/r+1/s=1$, by H\"{o}lder's inequality. Furthermore, because  $\alpha<1$, we
can guarantee $\left(  \alpha+\delta\right)  r<1$ by choosing $\delta>0$
sufficiently small. Next, Jensen's inequality implies that 
\[
E\left[\left(  \frac{N_{n}}{\bar{N}_{n}}\right)  ^{\left(  \alpha
+\delta\right)  r}\right]  \leq\left(  E\left[  \frac{N_{n}}{\bar{N}_{n}%
}\right]  \right)  ^{\left(  \alpha+\delta\right)  r} = 1 .
\]
Therefore, we obtain that
\[
E\left[\frac{m\left(  N_{n}\right)  }{m\left(  \bar{N}_{n}\right)  }I_{\{ N_{n}-\bar{N}_{n} >\varepsilon\bar{N}_{n}\}}  \right]  \rightarrow0  ~~{\rm as}  ~~ n\rightarrow\infty,
\]
which completes the proof.

\noindent{\textbf{Proof of Lemma~\ref{lemma:bound}}}

For any arrivals-only policy, $B_{n}\left( \infty \right) $ is stochastically bounded
by a system where no matches occur, i.e., agents leave only upon abandonment. In this case, each side of the market can be modeled as a $M/M/\infty$ queue, which has a $\text{Poisson}\left( \lambda n/\eta
\right) $ stationary queue length distribution. Because the total arrival rate of agents is $2\lambda n$ and two agents exit upon each match, the maximum long-run rate for matches under any policy is $\lambda n$. Hence, for any arrival-only policy we have that $U_{n}\leq \lambda nE\left[ m\left( P_n\right) \right] $, where $P_n$ is a $\text{Poisson}\left( \lambda n/\eta
\right) $ random variable.
Because Lemma \ref{lemma:Expectation_Max} applies for $P_n$, it follows that $U_{n}\leq \lambda nm\left( \lambda n/\eta
  \right) $. 

\begin{lemma}\label{lemma:VerySlowlyVarying} 
For $\alpha = 0$, there exists an $o(n)$ sequence $z_n$ such that $\lim_{n\to \infty} \frac{m(z_n)}{m(n)} = 1$.
\end{lemma}

\noindent{\textbf{Proof of Lemma \ref{lemma:VerySlowlyVarying}}}

Recall that $m(t) =0$ for $t\in[0,1)$ and let us define $$x_n =  \inf\left(x\in (0,1]:  1- \frac{m(x n)}{m(n)}    \le x \right).$$   We show that $x_n$ is $o(1)$. Suppose that this is not true. Then, we have that $\lim\sup_n x_n > 0$. Let $\epsilon = \lim\sup_n x_n$. By Assumption~\ref{as:Regularly_Varying}, for all $a \in (0,1]$ and for all $\epsilon'>0$ there exists $n_0(a, \epsilon')$ such that for all $n \ge n_0$ we have $1 - \frac{m(an)}{m(n)} \le \epsilon'$. Now, we pick $a=\epsilon/2$ and $\epsilon' = \epsilon/2$. Then, for each $n\ge n_0(\epsilon/2,\epsilon/2)$, we have that  $1 - \frac{m(\epsilon n/2)}{m(n)} \le \epsilon/2$. Hence, $x_n \le \epsilon/2$ for each $n \ge  n_0(\epsilon/2,\epsilon/2)$. This implies $ \lim\sup_n x_n \le \epsilon/2$ which is a contradiction. 

Let $z_n = n x_n$. Because $x_n$ is $o(1)$, it follows that $z_n$ is $o(n)$. Further, by construction, $\lim_{n\to \infty} \frac{m(z_n)}{m(n)} =1$. This completes the proof of the lemma.

\begin{lemma}\label{lemma:Proc_Limit}
Suppose that $\left\{ X_{n}\left( \cdot \right) \right\}
_{n\geq 1}$ is a sequence of stochastic processes in $\mathbb{R}^{d}$ and
define $X_{n}^{\ast }\left( t\right) =\sup_{0\leq s\leq t}\left\vert
X_{n}\left( s\right) \right\vert $. Assume that for each $t>0$, the sequence 
$\left\{ X_{n}^{\ast }\left( t\right) \right\} _{n\geq 1}$ is tight.
Moreover, suppose that 
\[
X_{n}\left( t\right) =X_{n}\left( 0\right) +\int_{0}^{t}b_{n}\left(
X_{n}\left( s\right) \right) ds+M_{n}\left( t\right),
\]%
where $\left\{ b_{n}\left( \cdot \right) \right\} _{n\geq 0}$ is a
sequence of continuous functions such that $b_{n}\left( \cdot \right)
\rightarrow b\left( \cdot \right) $ uniformly on compact sets and $%
b\left( \cdot \right) $ is locally Lipschitz, and $M_n(\cdot)$ is a martingale with quadratic variation $\left[
M_{n}\right] \left( \cdot \right) $ satisfying 
\[
E\left[ \left[ M_{n}\right] \left( t\right) \right] \rightarrow 0 
\]%
as $n\rightarrow \infty $ for each $t>0$. Finally, suppose that $X_{n}\left(
0\right) \rightarrow X\left( 0\right) $ in probability as $n\rightarrow
\infty $. Then $X_{n}\left( \cdot \right) \rightarrow X\left( \cdot \right) 
$ in probability in the uniform topology on compact sets, where $X\left(
\cdot \right) $ is the unique solution to 
\[
X\left( t\right) =X\left( 0\right) +\int_{0}^{t}b\left( X\left( s\right)
\right) ds. 
\]
\end{lemma}

\noindent{\textbf{Proof of Lemma \ref{lemma:Proc_Limit}}}

We have that for any $c>0$, on the set $X_{n}^{\ast}\left(  1\right)  \leq
c<\infty$, the process
\[
W_{n}\left(  t\right)  =\int_{0}^{t}b_{n}\left(  X_{n}\left(  s\right)
\right)  ds
\]
satisfies the following: there exists $n_{0}:=n_{0}\left(  c\right)  <\infty$ such that for
each $n\geq n_{0}$ and for any $0\leq s<t\leq1$, we have 
\begin{align*}
\left\vert W_{n}\left( t\right) -W_{n}\left( s\right) \right\vert & \leq
\int_{s}^{t}\left\vert b_{n}\left( X_{n}\left( r\right) \right) \right\vert
dr, \\
& \leq \sup_{0\leq r\leq 1}\left\vert b_{n}\left( X_{n}\left( r\right)
\right) \right\vert \left\vert t-s\right\vert, \\ & \leq (1+\sup_{0\leq r\leq
1}\left\vert b\left( X_{n}\left( r\right) \right) \right\vert )(t-s).
\end{align*}

The last inequality follows because $b_{n}\left(  \cdot\right)  \rightarrow
b\left(  \cdot\right)  $ uniformly on compact sets and because of the tightness of
$X_{n}^{\ast}\left(  1\right)  $.
The Arzela-Ascoli
theorem implies that $W_{n}\left( \cdot \right) $ is tight in the uniform topology.
The Burkholder-Davis-Gundy inequality implies that the martingale sequence converges to zero uniformly on
compact sets and therefore is tight. We
conclude that $X_{n}\left( \cdot \right) $ must be tight in the uniform
topology. Consequently, every subsequence contains a further sub-subsequence
converging to the solution to the dynamical system 
\[
X\left( t\right) =X\left( 0\right) +\int_{0}^{t}b\left( X\left( s\right)
\right) ds,
\]%
which in turn has a unique solution because $b\left( \cdot \right) $ is
locally Lipschitz.

\subsection{\textbf{Proof of Theorem \ref{thm:Population_Based}}}
\label{ssec-proofs-population}

For now, let us assume that there exists a real $z>0$ such that $z_{n}=z n.$

We define $\bar{B}_{n}\left(  t\right)  =n^{-1}B_{n}\left(  t\right)$ and $\bar{S}_{n}\left(  t\right)  =n^{-1}S_{n}\left(  t\right)$. By~(\ref{eq1b})-(\ref{eq2b}) in the main text, we can write
\begin{align}
\bar{B}_{n}\left(  t\right)   &  =\bar{B}_{n}\left(  0\right)  +\frac
{N_{B}^{+}\left(  n\lambda t\right)  }{n}-\frac{N_{B}^{-}\left(  \eta\int%
_{0}^{t}B_{n}\left(  r\right)  dr\right)  }{n}\nonumber\\
&  -\frac{\int_{0}^{t}I_{\{\bar{B}_{n}\left(  r_{-}\right)  \geq z\}}
dN_{S}^{+}\left(  \lambda nr\right)  }{n}-\frac{\int_{0}^{t}I_{\{\bar
{S}_{n}\left(  r_{-}\right)  \geq z\}}  dN_{B}^{+}\left(  n\lambda
r\right) }{n},\label{eq:barB_n}\\
\bar{S}_{n}\left(  t\right)   &  =\bar{S}_{n}\left(  0\right)  +\frac
{N_{S}^{+}\left(  n\lambda r\right)  }{n}-\frac{N_{S}^{-}\left(  \eta\int%
_{0}^{t}S_{n}\left(  r\right)  dr\right)  }{n}\nonumber\\
&  -\frac{\int_{0}^{t}I_{\{\bar{S}_{n}\left(  r_{-}\right)  \geq z\}}
d\tilde{N}_{B}^{+}\left(  \lambda nr\right)  }{n}-\frac{\int_{0}^{t}I_{\{
\bar{B}_{n}\left(  r_{-}\right)  \geq z\}}  dN_{S}^{+}\left(  \lambda
nr\right)  }{n}.\label{eq:barS_n}%
\end{align}
We assume that $\bar{B}_{n}\left(  0\right)  \rightarrow\bar{B}\left(
0\right)  $ and $\bar{S}_{n}\left(  0\right)  \rightarrow S\left(  0\right)
$. Although formally the process in~(\ref{eq:barB_n})-(\ref{eq:barS_n})
converges to%
\begin{align}
\bar{B}\left(  t\right)    & =\bar{B}\left(  0\right)  +\lambda t-\eta\int%
_{0}^{t}\bar{B}\left(  r\right)  dr-\lambda\int_{0}^{t}I_{\{\bar{B}\left(
r\right)  \geq z\}}  dr-\lambda\int_{0}^{t}I_{\{\bar{S}\left(
r\right)  \geq z\}}  dr, \label{eq64a} \\
\bar{S}\left(  t\right)    & =\bar{S}\left(  0\right)  +\lambda t-\eta\int%
_{0}^{t}\bar{S}\left(  r\right)  dr-\lambda\int_{0}^{t}I_{\{\bar{B}\left(
r\right)  \geq z\}}  dr-\lambda\int_{0}^{t}I_{\{\bar{S}\left(
r\right)  \geq z\}}  dr, \label{eq64b}
\end{align}
this dynamical system is non-standard because the indicator functions are not
continuous. Hence, we need to study this system as the solution to a certain
Skorokhod problem. In particular, we can write%
\begin{align}
\bar{B}\left(  t\right)    & =\bar{B}\left(  0\right)  +\lambda t-\eta\int%
_{0}^{t}\bar{B}\left(  r\right)  dr-L_{z}^{\bar{B}}\left(  t\right)
-L_{z}^{\bar{S}}\left(  t\right)  , \label{eq62aa}\\
\bar{S}\left(  t\right)    & =\bar{S}\left(  0\right)  +\lambda t-\eta\int%
_{0}^{t}\bar{S}\left(  r\right)  dr-L_{z}^{\bar{B}}\left(  t\right)
-L_{z}^{\bar{S}}\left(  t\right)  , \label{eq62bb}
\end{align}
where $L_{z}^{\bar{B}}\left(  \cdot\right)  ,$ $L_{z}^{\bar{S}}\left(
\cdot\right)  $ are nondecreasing processes such that $L_{z}^{\bar{B}}\left(
0\right)  =L_{z}^{\bar{S}}\left(  0\right)  =0$ and%
\[
\int_{0}^{t}\left(  \bar{B}\left(  r\right)  -z\right)  dL_{z}^{\bar{B}%
}\left(  r\right)  =\int_{0}^{t}\left(  \bar{S}\left(  r\right)  -z\right)
dL_{z}^{\bar{S}}\left(  r\right) =0 ,
\]
and $\bar{B}\left(  t\right)  ,\bar{S}\left(  t\right)  \leq z$. Hence,
$L_{z}^{\bar{B}}\left(  \cdot\right)  $ and $L_{z}^{\bar{S}}\left(
\cdot\right)  $ are minimal nondecreasing processes that constrain the
dynamics of $\bar{B}\left(  \cdot\right)  $ and $\bar{S}\left(  \cdot\right)
$ to stay below $z$. 

The existence and uniqueness of the solution to the dynamical system in~(\ref{eq62aa})-(\ref{eq62bb}) is 
studied in~\S\ref{Section_Skorokhod}, which appears at the end of the
proof of this theorem. Also, we see that if $z<\lambda/\eta$, the equilibrium point of~(\ref{eq62aa})-(\ref{eq62bb}) is
$(\bar{B} \left(  \infty\right),\bar{S} \left(\infty\right))  =\left(  z,z\right).$

Note that
\begin{eqnarray}
\bar{B}_{n}\left(  t\right)   &=& \bar{B}_{n}\left(  0\right)  +\lambda
t-\eta\int_{0}^{t}\bar{B}_{n}\left(  r\right)  dr\nonumber\\
&&  -\frac{\int_{0}^{t}I_{\{\bar{B}_{n}\left(  r_{-}\right)  \geq z\}}
dN_{B}^{+}(\lambda nr)}{n}-\frac{\int_{0}^{t}I_{\{\bar{S}_{n}\left(
r_{-}\right)  \geq z\}}  N_{S}^{+}\left(  \lambda nr\right)}{n}  \nonumber\\
&&  +\frac{N_{B}^{+}\left(  n\lambda t\right)  -n\lambda t}{n}+\frac{\eta
\int_{0}^{t}n\bar{B}_{n}\left(  r\right)  dr-N_{B}^{-}\left(  \eta\int_{0}%
^{t}n\bar{B}_{n}\left(  r\right)  dr\right)}{n}\nonumber\\
&&  +\frac{\int_{0}^{t}I_{\{\bar{B}_{n}\left(  r_{-}\right)  \geq
z\}}  dN_{B}^{+}(\lambda nr)}{n}-\frac{\int_{0}^{t}I_{\{\bar{B}_{n}\left(
r_{-}\right)  \geq z\}}  dN_{S}^{+}(\lambda nr)}{n}\nonumber\\
&&  +\frac{\int_{0}^{t}I_{\{\bar{S}_{n}\left(  r_{-}\right)  \geq
z\}}  dN_{S}^{+}\left(  \lambda nr\right)}{n}  -\frac{\int_{0}^{t}I_{\{  \bar
{S}_{n}\left(  r_{-}\right)  \geq z\}}  dN_{B}^{+}\left(  n\lambda
r\right)}{n},\nonumber\\
&=& \bar{B}_{n}\left(  0\right)  +\lambda t-\eta\int_{0}^{t}\bar{B}_{n}\left(
r\right)  dr-L_{z}^{\bar{B}_{n}}\left(  t\right)  -L_{z}^{\bar{S}_{n}}\left(
t\right)  \nonumber\\
&&  +M_{n,1}^B\left(  t\right)  +M_{n,2}^B\left(  t\right)  +M_{n,3}^B\left(
t\right)  +M_{n,4}^B\left(  t\right),  \label{eq:Mart_Bn}%
\end{eqnarray}
where%
\begin{align*}
M_{n,1}^B\left(  t\right)   &  = \frac{N_{B}^{+}\left(  n\lambda t\right)
-n\lambda t}{n},\\
M_{n,2}^B\left(  t\right)   &  = \frac{\eta\int_{0}^{t}n\bar{B}_{n}\left(
r\right)  dr-N_{B}^{-}\left(  \eta\int_{0}^{t}n\bar{B}_{n}\left(  r\right)
dr\right)}{n},\\
M_{n,3}^B\left(  t\right)   &  =\frac{\int_{0}^{t}I_{\left\{\bar{B}_{n}\left(
r_{-}\right)  \geq z\right\}}  dN_{B}^{+}(\lambda nr)}{n}-\frac{\int_{0}^{t}I_{\left\{
\bar{B}_{n}\left(  r_{-}\right)  \geq z\right\}}  dN_{S}^{+}(\lambda nr)}{n},\\
M_{n,4}^B\left(  t\right)   &  =\frac{\int_{0}^{t}I_{\left\{  \bar{S}_{n}\left(
r_{-}\right)  \geq z\right\}}  dN_{S}^{+}\left(  \lambda nr\right)}{n}  -\frac{\int%
_{0}^{t}I_{\left\{  \bar{S}_{n}\left(  r_{-}\right)  \geq z\right\}}  dN_{B}%
^{+}\left(  n\lambda r\right)}{n}.
\end{align*}
Similarly, 
\begin{eqnarray}
\bar{S}_{n} \left( t\right) &=&\bar{S}_{n}\left( 0\right)
+\lambda t-\eta \int_{0}^{t}\bar{S}_{n} \left( r\right)
dr-  L_{z}^{\bar{B}_{n} }\left(
t\right) -L_{z}^{\bar{S}_{n} }\left( t\right) \nonumber \\
&&+M_{n,1}^S\left( t\right) +M_{n,2}^S\left( t\right)
+M_{n,3}^S\left( t\right) +M_{n,4}^S\left( t\right) .\label{eq:Mart_Sn}
\end{eqnarray}%
In~(\ref{eq:Mart_Bn})-(\ref{eq:Mart_Sn}), the processes $\left\{ M_{n,i}^B\right\} _{n\geq 1}$ and $\left\{
M_{n,i}^S\right\} _{n\geq 1}$ are martingales for $i=1,2,3,4$ such
that%
\[
E\left[\sup_{0\leq r\leq t}\left\vert M_{n,i}^B\left( r\right) \right\vert
^{2}\right]+E\left[\sup_{0\leq r\leq t}\left\vert M_{n,i}^S\left( r\right)
\right\vert ^{2}\right]=O\left( n^{-1}\right),
\]%
which is obtained by upper bounding $\bar{B}_{n} \left( \cdot
\right) $ by that of a system where no matches occur, i.e., agents leave only upon abandonment so that either side of the market can be modeled as an $M/M/\infty$ queue, and then applying Doob's maximal inequalities. Therefore, we have that 
\[
M_{n,i}^B,M_{n,i}^S\rightarrow 0 
\]%
uniformly on compact sets as $n\rightarrow \infty $ in probability for $%
i=1,2,3,4$.

The Lipschitz continuity property of the Skorokhod map, which is established in Section \S\ref{Section_Skorokhod}, implies that 
\[
\bar{S}_{n} \left( \cdot \right) \rightarrow \bar{S}%
 \left( \cdot \right) ,\text{ }\bar{B}_{n}\left( \cdot \right) \rightarrow \bar{B} \left( \cdot \right) 
\]%
uniformly on compact sets in probability.

The dynamical system describing $\left( \bar{B} ,\bar{S}%
 \right) $ has a unique attractor, which is the point $\left(
z,z\right) $ if $\lambda /\eta \geq z$, given the initial condition $\bar{B}%
 \left( 0\right) \leq z$, $\bar{S} \left(
0\right) \leq z$.

To show that the limit interchange ($t\rightarrow \infty $ and $n\rightarrow \infty $)
holds, we begin by upper bounding with a system with no matching, which implies that the steady-state number of buyers in the 
system is less than or equal to a Poisson random variable with rate $\lambda n/\eta $. It follows that $E\left[ \bar{B}%
_{n} \left( \infty \right) \right] \leq \lambda /\eta $, which in turn implies the uniformity property, 
\begin{equation}
\sup_{n\geq 1}E\left[ \bar{B}_{n} \left( \infty \right)
\right] <\infty .  \label{Bcinf}
\end{equation}%

Tightness of $\bar{B}_{n}\left( \infty \right) $ over $n$ follows from~(\ref{Bcinf}). Therefore, by Prohorov's theorem, every subsequence admits a further sub-subsequence which converges weakly; so, by selecting such sub-subsequence we may assume that $\bar{B}%
_{n} \left( 0\right) \rightarrow Z$. Moreover, by the Skorokhod embedding, we
can assume that this convergence occurs almost surely. 

Next, we have that $\bar{B}_{n} \left( \cdot \right)
\rightarrow \bar{B} \left( \cdot \right) $ on compact sets if $%
\bar{B}_{n} \left( 0\right) \rightarrow Z$ as $n\rightarrow
\infty $. Let us select $\bar{B}_{n} \left( 0\right) $ in
stationarity (i.e. equal in distribution to $\bar{B}_{n}\left( \infty \right) $).  By stationarity, for any fixed $t>0$, $%
\bar{B}_{n} \left( t\right) \rightarrow \bar{B}\left( t\right) =Z$ in probability and therefore $Z$ must be a stationary distribution of
the dynamic system $\bar{B} \left( \cdot \right) $. But the
stability point of $\bar{B} \left( \cdot \right) $ is unique and therefore $Z=z$. Thus, we must have that $\left( \bar{B}%
_{n}\left( \infty \right) ,\bar{S}_{n}\left( \infty \right) \right)
\rightarrow \left( z,z\right) $ almost surely as $n\rightarrow \infty $. Consequently, since the limit is independent of the subsequence, we conclude that we can exchange limits and expectations.

Our next goal is to compute the utility rate. Note that if $\left( \bar{B}_{n}\left( 0\right) ,\bar{S}_{n}\left( 0\right)
\right) $ follows the stationary distribution then 
taking expectations on both sides of equation~(\ref{eq1b}) of the main text yields 
\begin{equation}
\label{eq62a}
\eta E\left( \bar{B}_{n}\left( \infty \right) \right) =\lambda \{P\left( 
\bar{S}_{n}\left( \infty \right) <z\right) -P\left( \bar{B}_{n}\left( \infty
\right) \geq z\right) \}. 
\end{equation}
Observe that 
\begin{eqnarray}
I_{\{\bar{S}_{n}\left( \infty \right) <z\}} &=&I_{\{\bar{S}%
_{n}\left( \infty \right) <z,\bar{B}_{n}\left( \infty \right) \geq z\}}
+I_{\{\bar{S}_{n}\left( \infty \right) <z,\bar{B}_{n}\left( \infty
\right) <z\}} ,  \nonumber  \\
&=&I_{\{\bar{B}_{n}\left( \infty \right) \geq z\}} +I_{\{\bar{B}%
_{n}\left( \infty \right) <z,\bar{S}_{n}\left( \infty \right) <z\}} . 
\label{Lim_P}
\end{eqnarray}%
Equations~(\ref{eq62a})-(\ref{Lim_P}) imply that  
\begin{equation}
\label{eq63a}
\eta E\left( \bar{B}_{n}\left( \infty \right) \right) = \lambda P\left( \bar{B}%
_{n}\left( \infty \right) <z,\bar{S}_{n}\left( \infty \right) <z\right) . 
\end{equation}%
Taking the limit in~(\ref{eq63a}) as $n\rightarrow \infty$, we conclude that 
\begin{equation}
\label{eq69a}
\frac{\eta z}{\lambda }=\lim_{n\rightarrow \infty }P\left( \bar{B}_{n}\left(
\infty \right) <z,\bar{S}_{n}\left( \infty \right) <z\right) . 
\end{equation}%
Equation~(\ref{Lim_P}) also implies that
\begin{equation}
\label{eq69b}
I_{\{\bar{B}_{n}\left( \infty \right) \geq z\}} +I_{\{\bar{S}%
_{n}\left( \infty \right) \geq z\}} =1-I_{\{\bar{B}_{n}\left( \infty
\right) <z,\bar{S}_{n}\left( \infty \right) <z\}} . 
\end{equation}%
By symmetry, we conclude from~(\ref{eq69a})-(\ref{eq69b}) that
\begin{equation} \label{eq:Fluid_Stationary}
\lim_{n\rightarrow \infty }P\left( \bar{B}_{n}\left( \infty \right) \geq
z\right) =\frac{1}{2}\left( 1-\frac{\eta z}{\lambda }\right) . 
\end{equation}%
Equation~(\ref{eq:Fluid_Stationary}) allows us to compute the utility rate: 
\begin{eqnarray}
U_{n}^p(z_n) &=&2\lambda nE\left[ m\left( B_{n}\left( \infty \right) \right)
I_{\{B_{n}\left( \infty \right) \geq zn\}} \right], \nonumber\\
&=&2\lambda nE\left[ m\left( B_{n}\left( \infty \right) \right) |B_{n}\left(
\infty \right) \geq zn\right] P\left( B_{n}\left( \infty \right) \geq
zn\right), \nonumber\\
&=&\lambda nm\left( zn\right) \left( 1-\frac{\eta z}{\lambda }\right) \left(
1+o\left( 1\right) \right) \label{eq70a}
\end{eqnarray}%
as $n\rightarrow \infty $, where the last equality follows from the use of Lemma \ref{lemma:Expectation_Max} and the observation that 
\[
P\left( \left\vert B_{n}\left( \infty \right) -zn\right\vert \geq
\varepsilon n|B_{n}\left( \infty \right) \geq zn\right) \rightarrow 0 ~~{\rm as} ~~n\rightarrow \infty.
\]%

Thus, by equation~(\ref{RV_max}) in the main text and~(\ref{eq70a}), 
\begin{equation} \label{eq:Utility_Limit}
\frac{U_n^p(z_n)}{n m(n)} = \lambda z^\alpha \left( 1-\frac{\eta z}{\lambda }\right) \left(
1+o\left( 1\right) \right).
\end{equation}

We first consider the case where $\alpha > 0$. Maximizing the right side of \eqref{eq:Utility_Limit} with respect to $z$ implies that among all policies such that $z_n=\Omega(n)$, 
setting $z_n = z_* n$ is asymptotically optimal. 
Further, for $z_{n}=o\left( n\right) $, using the same technique based on the fluid analysis, one can show that $\frac{U_n^p(z_n)}{n m(n)} = o(t)$; details are omitted for brevity. This completes the proof for $\alpha > 0$.

%
Finally, we consider the case where $\alpha = 0$. To make the dependence of $B_n(\infty)$ on $z_n$ more explicit, for the rest of the proof we denote $B_n(\infty)$ as $B_n^{(z_n)}(\infty)$. 
Using the fluid limit analysis similar to above, for any sequence of thresholds $z_{n}$ that is $o\left( n\right) $, we can show that $E[B_n^{(z_n)}(\infty)]$ is $o(n)$, as follows. In the pre-limit, in particular in equations~\eqref{eq:barB_n}-\eqref{eq:barS_n} and~\eqref{eq:Mart_Bn}-\eqref{eq:Mart_Sn},  we replace $z$ with $z_n/n$, and then using essentially the same arguments as above we obtain the fluid limit where $\left( \bar{B}^{(z_n)}
_{n}\left( \infty \right) ,\bar{S}^{(z_n)}_{n}\left( \infty \right) \right)
\rightarrow \left( 0,0\right) $ almost surely. Further, by using the same arguments as those used to obtain \eqref{eq:Fluid_Stationary}, we can show that 
\begin{equation}\label{eq:ProbAlpha0}
\lim_{n\rightarrow \infty }P\left( {B}_{n}^{(z_n)}\left( \infty \right) \geq
z_n \right) =\frac{1}{2}.
\end{equation}
By symmetry, we also have
\begin{equation}
\lim_{n\rightarrow \infty }P\left( {S}_{n}^{(z_n)}\left( \infty \right) \geq
z_n \right) =\frac{1}{2}.
\end{equation}

Now consider the $n^{\text{th}}$ system, i.e., the system with the arrival rate of buyers equal to $\lambda n$. For $k=1,2,\ldots,$ define $I_k^{(n)}$ as follows: $I_k^{(n)}$ is equal to $1$ if the $k^{\text{th}}$ arrival of buyers sees at least $z_n$ sellers upon arrival and is equal to $0$ otherwise. Then PASTA (Poisson Arrivals See Time Averages) implies that 
$$\lim_{k\to \infty} \frac{1}{k} \sum_{j=1}^{k} I_j^{(n)} = \frac{1}{2} + o(1).$$ 
For each $k=1,2,\ldots,$ and for each $n = 1,2,\ldots$, let $R_k^{(n)}$ be an independent random variable with distribution $F_{z_n}$. Then, again by PASTA, we have
$$\lim_{k\to \infty} \frac{1}{k} \sum_{j=1}^{k} I_j^{(n)}R_j^{(n)} = m(z_n) \left( \frac{1}{2} + o(1) \right).$$ 
However, by Assumption~\ref{as:UtilityDependence} and by symmetry, with probability $1$ we have 
\begin{equation}
U_n^p(z_n) \ge 2 \lambda n \lim_{k\to \infty} \frac{1}{k} \sum_{j=1}^{k} I_j^{(n)}R_j^{(n)}.
\end{equation}
Thus, we have $U_n^p(z_n) \ge \lambda n m(z_n) \left( 1 + o(1) \right).$ Consequently, for any sequence $z_n=o(n)$ such that $$\lim_{n\to \infty} \frac{m(z_n)}{m(n)} = 1,$$ we would have that  $U_n^p(z_n) \ge \lambda n m(n) \left( 1 + o(1) \right)$. Lemma~\ref{lemma:VerySlowlyVarying} guarantees that such a sequence exists. 

Combining these results with the upper bound in Lemma~\ref{lemma:bound} completes the proof for $\alpha = 0$, and thus also the overall proof of Theorem \ref{thm:Population_Based}.

\subsubsection{Skorokhod Problem}\label{Section_Skorokhod}

In this subsection, we consider the existence and uniqueness of the system of differential equations%
\begin{align*}
\bar{B}\left(  t\right)   &  =\bar{B}\left(  0\right)  -\eta\int_{0}%
^{t}\left(  \bar{B}\left(  r\right)  -\frac{\lambda}{\eta}\right)  dr-L_{z}^{\bar{B}%
}\left(  t\right)  -L_{z}^{\bar{S}}\left(  t\right)  ,\\
\bar{S}\left(  t\right)   &  =\bar{S}\left(  0\right)  -\eta\int_{0}%
^{t}\left(  \bar{S}\left(  r\right)  -\frac{\lambda}{\eta}\right)  dr-L_{z}^{\bar{B}%
}\left(  t\right)  -L_{z}^{\bar{S}}\left(  t\right)  ,
\end{align*}
where $L_{z}^{\bar{B}}\left(  \cdot\right)  ,$ $L_{z}^{\bar{S}}\left(
\cdot\right)  $ are nondecreasing processes such that $L_{z}^{\bar{B}}\left(
0\right)  =L_{z}^{\bar{S}}\left(  0\right) =0$, 
\[
\int_{0}^{t}\left(  \bar{B}\left(  r\right)  -z\right)  dL_{z}^{\bar{B}%
}\left(  r\right)  =\int_{0}^{t}\left(  \bar{S}\left(  r\right)  -z\right)
dL_{z}^{\bar{S}}\left(  r\right) =0 ,
\]
and $\bar{B}\left(  t\right)  ,\bar{S}\left(  t\right)  \leq z$. Hence,
$L_{z}^{\bar{B}}\left(  \cdot\right)  $ and $L_{z}^{\bar{S}}\left(
\cdot\right)  $ are minimal nondecreasing processes that constrain the
dynamics of $\bar{B}\left(  \cdot\right)  $ and $\bar{S}\left(  \cdot\right)
$ to stay below $z$. 

In order to use explicit expressions for Skorokhod problems studied in the positive orthant, we introduce a change of coordinates. Defining $\bar{B}_{z}\left(  t\right)  =z-\bar{B}\left(  t\right)$, $\bar{S}_{z}\left(  t\right)  =z-\bar{S}\left(  t\right) $ and 
$\bar{\lambda}_{z}=\lambda/\eta-z\geq0$, we have that%
\begin{align*}
&  \bar{B}_{z}\left(  t\right)  =\bar{B}_{z}\left(  0\right)  -\eta
\int_{0}^{t}\left(  \bar{B}_{z}\left(  r\right)  +\bar{\lambda}_{z}\right)
dr+L\left(  t\right)  \ge 0,\\
&  \bar{S}_{z}\left(  t\right)  =\bar{S}_{z}\left(  0\right)  -\eta
\int_{0}^{t}\left(  \bar{S}_{z}\left(  r\right)  +\bar{\lambda}_{z}\right)
dr+L\left(  t\right)  \ge 0,
\end{align*}
where
\[
\int_{0}^{t}\min\left(  \bar{B}_{z}\left(  r\right)  ,\bar{S}_{z}\left(
r\right)  \right)  dL\left(  r\right)  =0,\text{ \ }L\left(  0\right)  =0.
\]
Further, if we define 
\[
Z(t) = \min\left(  \bar{B}_{z}\left(  0\right)  -\eta\int_{0}^{t}\left(  \bar
{B}_{z}\left(  r\right)  +\bar{\lambda}_{z}\right)  dr,\bar{S}_{z}\left(
0\right)  -\eta\int_{0}^{t}\left(  \bar{S}_{z}\left(  r\right)  +\bar{\lambda
}_{z}\right)  dr\right),
\]
then 
\begin{align*}
Y\left(  t\right)   &  :=\min\left(  \bar{B}_{z}\left(  t\right)  ,\bar{S}%
_{z}\left(  t\right)  \right) \\
&  =Z\left(  t\right)  +L\left(  t\right)  ,
\end{align*}
which implies that
\begin{align*}
Y\left(  t\right)   &  =Z\left(  t\right)  +\max_{0\leq s\leq t}\left\{
-Z\left(  s\right)  ,0\right\} , \\
L\left(  t\right)   &  =\max_{0\leq s\leq t}\left\{  -\bar{B}_{z}\left(
0\right)  +\eta\int_{0}^{s}\left(  \bar{B}_{z}\left(  r\right)  +\bar{\lambda
}_{z}\right)  dr,-\bar{S}_{z}\left(  0\right)  +\eta\int_{0}^{s}\left(
\bar{S}_{z}\left(  r\right)  +\bar{\lambda}_{z}\right)  dr,0\right\}  .
\end{align*}
We then obtain $L\left(  t;\bar{B}_{z},\bar{S}_{z}\right)  :=L\left(
t\right)  $, to emphasize the dependence on $\left(  \bar{B}_{z},\bar{S}%
_{z}\right)  $, yielding 
\begin{align}
\bar{B}_{z}\left(  t\right)   &  =\bar{B}_{z}\left(  0\right)  +\eta\int%
_{0}^{t}\left(  \bar{B}_{z}\left(  r\right)  +\bar{\lambda}_{z}\right)
dr+L\left(  t;\bar{B}_{z},\bar{S}_{z}\right)  ,\label{Sys_a}\\
\bar{S}_{z}\left(  t\right)   &  =\bar{S}_{z}\left(  0\right)  +\eta\int%
_{0}^{t}\left(  \bar{S}_{z}\left(  r\right)  +\bar{\lambda}_{z}\right)
dr+L\left(  t;\bar{B}_{z},\bar{S}_{z}\right)  .\label{eq75a}
\end{align}
 We need to show that (\ref{Sys_a})-(\ref{eq75a}) has a unique solution. We first argue
uniqueness. Assume that there exists another solution that we shall denote as
$\left(  \bar{B}_{z}^{\prime},\bar{S}_{z}^{\prime}\right)  $ and consider
$\Delta_{B}=\bar{B}_{z}-\bar{B}_{z}^{\prime}$ and $\Delta_{S}=\bar{S}_{z}%
-\bar{S}_{z}^{\prime}$. Suppose that $\Delta_{S}\left(  0\right)
=\Delta_{B}\left(  0\right) =0 $. Then%
\begin{align}
\Delta_{B}\left(  t\right)   &  =\eta\int_{0}^{t}\Delta_{B}\left(  r\right)
dr+L\left(  t;\bar{B}_{z},\bar{S}_{z}\right)  -L\left(  t;\bar{B}_{z}^{\prime
},\bar{S}_{z}^{\prime}\right), \label{DelRev1} \\
\Delta_{S}\left(  t\right)   &  =\eta\int_{0}^{t}\Delta_{S}\left(  r\right)
dr+L\left(  t;\bar{B}_{z},\bar{S}_{z}\right)  -L\left(  t;\bar{B}_{z}^{\prime
},\bar{S}_{z}^{\prime}\right). \label{DelRev2} 
\end{align}
Now, consider any real numbers $a,b,c,a^{\prime},b^{\prime},c^{\prime}$.
Suppose without loss of generality that $a=\max\left(  a,b,c\right)  \geq
\max\left(  a^{\prime},b^{\prime},c^{\prime}\right)  $. Then, since
$\max\left(  a^{\prime},b^{\prime},c^{\prime}\right)  \geq a^{\prime}$, we
conclude that
\[
0\leq a-\max\left(  a^{\prime},b^{\prime},c^{\prime}\right)  \leq a-a^{\prime
},
\]
which implies
\[
\left\vert \max\left(  a,b,c\right)  -\max\left(  a^{\prime},b^{\prime
},c^{\prime}\right)  \right\vert \leq\left\vert a-a^{\prime}\right\vert
+\left\vert b-b^{\prime}\right\vert +\left\vert c-c^{\prime}\right\vert .
\]
Therefore, if $D\left(  t\right)  =\left\vert \Delta_{B}\left(  t\right)
\right\vert +\left\vert \Delta_{S}\left(  t\right)  \right\vert $, we have, for example, that
\[
\left\vert L\left(  t;\bar{B}_{z},\bar{S}_{z}\right)  -L\left(  t;\bar{B}_{z}^{\prime
},\bar{S}_{z}^{\prime}\right)\right\vert
\leq   \eta\int_{0}^{s}\left\vert \Delta_{B}\left(  r\right)
\right\vert  dr + \eta\int_{0}^{s}\left\vert \Delta_{S}\left(  r\right)
\right\vert  dr. 
\]
As a result, adding together (\ref{DelRev1}) and (\ref{DelRev2}), and using the triangle inequality we conclude that 
\[
D\left(  t\right)  \leq2\eta\int_{0}^{t}D\left(  r\right)  dr.
\]
Because $D\left(  0\right)  =0$, a direct application of Gronwall's
inequality yields that $D\left(  t\right)  =0$ for all $t>0$, and uniqueness follows.

Now we argue the existence of a solution to (\ref{Sys_a})-(\ref{eq75a}). The construction
follows by applying a standard Picard iteration. Let
\begin{align*}
\mathfrak{B}_{t}\left(  B,S\right)   &  =B\left(  0\right)  +\eta\int_{0}%
^{t}\left(  B\left(  r\right)  +\bar{\lambda}_{z}\right)  dr+L\left(
t;B,S\right)  ,\\
\mathfrak{S}_{t}\left(  B,S\right)   &  =S\left(  0\right)  +\eta\int_{0}%
^{t}\left(  S\left(  r\right)  +\bar{\lambda}_{z}\right)  dr+L\left(
t;B,S\right)  .
\end{align*}
Observe that the map $\left(  B,S\right)  \rightarrow\left(  \mathfrak{B}%
,\mathfrak{S}\right)  $ is Lipschitz continuous with respect to the uniform
topology over any compact time interval $[0,T]$, using the corresponding
uniform metric 
\[
\left\Vert \left(  B,S\right)  \right\Vert _{[0,T]}=\sup_{0\leq t\leq
T}\left(  \left\vert B\left(  t\right)  \right\vert +\left\vert S\left(
t\right)  \right\vert \right)  .
\]
Define $B_{z}^{\left(  0\right)  }\left(  t\right)  =B_{z}^{\left(
0\right)  }\left(  0\right)  ,$ $S_{z}^{\left(  0\right)  }\left(  t\right)
=S_{z}^{\left(  0\right)  }\left(  0\right)  $, and iteratively, for $m\geq1$,
\[
B_{z}^{\left(  m\right)  }\left(  t\right)  =\mathfrak{B}_{t}(B_{z}^{\left(
m-1\right)  },S_{z}^{\left(  m-1\right)  });\text{ }S_{z}^{\left(  m\right)
}\left(  t\right)  =\mathfrak{S}_{t}(B_{z}^{\left(  m-1\right)  }%
,S_{z}^{\left(  m-1\right)  }).
\]
Noting that $B_{z}^{\left(  m\right)  }\left(  0\right)  =B_{z}^{\left(
m-1\right)  }\left(  0\right)  $ and $S_{z}^{\left(  m\right)  }\left(
0\right)  =S_{z}^{\left(  m-1\right)  }\left(  0\right)  $, we have 
\begin{align*}
&  |B_{z}^{\left(  m\right)  }\left(  t\right)  -B_{z}^{\left(  m-1\right)
}\left(  t\right)  |\smallskip+\smallskip|S_{z}^{\left(  m\right)  }\left(
t\right)  -S_{z}^{\left(  m-1\right)  }\left(  t\right)  |\\
&  \leq3\eta\int_{0}^{t}|B_{z}^{\left(  m-1\right)  }\left(  r\right)
-B_{z}^{\left(  m-2\right)  }\left(  r\right)  |+|S_{z}^{\left(  m-1\right)
}\left(  r\right)  -S_{z}^{\left(  m-2\right)  }\left(  r\right)  |dr.
\end{align*}
Consequently, 
we conclude that
\[
\left\Vert (B_{z}^{\left(  m\right)  },S_{z}^{\left(  m\right)  }%
)-(B_{z}^{\left(  m-1\right)  },S_{z}^{\left(  m-1\right)  })\right\Vert
_{[0,T]}\leq3\eta T\left\Vert (B_{z}^{\left(  m-1\right)  },S_{z}^{\left(
m-1\right)  })-(B_{z}^{\left(  m-2\right)  },S_{z}^{\left(  m-2\right)
})\right\Vert _{[0,T]}.
\]
Choosing $0<T<1/(3\eta)$ we can deduce -- by applying successive iterations and the
triangle inequality -- that $\{(B_{z}^{\left(  m\right)  },S_{z}^{\left(
m\right)  }):m\geq1\}$ forms a Cauchy sequence in the space of continuous
functions endowed with the uniform topology, which is a complete separable
metric space. Therefore, by continuity of the map $\left(  B,S\right)
\rightarrow\left(  \mathfrak{B},\mathfrak{S}\right)  $, the limiting sequence
must satisfy (\ref{Sys_a}). The construction can be applied sequentially to
consecutive intervals of size less than $1/(3\eta)$.

\subsection{\textbf{Proof of Theorem \ref{thm:utility_based}}}
\label{ssec-proofs-utility}

For now, we assume that the thresholds satisfy 
$$\frac{v_{n}}{m(n)} \to v \text{ for some } v \ge 0.$$
We consider a Poisson-flow representation of the scaled utility-based dynamics, whose validity is demonstrated in~\S\ref{Section_Markov}. In particular, it suffices to study the scaled processes $\bar{B}_{n}\left( t\right) =n^{-1}B_{n}\left( t\right)$ and $\bar{S}_{n}\left( t\right) =n^{-1}S_{n}\left( t\right)$, which give rise to the representation 
\begin{eqnarray}
\bar{B}_{n}\left( t\right) &=& \bar{B}_{n}\left( 0\right) +\frac{N_B^+\left(
\lambda n\int_{0}^{t}F_{S_{n}\left( r\right) }\left( v_{n}\right) dr\right)}{n}
-\frac{N_B^-\left( \eta \int_{0}^{t}B_{n}\left( r\right) dr\right)}{n} \nonumber
\\ && -\frac{\tilde{N}_S^+\left( \lambda n\int_{0}^{t}\left( 1-F_{B_{n}\left(
r\right) }\left( v_{n}\right) \right) dr\right)}{n} , \label{eq76a}
\end{eqnarray}

\begin{eqnarray}
\bar{S}_{n}\left( t\right) &=& \bar{S}_{n}\left( 0\right)
+ \frac{N_S^+\left( \lambda n\int_{0}^{t}F_{B_{n}\left( r\right)
}\left( v_{n}\right) dr\right)}{n} - \frac{N_S^-\left( \eta
\int_{0}^{t}S_{n}\left( r\right) dr\right)}{n} \nonumber \\
&& - \frac{\tilde{N}_B^+\left( \lambda n\int_{0}^{t}\left( 1-F_{S_{n}\left( r\right)
}\left( v_{n}\right) \right) dr\right)}{n}, \label{eq76b}
\end{eqnarray}%
where $N_B^+\left( \cdot \right)
,\tilde{N}_B^+\left( \cdot \right)
,N_B^-\left( \cdot \right) ,N_S^+\left(\cdot \right)
,\tilde{N}_S^+\left(\cdot \right),
N_S^-\left( \cdot \right)$ are independent Poisson processes with unit mean.

Note that under Assumption \ref{as:Max_Pareto}, we have 
\begin{eqnarray*}
F_{nw}\left( m\left( n\right) t\right) &=&P\left( M\left( nw\right) \leq
m\left( n\right) t\right), \\
&=&P\left( Xm\left( nw\right) \leq m\left( n\right) t\right) \left(
1+o\left( 1\right) \right), \\
&=&e^{-\kappa w/t^{1/\alpha}} \left( 1+o\left( 1\right)
\right) .
\end{eqnarray*}
Using this result, the putative fluid limit of~(\ref{eq76a})-(\ref{eq76b}) is given by 

$$
\bar{B}\left( t\right) =\bar{B}\left( 0\right)  +\lambda \int_{0}^{t} e^{-\kappa 
\bar{S}\left( r\right) /v^{1/\alpha }} -\eta \int_{0}^{t}%
\bar{B}\left( r\right) dr - \lambda \int_{0}^{t} \left( 1 -e^{-\kappa \bar{B}%
\left( r\right) /v^{1/\alpha }}\right)  dr, 
$$

$$
\bar{S}\left( t\right) = \bar{S}\left( 0\right) +\lambda \int_{0}^{t}e^{-\kappa 
\bar{B}\left( r\right) /v^{1/\alpha }} -\eta \int_{0}^{t} 
\bar{S}\left( r\right) dr  - \lambda \int_{0}^{t} \left(1-e^{-\kappa \bar{S}%
\left( r\right) /v^{1/\alpha }} \right) dr.
$$

We proceed in four steps. The first step is to obtain a martingale
decomposition similar to that given in the proof of Theorem \ref{thm:Population_Based} part ii).  The martingales will converge to zero on compact sets. The second step is to show that $\left\{\left(\bar{B}_n(\cdot),\bar{S}_n(\cdot)\right) \right\}_{n\ge 1}$ is tight in the Skorokhod topology using the technique developed in Ethier and Kurtz (2005). The third step is to show that the putative fluid limit has a unique solution. The fourth step is to show that the ordinary differential equation (ODE) describing the fluid limit has a unique stationary point. 

Together, these four steps imply that any subsequence of the sequence $\left\{\left(\bar{B}_n(\cdot),\bar{S}_n(\cdot)\right) \right\}_{n\ge 1}$ will contain a subsequence that converges (by tightness) to the unique solution of the above ODE, and therefore the fluid limit convergence holds. Further, since $\sup_n E[\bar{B}_n (\infty)] < \infty$ follows easily from the upper bound where no matching happens, we have that any subsequence of $\left\{(\bar{B}_n (\infty),\bar{S}_n(\infty))\right\}_{n\ge 1}$ will contain a subsequence that converges. 

Now, consider the stationary versions of the process $(B_n(.),S_n(.))$. By stationarity,  we have that if any subsequence of $\left\{(\bar{B}_n (\infty),\bar{S}_n(\infty))\right\}_{n\ge 1}$ converges then it has to converge to $(\bar z,\bar z)$, which is the unique stationary point of the dynamical system describing $(\bar{B} (t),\bar{S}(t))$. Thus, $\left\{(\bar{B}_n (\infty),\bar{S}_n(\infty))\right\}_{n\ge 1}$ converges to $(\bar z, \bar z)$.

{\bf Step 1:} Consider the martingales 
\begin{eqnarray}
M_{n,1}^B(t) &=& \frac{N_B^+\left(
\lambda n\int_{0}^{t}F_{S_{n}\left( r\right) }\left( v_{n}\right) dr\right)}{n} - \lambda \int_{0}^{t}F_{S_{n}\left( r\right) }\left( v_{n}\right) dr ,\nonumber \\
M_{n,2}^B(t) &=& \frac{N_B^-\left( \eta \int_{0}^{t}B_{n}\left( r\right) dr\right)}{n} -\eta \int_{0}^{t}%
\bar{B}_n\left( r\right) dr , \nonumber \\
M_{n,3}^B(t) &=&  \frac{\tilde{N}_S^+\left( \lambda n\int_{0}^{t}\left( 1-F_{B_{n}\left(
r\right) }\left( v_{n}\right) \right) dr\right)}{n} - \lambda \int_{0}^{t}\left( 1-F_{B_{n}\left(
r\right) }\left( v_{n}\right) \right) dr ,\nonumber
\end{eqnarray}
so that~(\ref{eq76a}) can be expressed as 

\begin{eqnarray} 
\bar B_n(t) &=& \bar B_n(0) + \lambda \int_{0}^{t}F_{S_{n}\left( r\right) }\left( v_{n}\right)  - \eta \int_{0}^{t}
\bar{B}_n\left( r\right) dr  - \lambda \int_{0}^{t}\left( 1-F_{B_{n}\left(
r\right) }\left( v_{n}\right) \right) dr \nonumber \\
&& + M_{n,1}^B(t) - M_{n,2}^B(t) - M_{n,3}^B(t).\nonumber
\end{eqnarray}

As argued in the proof of Theorem \ref{thm:Population_Based} part ii), using the upper bound on $B_n(t)$ vis-a-vis no matching and Doob's maximal inequality, we get 
$$E\left[ \sup_{0\leq r\leq t}\left\vert M_{n,i}\left( r\right) \right\vert
^{2}\right]=O\left( n^{-1}\right) .$$
Thus, for $i=1,2,3$, $M_{n,i}^B \to 0$ as $n \to \infty$ uniformly on compact sets. 

Further, because $F_{S_{n}\left( r\right) } $ lies in $[0,1]$ w.p.\ 1, from Assumption \ref{as:Max_Pareto} the process $ \int_{0}^{t}F_{S_{n}\left( r\right) }\left( v_{n}\right) dr - \int_{0}^{t} \exp \left( -\kappa 
\bar{S}_n\left( r\right) /v^{1/\alpha }\right)dr$ converges to $0$ uniformly over compact sets. Similarly, $$ \int_{0}^{t}\left( 1-F_{B_{n}\left(r\right) }\left( v_{n}\right) \right) dr - \int_{0}^{t} \left( 1 -e^{-\kappa \bar{B}_n\left( r\right) /v^{1/\alpha }} \right) dr$$ converges to $0$ uniformly over compact sets.

Analogously define martingales $M^S_{n,i}$ for $i=1,2,3$ for $\bar S_n$. Again, for $i=1,2,3$, $M^S_{n,i} \to 0$ as $n \to \infty$ uniformly on compact sets. 

Because $F_{S_{n}(r)}$ lies in $[0,1]$ w.p.1, the quadratic variations $[M_{n,1}(t)]$ and $[M_{n,3}(t)]$ are bounded from above by the quadratic variation of $  n^{-1}N_B^+\left(
\lambda nt \right)$, which tends to $0$ as $n\to \infty$. We now show that $[M_{n,2}(t)]$ also tends to $0$ as $n\to \infty$ for a given $t$. Note that $B_{n}(.)$ is bounded from above by the process corresponding to no matching, which has $O(n)$ mean. Thus, the mean number of jumps in $n^{-1}N_B^-\left( \eta \int_{0}^{t}B_{n}\left( r\right) dr\right)$ is $O(n)$. Further, the size of each jump is $1/n$ with probability $1$. Hence,  the quadratic variation $[M_{n,2}(t)]$ tends to $0$ as $n \to \infty$. 

Similarly, the quadratic variations of $M^S_{n,i}(t)$ for $i=1,2,3$ also tend to $0$ as $n\to \infty$. 

Thus, from Lemma~\ref{lemma:Proc_Limit} in~\S\ref{ssec-proofs-lemmas} we have that  
\[
\bar{B}_{n}\left( \cdot \right) \rightarrow \bar{B}\left( \cdot \right)  \text{ and }\bar{S}_{n}\left( \cdot \right) \rightarrow \bar{S}%
\left( \cdot \right) 
\]%
uniformly on compact sets in probability.

{\bf Step 2:} Using Theorem 7.2 of Chapter 3 in Ethier and Kurtz (2005), we show that the family $\left\{\left((\bar{B}_n(t): t\ge 0),(\bar{S}_n(t):t\ge 0)\right)\right\}_{n\ge 1}$ is tight in the Skorokhod topology. The first condition of Theorem 7.2 of Ethier and Kurtz (2005) holds easily since for each $t$ we have that $\bar B_n(t)$ and $\bar S_n(t)$ are both stochastically bounded from above by $1/n$ times Poisson distributed random variables with mean $\lambda n t$, each of which concentrates as $n \to \infty$. 

We now show that the second condition of Theorem 7.2 of Chapter 3 in Ethier and Kurtz (2005) holds as well. Note that, for each $n$, the times of positive jumps in $B_n$ are a subset of jump times in a Poisson  process of rate $\lambda n$. Also, the departures from $B_n$ may occur  either when a customer in $B_n$ abandons or when a customer arrives in $S_n$. Further, the times of customer abandonment is a subset of departure times in a $M/M/\infty$ queue, which -- due to the time-reversibility of the $M/M/\infty$ queue -- form a Poisson process with rate $\lambda n$ (plus a finite number of departures due to finite initial conditions). Thus, the times of negative jumps in $B_n$ are a subset of jump times in a Poisson process of rate $\lambda n$. Also, w.p.1, each jump is of size $1$. 

Thus, the modulus of continuity (see page 122 of Ethier and Kurtz 2005) of $\bar B_n$ is less than that of $n^{-1}A_n^B$ where $A_n^B$ is a Poisson process of rate $3 \lambda n$.  Similarly,  the modulus of continuity of $\bar S_n$ is less than that of $n^{-1}A^S_n$ where $A^S_n$ is a Poisson process of rate $3 \lambda n$. Therefore, it is enough to verify the second condition of Theorem 7.2 of Chapter 3 in Ethier and Kurtz (2005) for $\left\{\left((n^{-1}\bar{A}_n^B(t): t\ge 0),(n^{-1}\bar{A}^S_n(t):t\ge 0)\right)\right\}_{n\ge 1}$, which is easy to do.

{\bf Step 3:} Recall that $\bar{B}_{n}\left( 0\right) \rightarrow \bar{B}\left(
0\right) $ and $\bar{S}_{n}\left( 0\right) \rightarrow S\left( 0\right) $. Because the dynamical system describing $\bar{B}\left( t\right)$ and $\bar{S}\left( t\right)$ is an ODE with Lipschitz coefficients, it has a unique solution. 

\textbf{Step 4:} Using the fluid limit characterization, the fundamental
theorem of calculus, and symmetry, for each stationary solution to the ODE we have%

\begin{eqnarray}
0 &=& \lambda e^{-\kappa\bar{x} / v^{1/\alpha}} - \eta\bar{x} - \lambda
(1-e^{-\kappa\bar{x}/v^{1/\alpha}}),\nonumber \\
&=& - \eta\bar{x} - \lambda+ 2 \lambda e^{-\kappa\bar{x} /
v^{1/\alpha}}.\label{eq:zBar}
\end{eqnarray}
Solving~(\ref{eq:zBar}) for $v$, we define for $\bar{x}\in\left(  0,\lambda/\eta\right)
$,
\begin{equation}
v\left(  \bar{x}\right) = \left(\frac{\kappa \bar x}{\ln\left(\frac{2\lambda}{\eta \bar x+\lambda}\right)}\right)^{\alpha}.
\label{Eq_v_zbar}%
\end{equation}
We can also uniquely solve for $\bar{x}\left(  v\right)  $ for $v\in\left(
0,\infty\right)  $ by finding the inverse of (\ref{Eq_v_zbar}), which yields 
\begin{equation}
\bar{x}\left(  v\right)  =-\frac{\lambda}{\eta}+\frac{v^{1/\alpha}%
}{\kappa}W\left(  \frac{2\lambda\kappa}{\eta v^{1/\alpha}}\exp\left(
\frac{\lambda\kappa}{\eta v^{1/\alpha}}\right)  \right)  , \label{Pr_Eq_W}%
\end{equation}
where $W(\cdot)$ is the Lambert W function. Although we can work with either $v\in\left(  0,\infty\right)  $ or $\bar{x}%
\in\left(  0,\lambda/\eta\right)  $ when optimizing the asymptotic 
utility rate, it will be more convenient to optimize in terms of $\bar{x}$ and then
find the optimal utility threshold using (\ref{Eq_v_zbar}), as we now explain.

For each $v\geq0$, we have established that $\left\{(\bar{B}_{n}%
(\infty),\bar{S}_{n}(\infty))\right\}_{n\ge 1}$ converges to a unique $(\bar{x},\bar{x})$,
which can be characterized as above. Observe that $\bar{x}(v)\rightarrow
\lambda/\eta$ as $v\rightarrow\infty$. Also, again using \eqref{eq:zBar}, we
have that $\bar{x}(v)\rightarrow0$ as $v\rightarrow0$, which is same as
the limit under the greedy policy. Thus, if $v_{n}$ is $o(n)$ then
$\bar{B}_{n}(\infty)$ is also $o(n)$.



Assumption~\ref{as:Max_Pareto} implies that
\begin{align*}
&  E\left[  E[M\left(  B_{n}\left(  \infty\right)  \right)  I_{\left\{M\left(
B_{n}\left(  \infty\right)  \right)  \geq v_{n}\right\}}  |B_{n}\left(
\infty\right)  ] \right] \\
&  \sim E\left[  E[m\left(  B_{n}\left(  \infty\right)  \right)  X I_{\left\{
m\left(  B_{n}\left(  \infty\right)  \right)  X \geq v_{n}\right\}}
|B_{n}\left(  \infty\right)  ] \right]  ,\\
&  \sim m\left(  n\bar{x}(v)\right)  E\left[  X I_{\left\{m\left(  n\bar
{x}(v)\right)  X \geq v_{n}\right\}} \right], \\
&  \sim m\left(  n\bar{x}(v)\right)  E\left[  X I_{\left\{\bar{x}(v)^{\alpha} X
\geq v\right\}} \right].
\end{align*}
Thus, by~(\ref{eq6b}) in the main text, the utility rate satisfies
\begin{align*}
U_{n}^u(v_n)  &  =2\lambda nE\left[  E[M\left(  B_{n}\left(  \infty\right)  \right)
I_{\left\{M\left(  B_{n}\left(  \infty\right)  \right)  \geq v_{n}\right\}}
|B_{n}\left(  \infty\right)  ]\right] , \\
&  =2\lambda n m\left(  n\bar{x}(v)\right)  E\left[XI_{\left\{X\geq\frac{v}{\bar
{x}(v)^{\alpha}}\right\}}\right]\left(  1+o\left(  1\right)  \right)  .
\end{align*}
Assumption \ref{as:Regularly_Varying} implies that
$U_{n}^u(v_n)/m(n)$ is $o(n)$ if $v_{n}=o(n)$, and is $\Theta(n)$ if $v_{n}$ is
$\Theta(n)$. 
Furthermore, $E\left[XI_{\left\{X\geq v/\bar{x}(v)^{\alpha}\right\}}
\right]\rightarrow0$ if $v\rightarrow\infty$. Thus, $U^u_{n}(v_n)/m(n)$ is $o(n)$ if
$v_{n}=\omega(n)$ ($f(n)$ is $\omega(n)$ if there exist $c>0$ and an integer $n_o\ge 1$ such that $f(n)> cn$ for all integers $n\ge n_o$). Consequently, the optimal policy can be computed either as
\[
\sup_{v\in\left(  0,\infty\right)  }2\lambda\bar{x}\left(  v\right)  ^{\alpha
}E\left[ XI_{\left\{  X\geq\frac{v}{\bar{x}(v)^{\alpha}}\right\}}  \right]  ,
\]
or, in terms of $\bar{x}$, as
\begin{equation}
\sup_{\bar{x}\in\left(  0,\lambda/\eta\right)  }2\lambda\bar{x}^{\alpha
}E\left[XI_{\left\{X\geq\frac{v(\bar{x})}{\bar{x}^\alpha}\right\}}  \right]  ,
\label{Opt_z}%
\end{equation}
and we will solve~(\ref{Opt_z}). 

Recall that by Assumption~\ref{as:Max_Pareto}, $X=\left(  \kappa^{-1}T\right)
^{-\alpha}$, where $T$ is an exponential random variable with mean one. It follows that 
\begin{eqnarray}
E\left[  XI_{\left\{  X\geq\frac{v\left(  \bar{x}\right)  }{\bar{x}^{\alpha}%
}\right\}}  \right] &=& E\left[ \left(  \kappa^{-1}T\right)  ^{-\alpha
}I_{\left
\{  \frac{\gamma\bar{x}}{v\left(  \bar{x}\right)  ^{1/\alpha}}\geq
T\right\}}  \right] , \label{eq72a}\\
&=& \kappa^\alpha\int_{0}^{\kappa\bar{x}/v\left(  \bar{x}\right)  ^{1/\alpha}%
}t^{-\alpha}e^{-t} dt ,\label{eq72b}
\end{eqnarray}
where by~(\ref{Eq_v_zbar}) the upper integration limit in~(\ref{eq72b}) is 
\begin{equation}
\frac{\kappa\bar{x}}{v\left(  \bar{x}\right)  ^{1/\alpha}} = \ln\left(\frac{2\lambda}{\eta \bar x+\lambda}\right)  .
\label{eq72c}
\end{equation}

Note that the objective function of (\ref{Opt_z}) is zero at $\bar{x}=0$
and at $\bar{x}=\lambda/\eta$ (taking the limits from the left and right,
respectively). Therefore, since the right side of~(\ref{Opt_z}) is positive for $\bar
{x}\in\left(  0,\lambda/\eta\right)  $, any global maximizer in (\ref{Opt_z})
must be a stationary point. This, in turn, implies that any global maximizer
must satisfy
\begin{align*}
0 &  =\frac{d}{d\bar{x}}\bar{x}^{\alpha}\int_{0}^{\kappa\bar{x}/v\left(
\bar{x}\right)  ^{1/\alpha}}t^{-\alpha}e^{-t} dt ,\\
&  =\bar{x}^{\alpha}\frac{d}{d\bar{x}}\int_{0}^{\kappa\bar{x}/v\left(  \bar
{x}\right)  ^{1/\alpha}}t^{-\alpha}e^{-t} dt+\alpha\bar
{x}^{\alpha-1}\int_{0}^{\kappa\bar{x}/v\left(  \bar{x}\right)  ^{1/\alpha}%
}t^{-\alpha}e^{-t} dt.
\end{align*}
Equation~(\ref{eq72c}) implies that 
\begin{equation}
\frac{d}{d\bar{x}}\frac{\kappa\bar{x}}{v\left(  \bar{x}\right)  ^{1/\alpha}}=\frac
{d}{d\bar{x}}\ln\left(\frac{2\lambda}{\eta \bar x+\lambda}\right)  = -\frac{\eta}{\left(  \eta\bar{x}+\lambda\right)
}.\label{Der}%
\end{equation}
Therefore,%
\begin{eqnarray}
\frac{d}{d\bar{x}}\int_{0}^{\kappa\bar{x}/v\left(  \bar{x}\right)
^{1/\alpha}}t^{-\alpha}e^{-t} dt 
&=& \left( \frac{\kappa\bar{x}}{v\left(  \bar{x}\right)  ^{1/\alpha}}\right)
^{-\alpha}\exp\left(\ln\left(\frac{\eta \bar x+\lambda}{2\lambda}\right)\right)   \frac{d}{d\bar{x}}\frac{\kappa\bar{x}}{v\left(  \bar{x}\right)  ^{1/\alpha}} , \nonumber \\
&=& -\frac{v(\bar x)\eta}{2\lambda(\kappa\bar x)^{\alpha}} , \nonumber
\end{eqnarray}
which implies that
\begin{align*}
\alpha\bar{x}^{-1}\int_{0}^{\kappa\bar{x}/v\left(  \bar{x}\right)  ^{1/\alpha
}}t^{-\alpha}e^{-t} dt  & =\frac{d}{d\bar{x}}\int_{0}%
^{\kappa\bar{x}/v\left(  \bar{x}\right)  ^{1/\alpha}}t^{-\alpha}e^{-t} dt ,\\
& =\frac{v(\bar x)\eta}{2\lambda(\kappa\bar x)^{\alpha}}  ,
\end{align*}
or equivalently, 
\begin{equation}
 \frac{v\left(  \bar{x}\right) \eta}{2\lambda\kappa^{\alpha}}=\alpha
\bar{x}^{\alpha-1}\int_{0}^{\kappa\bar{x}/v\left(  \bar{x}\right)  ^{1/\alpha
}}t^{-\alpha}e^{-t}dt.\label{Eq_z_star}%
\end{equation}
Because $\alpha\in\left(  0,1\right)  $, the right side of~(\ref{Eq_z_star})  is decreasing and continuous in $\left(  0,\lambda
/\eta\right)  $, whereas by~(\ref{Der}) the left side of~(\ref{Eq_z_star}) is increasing in the same range. Moreover, the left side of~(\ref{Eq_z_star}) 
vanishes at zero and the right  side vanishes at $\lambda/\eta$. We
conclude that there a unique solution $x_{\ast}$ to (\ref{Eq_z_star}). Finally, the optimal policy is given by $v\left(  x_{\ast}\right)  $ in~(\ref{Eq_v_zbar}) and the
conclusion of the theorem follows directly by substituting in the
expresson for $x_{\ast}$ into (\ref{Opt_z}).

\subsubsection{Markov Dynamics of Utility-Based Process}\label{Section_Markov}

In this subsection, we show that the actual dynamics of the utility threshold policy are equivalent to the Poisson-flow representation given by
\begin{eqnarray}
B_{n}\left( t\right) &=&B_{n}\left( 0\right) +N_B^+\left( \lambda
n\int_{0}^{t}F_{S_{n}\left( r\right) }\left( v_{n}\right) dr\right)
-N_B^-\left( \eta \int_{0}^{t}B_{n}\left( r\right) dr\right) \nonumber \\
&&-\tilde{N}_S^+\left( \lambda n\int_{0}^{t}\left( 1-F_{B_{n}\left(
r\right) }\left( v_{n}\right) \right) dr\right) , \label{eq:BnUtilityBased} \\
S_{n}\left( t\right) &=&S_{n}\left( 0\right) +N_S^+\left( \lambda
n\int_{0}^{t}F_{B_{n}\left( r\right) }\left( v_{n}\right) dr\right)
-N_S^-\left( \eta \int_{0}^{t}S_{n}\left( r\right) dr\right) \nonumber \\
&&-\tilde{N}_B^+\left( \lambda n\int_{0}^{t}\left( 1-F_{S_{n}\left( r\right) }\left(
v_{n}\right) \right) dr\right), \label{eq:SnUtilityBased}
\end{eqnarray}%
where $N_B^+,\tilde{N}_B^+,N_S^+,\tilde{N}_S^+,N_B^-,N_S^-$ are all
independent Poisson processes with unit mean. For simplicity, we shall let $n=1$ and $\lambda=\eta=1$. 

Recall that the actual dynamics under the utility threshold policy are governed by the equations
\begin{align}
\widetilde{B}\left(  t\right)   &  =\widetilde{B}\left(  0\right)  +\sum
_{j=1}^{N_B^+\left(  t\right)  }I_{\left\{  \max_{i=1}^{\widetilde{S}\left(
A^B_{j-}\right)  }V_{i,j}^B\leq v\right\}}  -\sum_{j=1}^{N_S^+\left(
t\right)  }I_{\left\{  \max_{i=1}^{\widetilde{B}(A_{j-}^{S})}V_{i,j}%
^{S}>v\right\}} \label{MP0}\\
&  -N_B^-\left(  \int_{0}^{t}\widetilde{B}\left(  r_{-}\right)  dr\right)
,\nonumber\\
\widetilde{S}\left(  t\right)   &  =\widetilde{S}\left(  0\right)  +\sum
_{j=1}^{N_S^+\left(  t\right)  }I_{\left\{  \max_{i=1}^{\widetilde{B}%
(A_{j-}^{S})}V_{i,j}^{S}\leq v\right\}}  -\sum_{j=1}^{N_B^+\left(
t\right)  }I_{\left\{  \max_{i=1}^{\widetilde{S}(A_{j-}^{S})}V_{i,j}^B%
>v\right\}} \nonumber\\
&  -N_S^-\left(  \int_{0}^{t}\widetilde{S}\left(  r\right)
dr\right)  ,\label{MPss}
\end{align}
where 
$\left\{  V_{i,j}^B:i\geq1,j\geq1\right\}  $ and $\left\{  V_{i,j}^{S}:i\geq1,j\geq1\right\} $ form two independent arrays of i.i.d. random variables with CDF $F\left(
\cdot\right)  $. In addition, $\left\{  A_{j}^B:j\geq1\right\}  $ is the sequence of arrival times associated with 
$N_B^+$ and $\left\{  A_{j}^{S}:j\geq1\right\}  $ is the sequence of arrival times 
associated with $N_S^+$. Because of the mutual independence among the
$V_{i,j}^B$s, the $V_{i,j}^{S}$s and all of the unit rate Poisson
processes, $N_B^+,N_S^+,N_B^-,N_S^-$, it follows that the
process $\left(  \widetilde{B},\widetilde{S}\right)  $ is Markovian and is
non-explosive because each of its coordinates (i.e. $\widetilde{B}$ and
$\widetilde{S}$, respectively) can be bounded by independent infinite server
queues, simply by setting $v=\infty$. Consequently, this Markov process is well defined.

We have introduced a slight inconsistency in the notation in this subsection only, since we are now using $(\tilde{B},\tilde{S})$ to denote the actual dynamics. Ultimately, this is not important because, as our analysis in this subsection demonstrates, these are equivalent representations. The strategy consists of showing that the generators (or rate matrices) of the processes coincide.


Let $f:\mathbb{Z}_+ \times \mathbb{Z}_+ \rightarrow\mathbb{R}$ 
be any bounded function and note
that (by standard properties of the Poisson process),

\begin{align}
&  E\left[  f\left(  B\left(  h\right)  ,S\left(  h\right)  \right)  -f\left(
B\left(  0\right)  ,S\left(  0\right)  \right)  |B\left(  0\right)  ,S\left(
0\right)  \right] \nonumber \\
&  =E\left[  \int_{0}^{h}\left[f\left(  B\left(  t_{-}\right)  +1,S\left(
t_{-}\right)  \right)  -f\left(  B\left(  t_{-}\right)  ,S\left(
t_{-}\right)  \right)\right]  dN_B^+\left(  \int_{0}^{t}F_{S\left(  r\right)
}\left(  v\right)  dr\right)  |B\left(  0\right)  ,S\left(  0\right)  \right]
\nonumber\\
&  +E\left[  \int_{0}^{h}\left[ f\left(  B\left(  t_{-}\right)  -1,S\left(
t_{-}\right)  \right)  -f\left(  B\left(  t_{-}\right)  ,S\left(
t_{-}\right)  \right)  \right]  dN_B^-\left(  \int_{0}^{t}B\left(  r\right)
dr\right)  |B\left(  0\right)  ,S\left(  0\right)  \right] \nonumber\\
&  +E\left[  \int_{0}^{h}\left[  f\left(  B\left(  t_{-}\right)  -1,S\left(
t_{-}\right)  \right)  -f\left(  B\left(  t_{-}\right)  ,S\left(
t_{-}\right)  \right)  \right]  d\tilde{N}_S^+\left(  \int_{0}%
^{t}\left(  1-F_{B\left(  r\right)  }\left(  v\right)  \right)  dr\right)
|B\left(  0\right)  ,S\left(  0\right)  \right] \nonumber\\
&  +E\left[  \int_{0}^{h}\left[f\left(  B\left(  t_{-}\right)  ,S\left(
t_{-}\right)  +1\right)  -f\left(  B\left(  t_{-}\right)  ,S\left(
t_{-}\right)  \right)\right]  dN_S^+\left(  \int_{0}^{t}F_{B\left(
r\right)  }\left(  v\right)  dr\right)  |B\left(  0\right)  ,S\left(
0\right)  \right] \nonumber\\
&  +E\left[  \int_{0}^{h}\left[  f\left(  B\left(  t_{-}\right),S\left(
t_{-}\right)-1  \right)  -f\left(  B\left(  t_{-}\right)  ,S\left(
t_{-}\right)  \right)  \right]  dN_S^-\left(  \int_{0}^{t}S\left(
r\right)  dr\right)  |B\left(  0\right)  ,S\left(  0\right)  \right]
\nonumber\\
&  +E\left[  \int_{0}^{h}\left[  f\left(  B\left(  t_{-}\right) ,S\left(
t_{-}\right) - 1 \right)  -f\left(  B\left(  t_{-}\right)  ,S\left(
t_{-}\right)  \right)  \right]  d\tilde{N}_B^+\left(  \int_{0}^{t}\left(
1-F_{S\left(  r\right)  }\left(  v\right)  \right)  dr\right)  |B\left(
0\right)  ,S\left(  0\right)  \right] \nonumber\\
&  +o\left(  h\right) ~~{\rm as} ~~ h\rightarrow0.\label{EG}
\end{align}

Since
\begin{equation}
\mathfrak{M}_{0}\left(  t\right)  =N_B^+\left(  \int_{0}^{t}F_{S\left(
r\right)  }\left(  v\right)  dr\right)  -\int_{0}^{t}F_{S\left(  r\right)
}\left(  v\right)  dr,
\label{eq92a}
\end{equation}
is a martingale, we have that
\begin{align*}
&  E\left[  \int_{0}^{h}\left[f\left(  B\left(  t_{-}\right)  +1,S\left(
t_{-}\right)  \right)  -f\left(  B\left(  t_{-}\right)  ,S\left(
t_{-}\right)  \right)\right]  dN_B^+\left(  \int_{0}^{t}F_{S\left(  r\right)
}\left(  v\right)  dr\right)  |B\left(  0\right), S(0)  \right] \\
&  =E\left[  \int_{0}^{h}\left[f\left(  B\left(  t_{-}\right)  +1,S\left(
t_{-}\right)  \right)  -f\left(  B\left(  t_{-}\right)  ,S\left(
t_{-}\right)  \right) \right] F_{S\left(  t\right)  }\left(  v\right)  dt|B\left(
0\right), S(0)  \right] , \\
&  =\left[  f\left(  B\left(  0\right)  +1,S\left(  0\right)  \right)
-f\left(  B\left(  0\right)  ,S\left(  0\right)  \right)  \right]  F_{S\left(
0\right)  }\left(  v\right)  h+o\left(  h\right)   ~~{\rm as} ~~ h\rightarrow0 .
\end{align*}
Similarly, we can evaluate each of the expectations appearing in the right
side of (\ref{EG}); e.g., the second and third expectations are %
\begin{align}
&  E\left[  \int_{0}^{h}\left[f\left(  B\left(  t_{-}\right)  -1,S\left(
t_{-}\right)  \right)  -f\left(  B\left(  t_{-}\right)  ,S\left(
t_{-}\right)  \right) \right] dN_B^-\left(  \int_{0}^{t}B\left(  r\right)  dr\right)
|B\left(  0\right), S(0)  \right] \label{Rate}\\
&  +E\left[  \int_{0}^{h}\left[f\left(  B\left(  t_{-}\right)  -1,S\left(
t_{-}\right)  \right)  -f\left(  B\left(  t_{-}\right)  ,S\left(
t_{-}\right)  \right) \right] d\tilde{N}_S^+\left(  \int_{0}^{t}\left(
1-F_{B\left(  r\right)  }\left(  v\right)  \right)  dr\right)  |B\left(
0\right), S(0)  \right] \nonumber\\
&  =\left[f\left(  B\left(  0\right)  -1,S\left(  0\right)  \right)  -f\left(
B\left(  0\right)  ,S\left(  0\right) \right) \right]  \left[B\left(  0\right)  +\left(
1-F_{B\left(  0\right)  }\left(  v\right)  \right)  \right] h+o\left(  h\right)   ~~{\rm as} ~~ h\rightarrow0 .\nonumber
\end{align}

The above calculations show that the key to verifying that two well-defined Markov jump
processes with unit-size jumps are identical in law (i.e., have the same
generator) is showing that the corresponding compensators of the
associated point processes agree (i.e., they depend on the associated processes
in the same way). The corresponding martingales that identify
the compensators are, in addition to $\mathfrak{M}_{0}$ in~(\ref{eq92a}), 
\begin{align*}
\mathfrak{M}_{1}\left(  t\right)   &  =N_B^-\left(  \int_{0}^{t}B\left(
r\right)  dr\right)  -\int_{0}^{t}B\left(  r\right)  dr,\\
\mathfrak{M}_{2}\left(  t\right)   &  =\tilde{N}_S^+\left(  \int%
_{0}^{t}\left(  1-F_{B\left(  r\right)  }\left(  v\right)  \right)  dr\right)
-\int_{0}^{t}\left(  1-F_{B\left(  r\right)  }\left(  v\right)  \right)  dr,\\
\mathfrak{M}_{3}\left(  t\right)   &  =N_S^+\left(  \int_{0}%
^{t}F_{B\left(  r\right)  }\left(  v\right)  dr\right)  -\int_{0}%
^{t}F_{B\left(  r\right)  }\left(  v\right)  dr,\\
\mathfrak{M}_{4}\left(  t\right)   &  =N_S^-\left(  \int_{0}%
^{t}S\left(  r\right)  dr\right)  -\int_{0}^{t}S\left(  r\right)  dr,\\
\mathfrak{M}_{5}\left(  t\right)   &  =\tilde{N}_B^+\left(  \int_{0}^{t}\left(
1-F_{S\left(  r\right)  }\left(  v\right)  \right)  dr\right)  -\int_{0}%
^{t}\left(  1-F_{S\left(  r\right)  }\left(  v\right)  \right)  dr.
\end{align*}

To identify the corresponding compensators of the actual dynamics of the
utility threshold policy, we express these dynamics by means of a point
process representation and then compute the corresponding compensators
with respect to the $\sigma$-field generated by the population processes
(buyers and sellers).

We need to study the compensator of the point processes
\begin{equation}
\label{eq91a}
\sum_{j=1}^{N_B^+\left(  t\right)  }I_{\left\{  \max_{i=1}^{\widetilde{S}%
\left(  A_{j-}^B\right)  }V_{i,j}^B\leq v\right\}}  ,
\end{equation}
\begin{equation}
\label{eq91b}
\sum_{j=1}^{N_{S}^+\left(  t\right)  }I_{\left\{  \max_{i=1}^{\widetilde{B}(A_{j-}^{S}%
)}V_{i,j}^{S}>v\right\}} ,
\end{equation}
\begin{equation}
\label{eq91c}
\sum_{j=1}^{N_S^+\left(  t\right)  }I_{\left\{  \max_{i=1}%
^{\widetilde{B}(A_{j}^{S}-)}V_{i,j}^{S}\leq v\right\}}  ,
\end{equation}
\begin{equation}
\label{eq91d}
\sum
_{j=1}^{N_B^+\left(  t\right)  }I_{\left\{  \max_{i=1}^{\widetilde{S}%
(A_{j}^{B}-)}V_{i,j}^B>v\right\}}  ,
\end{equation}
with respect to the filtration generated by the processes $\widetilde{B}%
\left(  \cdot\right)  $ and $\widetilde{S}\left(  \cdot\right)  $, which we
shall denote as $\mathcal{G=}\left\{  \mathcal{G}_{t}:t\geq0\right\}  $. We
claim that for any $t,r>0$, the conditional expectation of~(\ref{eq91a}) can be expressed as 
\[
E\left[ \sum_{j=N_B^+\left(  t\right)  +1}^{N_B^+\left(  t+r\right)
}I_{\left\{\max_{i=1}^{\widetilde{S}\left(  A_{j}^B-\right)  }V_{i,j}^B\leq
v\right\}}  |\mathcal{G}_{t}\right]  =E\left[  \int_{t}^{t+r}F_{\widetilde{S}%
\left(  u\right)  }\left(  v\right)  du|\mathcal{G}_{t}\right]  .
\]
By Fubini's theorem, because every term is nonnegative in the second equality in the following display,
we have that
\begin{eqnarray}
E\left[ \sum_{j=N_B^+\left(  t\right)  +1}^{N_B^+\left(  t+r\right)
}I_{\left\{  \max_{i=1}^{\widetilde{S}\left(  A_{j}^B-\right)  }V_{i,j}^B\leq
v\right\}}  |\mathcal{G}_{t}\right] 
&=& E\left[  \sum_{j=1}^{\infty}I_{\left\{  t<A_{j}^B\leq t+r\right\}}  I_{\left\{
\max_{i=1}^{\widetilde{S}\left(  A_{j}^B-\right)  }V_{i,j}^B\leq v\right\}}
|\mathcal{G}_{t}\right] , \nonumber \\
&=& \sum_{j=1}^{\infty}E\left[  I_{\left\{ t<A_{j}^B\leq t+r\right\}}  I_{\left\{
\max_{i=1}^{\widetilde{S}\left(  A_{j}^B-\right)  }V_{i,j}^B\leq v\right\}}
|\mathcal{G}_{t}\right] .\nonumber
\end{eqnarray}

Note that for each $j\geq1$, $A_{j}^B$ is a stopping time, and therefore, by the
tower property in the second equality,
\begin{align*}
&  E\left[  I_{\left\{ t<A_{j}^B\leq t+r\right\}} I_{\left\{  \max_{i=1}%
^{\widetilde{S}\left(  A_{j}^B-\right)  }V_{i,j}^B\leq v\right\}}  |\mathcal{G}%
_{t}\right] \\
&  =\int_{t}^{t+r}E\left[  I_{\left\{  \max_{i=1}^{\widetilde{S}\left(
u_{-}\right)  }V_{i,j}^B\leq v\right\}}  |\mathcal{G}_{t},A_{j}^B=u\right]  P\left(
A_{j}^B\in du|\mathcal{G}_{t}\right) , \\
&  =\int_{t}^{t+r}E\left[  E\left[  I_{\left\{ \max_{i=1}^{\widetilde{S}\left(
u_{-}\right)  }V_{i,j}^B\leq v\right\}}  |\mathcal{G}_{u-},A_{j}^B=u\right]
|\mathcal{G}_{t},A_{j}^B=u\right]  P\left(  A_{j}^B\in du|\mathcal{G}_{t}\right) ,
\\
&  =\int_{t}^{t+r}E\left[  F_{\widetilde{S}\left(  u_{-}\right)  }\left(
v\right)  |\mathcal{G}_{t},A_{j}^B=u\right]  P\left(  A_{j}^B\in du|\mathcal{G}%
_{t}\right) ,\\
&  =E\left[  I_{\left\{  t<A_{j}^B\leq t+r\right\}}  F_{\widetilde{S}\left(
A_{j-}^B\right)  }\left(  v\right)  |\mathcal{G}_{t}\right]  .
\end{align*}
Applying Fubini's theorem again and summing over $j$, we conclude that
\[
E\left[  \sum_{j=N_B^+\left(  t\right)  +1}^{N_B^+\left(  t+r\right)
}I_{\left\{ \max_{i=1}^{\widetilde{S}\left(  A_{j}^B-\right)  }V_{i,j}^B\leq
v\right\}}  |\mathcal{G}_{t}\right]  = E\left[  \int_{t}^{t+r}F_{\widetilde{S}%
\left(  u_{-}\right)  }\left(  v\right)  dN_B^+\left(  u\right)
|\mathcal{G}_{t}\right]  .
\]
However, we have that the compensator of $N_B^+\left(  \cdot\right)  $ is the
identity and therefore%
\[
E\left[  \sum_{j=N_B^+\left(  t\right)  +1}^{N_B^+\left(  t+r\right)
}I_{\left\{ \max_{i=1}^{\widetilde{S}\left(  A_{j}^B-\right)  }V_{i,j}^B\leq
v\right\}}  |\mathcal{G}_{t}\right] = E\left[  \int_{t}^{t+r}F_{\widetilde{S}%
\left(  u_{-}\right)  }\left(  v\right)  du|\mathcal{G}_{t}\right]  .
\]
We conclude that
\[
\widetilde{\mathfrak{M}}_{0}\left(  t\right)  =\sum_{j=1}^{N_B^+\left(
t\right)  }I_{\left\{ \max_{i=1}^{\widetilde{S}\left(  A_{j-}^B\right)  }%
V_{i,j}^B\leq v\right\}}  -\int_{0}^{t}F_{\widetilde{S}\left(  u\right)  }\left(
v\right)  du
\]
is a martingale and it proves the corresponding compensator with respect to
$\mathcal{G}$. A completely analogous development can be obtained for the
point processes in~(\ref{eq91b})-(\ref{eq91d}), resulting in the martingales 
\begin{align*}
\widetilde{\mathfrak{M}}_{1}\left(  t\right)   &  =N_B^-\left(  \int_{0}%
^{t}\widetilde{B}\left(  r\right)  dr\right)  -\int_{0}^{t}\widetilde{B}%
\left(  r\right)  dr ,\\
\widetilde{\mathfrak{M}}_{2}\left(  t\right)   &  =\sum_{j=1}^{N_S^+\left(  t\right)  }I_{\left\{ \max_{i=1}^{\widetilde{B}(A_{j-}^{S}%
)}V_{i,j}^{S}>v\right\}}  -\int_{0}^{t}\left(  1-F_{\widetilde{B}\left(
r\right)  }\left(  v\right)  \right)  dr ,\\
\widetilde{\mathfrak{M}}_{3}\left(  t\right)   &  =\sum_{j=1}^{N_S^+\left(  t\right)  }I_{\left\{  \max_{i=1}^{\widetilde{B}(A_{j-}^{S}%
)}V_{i,j}^{S}\leq v\right\}}  -\int_{0}^{t}F_{\widetilde{B}\left(
r\right)  }\left(  v\right)  dr ,\\
\widetilde{\mathfrak{M}}_{4}\left(  t\right)   &  =N_S^-\left(
\int_{0}^{t}\widetilde{S}\left(  r\right)  dr\right)  -\int_{0}^{t}%
\widetilde{S}\left(  r\right)  dr ,\\
\widetilde{\mathfrak{M}}_{5}\left(  t\right)   &  =\sum_{j=1}^{N_B^+\left(
t\right)  }I_{\left\{ \max_{i=1}^{\widetilde{S}(A_{j-}^{B})}V_{i,j}^B%
>v\right\}}  -\int_{0}^{t}\left(  1-F_{\widetilde{S}\left(  r\right)  }\left(
v\right)  \right)  dr.
\end{align*}
This implies, by the reasoning given right after (\ref{Rate}) and comparing
$\widetilde{\mathfrak{M}}_{i}$ vs $\mathfrak{M}_{i}$ for $i\in\{0,1,\ldots,5\}$, 
that (\ref{eq:BnUtilityBased})-(\ref{eq:SnUtilityBased}) and (\ref{MP0})-(\ref{MPss}) are equivalent.

\subsection{\textbf{Proof of Theorem \ref{them}}}
\label{ssec-proofs-unbalanced}

We assume that the thresholds satisfy 
\begin{eqnarray}
	\frac{v_{n,b}}{m(n)}\rightarrow v_b, ~~\frac{v_{n,s}}{m(n)}\rightarrow v_s\text{ for some }v_b,v_s\geq0.\label{eq41a}
\end{eqnarray}
Note that $v_b,v_s$ could be infinity.
Define the scaled quantity $\bar{B}_n(t)=n^{-1}B_n(t),\bar{S}_n(t)=n^{-1}S_n(t).$ According to the analysis in~\S\ref{Section_Markov}, it suffices to express $\bar{B}_n(\cdot),\bar{S}_n(\cdot)$ in a Poisson-flow representation of the scaled system:
\begin{eqnarray}
\bar{B}_{n}\left( t\right) &=& \bar{B}_{n}\left( 0\right) +\frac{N_B^+\left(
\lambda_b n\int_{0}^{t}F_{S_{n}\left( r\right) }\left( v_{n,s}\right) dr\right)}{n}
-\frac{N_B^-\left( \eta_b \int_{0}^{t}B_{n}\left( r\right) dr\right)}{n} \nonumber
\\ && -\frac{\tilde{N}_S^+\left( \lambda_s n\int_{0}^{t}\left( 1-F_{B_{n}\left(
r\right) }\left( v_{n,b}\right) \right) dr\right)}{n} , \label{eq:B_bar_n}
\end{eqnarray}

\begin{eqnarray}
\bar{S}_{n}\left( t\right) &=& \bar{S}_{n}\left( 0\right)
+ \frac{N_S^+\left( \lambda_s n\int_{0}^{t}F_{B_{n}\left( r\right)
}\left( v_{n,b}\right) dr\right)}{n} - \frac{N_S^-\left( \eta_s
\int_{0}^{t}S_{n}\left( r\right) dr\right)}{n} \nonumber \\
&& - \frac{\tilde{N}_B^+\left( \lambda_b n\int_{0}^{t}\left( 1-F_{S_{n}\left( r\right)
}\left( v_{n,s}\right) \right) dr\right)}{n}, \label{eq:S_bar_n}
\end{eqnarray}%

Define the putative fluid limit 
\begin{equation}
\bar{B}\left( t\right) =\bar{B}\left( 0\right) +\lambda _{b}\int_{0}^{t}e^{-\kappa \bar{S}\left( r\right) /v_{s}^{1/\alpha }} -\eta
_{b}\int_{0}^{t}\bar{B}\left( r\right) dr 
-\lambda _{s}\int_{0}^{t}\left( 1-e^{-\kappa \bar{B}\left( r\right)
/v_{b}^{1/\alpha}}\right) dr,
\label{eq:fluid_b}
\end{equation}
\begin{equation}
\bar{S}\left( t\right) =\bar{S}\left( 0\right) +\lambda _{s}\int_{0}^{t}e^{-\kappa \bar{B}\left( r\right) /v_{b}^{1/\alpha }} -\eta
_{s}\int_{0}^{t}\bar{S}\left( r\right) dr \\
-\lambda _{b}\int_{0}^{t}\left( 1-e^{-\kappa \bar{S}\left( r\right)
/v_{s}^{1/\alpha }}\right) dr.
\label{eq:fluid_s}
\end{equation}
We want to show that $(\bar{B}_n,\bar{S}_n)$ converges weakly to $(\bar{B},\bar{S})$, and then by studying the stationary point of $(\bar{B},\bar{S})$, get the distribution of $\lim_{n\rightarrow\infty}(\bar{B}_n(\infty),\bar{S}(\infty))$ for use in the asymptotic analysis of the utility rate. We divide the proof into six steps. The first step is to show that $(\bar{B}_n(\cdot),\bar{S}_n(\cdot))\rightarrow (\bar{B}(\cdot),\bar{S}(\cdot))$ uniformly on compact sets in probability. The second step is to show the tightness of  $\{\bar{B}_n(\cdot),\bar{S}_n(\cdot)\}$ in the Skorokhod topology. The third step is the existence and uniqueness of the solution of the putative fluid limit in (\ref{eq:fluid_b})-(\ref{eq:fluid_s}). These first three steps follow directly from the proof of Theorem~\ref{thm:utility_based} in~\S\ref{ssec-proofs-utility}, and are omitted. The fourth step is to show that there exists a unique stationary point of the fluid limit and conclude that $(\bar{B}_n(\infty),\bar{S}_n(\infty))\rightarrow (b,s)$ in probability, where $(b,s)$ is the stationary point of the fluid limit (\ref{eq:fluid_b})-(\ref{eq:fluid_s}). The fifth step is to determine the  limit of the utility rate, and the last step is to compute the asymptotically optimal  thresholds.

\paragraph{Step 4} Any stationary solution $(b,s)$ of \eqref{eq:fluid_b}-\eqref{eq:fluid_s} must satisfy 
\begin{eqnarray}
\lambda _{b}e^{-\kappa s/v_{s}^{1/\alpha
}} +\lambda _{s}e^{-\kappa b/v_{b}^{1/\alpha }} &=& \eta _{b}b+\lambda _{s}, \label{eq:stationary point1} \\
\lambda _{b}e^{-\kappa s/v_{s}^{1/\alpha
}} +\lambda _{s}e^{-\kappa b/v_{b}^{1/\alpha }} &=& \eta _{s}s+\lambda _{b}. \label{eq:stationary point2}
\end{eqnarray}%
By~(\ref{eq:stationary point1})-(\ref{eq:stationary point2}), we know that 
\begin{align*}
	b=\frac{\eta _ss+\lambda_b-\lambda_s}{\eta_b}.
\end{align*}
It suffices to show that (\ref{eq:stationary point2}) has a unique solution. 
Because the right side of~(\ref{eq:stationary point2}) is a monotonically increasing function and the left side of~(\ref{eq:stationary point2}) is a monotonically decreasing function, as well as the facts that the left side exceeds the right side when $s=0$ and the right side exceeds the left side when $s=\lambda_s/\eta_s$, we conclude that the fluid limit has a unique stationary point.

We now show that
$(\bar{B}_n(\infty),\bar{S}_n(\infty))\rightarrow (b,s)$ in probability. First, by tightness, any subsequence of $(B(\cdot),S(\cdot))$ contains a subsubsequence that converges weakly. Because that subsubsequence converges to $(\bar{B},\bar{S})$ in probability uniformly on compact sets, the stationary point converges to $(b,s)$ in distribution. Hence, the equilibrium of the whole sequence $(\bar{B}_n(\infty),\bar{S}_n(\infty))\rightarrow (b,s)$ in probability.


\paragraph{Step 5} 
We want to show that the  limit of the utility rate is
\begin{equation}
\lim_{n\rightarrow\infty} \frac{U_n^u(v_{n,b},v_{n,s})}{nm(n)} 
		= \lambda_s b^\alpha E[XI_{\{b^\alpha X>v_b\}}]+\lambda_b s^\alpha E[XI_{\{s^\alpha X>v_s\}}], \label{eq:utilityrate_asymptotic}
\end{equation}
where the utility rate $U_n(v_{n,b},v_{n,s})$ is defined in~(\ref{unbalanced-utility}) in the main text.
To prove~(\ref{eq:utilityrate_asymptotic}), it suffices to show that there exists a subsubsequence of any subsequence of $U_n(v_{n,b},v_{n,s})/(nm(n))$ converging to that limit.
Because  $(\bar{B}_n(\infty),\bar{S}_n(\infty))\rightarrow (b,s)$ in probability, for any subsequence  there exists a subsubsequence that converges almost surely. Pick such a subsubsequence, by slightly abusing the notation, for further analysis, we denote such a subsubsequence as $(\bar{B}_n(\infty),\bar{S}_n(\infty)),$ which also implies that $(B_n(\infty),S_n(\infty))\rightarrow(\infty,\infty )$ almost surely. 
 By Assumption~\ref{as:Max_Pareto}, we have $\frac{M(B_n(\infty))}{m(B_n(\infty))}\Rightarrow X$. Combined with the fact that $\frac{m(B_n(\infty))}{m(nb)}\rightarrow 1$ almost surely, we obtain 
 \begin{align}\label{converge_z}
 	\frac{M(B_n(\infty))}{m(nb)}\Rightarrow X.
 \end{align}
Because the Lebesgue measure of the discontinuity point of the indicator function $I$ is $0$, we can apply the continuous mapping theorem and obtain
\begin{equation}
\frac{M(B_n(\infty))}{m(nb)}I_{\{M(B_n(\infty)\geq v_{n,b} \}}
=\frac{M(B_n(\infty))}{m(nb)}I_{\left\{\frac{M(B_n(\infty))}{m(nb)}\frac{m(nb)}{m(n)} \geq \frac{v_{n,b}}{m(n)} \right\} }\Rightarrow XI_{\{b^\alpha X\geq v_b\}}\label{converge_I}.
\end{equation}

Because
\begin{align*}
	E\left[\frac{M(B_n(\infty))}{m(nb)}\right]&=\frac{E[m(B_n(\infty)]}{m(E[B_n(\infty)])}\frac{m(E[B_n(\infty)])}{m(bn)},
\end{align*}
by Lemma~\ref{lemma:Expectation_Max}, it follows that $\frac{E[m(B_n(\infty)]}{m(E[B_n(\infty)])}\rightarrow 1.$ Because $\bar{B}_n(\infty)$ is bounded above by the queueing system without matching, which has expectation $\lambda_b/\eta_b<	\infty,$ we have
$ E\left[\bar{B}_n(\infty)\right]\rightarrow b$ by dominated convergence, which implies that $\frac{m(E[B_n(\infty)])}{m(bn)}\rightarrow 1.$ Thus, we have $\lim_{n\rightarrow\infty} E\left[\frac{M(B_n(\infty))}{m(nb)}\right]=1$.
Because 
\[0\leq \frac{M(B_n(\infty))}{m(nb)}I_{\{M(B_n(\infty)\geq v_{n,b} \}}\leq \frac{M(B_n(\infty))}{m(nb)}, \] 
by dominated convergence, we have
\begin{align*}
	\lim_{n\rightarrow\infty} E\left[\frac{M(B_n(\infty))}{m(nb)}I_{\{M(B_n(\infty)\geq v_{n,b} \}}\right]= E[XI_{\{b^\alpha X\geq v_b\}}].
\end{align*}
Similarly, we can show 
\begin{align*}
	\lim_{n\rightarrow\infty} E\left[\frac{M(S_n(\infty))}{m(ns)}I_{\{M(S_n(\infty)\geq v_{n,s} \}}\right]= E[XI_{\{s^\alpha X\geq v_s\}}].
\end{align*}
Note that by Assumption~\ref{as:Regularly_Varying},  $\lim_{n\to \infty}\frac{m(nb)}{m(n)}=b^\alpha, \lim_{n\to\infty}\frac{m(ns)}{m(n)}=s^\alpha$. It follows that 
\[\lim_{n\rightarrow\infty} E\left[\frac{M(B_n(\infty))}{m(n)}I_{\{M(B_n(\infty)\geq v_{n,b} \}}\right]= b^\alpha E[XI_{\{b^\alpha X\geq v_b\}}], \] which implies that asymptotic utility rate in~(\ref{eq:utilityrate_asymptotic}).

\paragraph{Step 6} 
By steps 1 to 5, we have shown that the asymptotically optimal thresholds can be obtained by solving the optimization problem (\ref{opt:theorem})-(\ref{opt:condition_2}) in Theorem~\ref{them}. In this final step, we compute the optimal thresholds by reducing (\ref{opt:theorem})-(\ref{opt:condition_2}) to a one-dimensional optimization problem over a compact interval by making a change of variables, showing that $v_b=v_s$ for any feasible pair $(b,s)$, and solving for the optimal $(b,s)$. 

The case of no matching, i.e., $v_b=v_s=\infty$, provides an upper bound on the system's steady-state population of $\rho_b=\frac{\lambda_b}{\eta_b}$ and $\rho_s=\frac{\lambda_s}{\eta_s}$. Hence, $b\in \lbrack 0,\rho _{b}]$ and $s\in \lbrack 0,\rho _{s}]$. Recall by Assumption~\ref{as:Max_Pareto} that $X=\left( \kappa ^{-1}T\right) ^{-\alpha }$, where $T$ is exponentially distributed with unit mean. Defining
\[
G\left( x\right) =\int_{0}^{x}t^{-\alpha }e^{-t} dt
\]%
and using the change of variables 
\begin{equation}
\label{eq46a}
x=\frac{\kappa s}{v_s^{1/\alpha}} ~~{\rm and} ~~y=\frac{\kappa b}{v_b^{1/\alpha}},
\end{equation}
we use~(\ref{eq72a})-(\ref{eq72b}) to express problem~(\ref{opt:theorem})-(\ref{opt:condition_2}) as \begin{eqnarray}
\max_{s,b\in \left( 0,\rho _{s}\right) \times \left( 0,\rho _{b}\right)
,\left( x,y\right) \in R_{+}^{2}} & & \lambda _{b}s^{\alpha }G\left( x\right)+\lambda _{s}b^{\alpha }G\left( y\right)
 \label{eq47}\\
{\rm subject ~to} & & 
\lambda _{b}e^{-x} +\lambda _{s}e^{-y} =b\eta _{b}+\lambda _{s} =s\eta
_{s}+\lambda _{b}.\label{eq48}
\end{eqnarray}

Suppose for now that $b$ and $s$ are fixed (and feasible) and we are optimizing
over $\left( x,y\right)$ in~(\ref{eq46a}). Because $G\left( \cdot \right)$ in~(\ref{eq47}) is strictly
concave and increasing, problem~(\ref{eq47})-(\ref{eq48}) is given by 
\begin{eqnarray}
\max_{\left( x,y\right) \in R_{+}^{2}} & & \lambda _{b}s^{\alpha }G\left( x\right) + \lambda _{s}b^{\alpha }G\left(
y\right) ,\label{eq49}\\
{\rm subject ~to} & & 
\lambda _{b}e^{- x} +\lambda _{s}e^{- y} \geq
b\eta _{b}+\lambda _{s} = s\eta _{s}+\lambda _{b} ,\label{eq50}
\end{eqnarray}
which is a convex optimization problem because the constraints in~(\ref{eq50}) form a convex set. The
optimality conditions for~(\ref{eq49})-(\ref{eq50}) are 
\[
\lambda _{b}s^{\alpha }x^{-\alpha }e^{-x} =\beta \lambda _{b} e^{
-x} ,
\]
\[
\lambda _{s}b^{\alpha }y^{-\alpha }e^{-y} =\beta \lambda _{s} e^{
-y} ,
\]
where $\beta $ is a Lagrange multiplier, which yields
\[
\frac{b}{y}=\frac{s}{x}=\beta^{1/\alpha }.
\]%
Hence, the change of variables in~(\ref{eq46a}) implies that 
\begin{equation}
\label{eq50aa}
v_{b}=v_{s}.
\end{equation}

Now define the change of variable 
\begin{equation}
\label{eq50a}
\tau = \frac{\kappa}{v_{b}^{1/\alpha }},
\end{equation}
so that $x=\tau s$ and $y=\tau b$. For a given $s\in (0,\rho_s)$, (\ref{eq48}) implies that 
\begin{equation}
\label{eq51}
b\left( s\right) =\frac{s\eta _{s}+\lambda _{b}-\lambda _{s}}{\eta
_{b}}.
\end{equation}
Substituting~(\ref{eq51}) into~(\ref{eq48}), we define $\tau \left( s\right) $ to be the unique solution to
\begin{equation}
\lambda _{b}e^{-s\tau} +\lambda _{s}e^{-b(s) \tau} =s\eta _{s}+\lambda _{b}
\label{eq52}
\end{equation}
for any $s\in (0,\rho _{s})$ and $b\in \left( 0,\rho _{s}\right) $. Note
that $\tau \left( s\right) $ is well defined because the left side of~(\ref{eq52}) is
decreasing in $\tau $, given $s$ and therefore $b\left( s\right)$.
Moreover, the left side of~(\ref{eq52}) is larger than the right side when $\tau =0$
and the right side is larger than the left side if $\tau =\infty $.

Then the optimization problem~(\ref{eq47})-(\ref{eq48}) takes the form
\begin{equation}
\label{eq53}
\max_{s\in \lbrack 0,\rho _{s}]} ~\lambda_s b(s)^\alpha G(b(s)\tau(s))+\lambda_b s^\alpha G(s\tau(s)).
\end{equation}
Let $H(s)=\lambda_s b(s)^\alpha G(b(s)\tau(s))+\lambda_b s^\alpha G(s\tau(s))$ from~(\ref{eq53}). The remainder of this proof is devoted to showing that $H(s)$ is concave on $s\in (0,\lambda_s/\eta_s)$, which implies that it has a unique maximizer. To prove that $H(s)$ is concave, it suffices to show that the second derivative $H''(s)\leq 0$.  
The first derivative of $H(s)$ is
\begin{equation}
	H'(s) = ~\lambda_s \alpha b^{\alpha-1}b' G(b\tau)+\lambda_b \alpha s^{\alpha-1} G(s\tau) +\tau^{-\alpha}[\lambda_s e^{-b\tau}(b\tau'+\tau b')+\lambda_be^{-s\tau}(s\tau'+\tau)]. \label{eq:F'_1} 
\end{equation}
By~(\ref{eq50a}), equation~(\ref{eq48}) can be expressed as 
\begin{equation}
\eta_ss+\lambda_b = \lambda_be^{-s\tau}+\lambda_se^{-b\tau}. 
\label{eq50b}
\end{equation}
Taking the derivative of (\ref{eq50b}) with respect to $s$ yields 
\begin{equation}
	\eta_s = -\lambda_be^{-s\tau}(s\tau'+\tau)-\lambda_s e^{-b\tau}(b\tau'+\tau b'),\label{eq:tau'}
\end{equation}
and substituting (\ref{eq:tau'}) into (\ref{eq:F'_1}) gives 
\begin{equation}
    \label{eq50c}
	H'(s) = \lambda_s \alpha b^{\alpha-1} G(b\tau)b'+\lambda_b \alpha s^{\alpha-1} G(s\tau) -\eta_s\tau^{-\alpha}.
\end{equation}
Taking the derivative of~(\ref{eq50c}) gives the second derivative
\begin{align}
	H''(s)=& ~\lambda_s \alpha(\alpha-1)b^{\alpha-2}G(b\tau)(b')^2+\lambda_b\alpha(\alpha-1)s^{\alpha-2}G(s\tau)+\lambda_s\alpha b^{\alpha-1}b'G'(b\tau)(b'\tau+\tau'b)\nonumber\\
	&+\lambda_b \alpha s^{\alpha-1} G'(s\tau)(\tau+s\tau')+\alpha\eta_s\tau^{-\alpha-1}\tau', \nonumber\\
	=& ~\alpha\lambda_s b^{\alpha-2}(b')^2[(\alpha-1 )G(b\tau)+(b\tau)^{1-\alpha}e^{-b\tau}] + \alpha\lambda_b s^{\alpha-2}[(\alpha-1 )G(s\tau)+(s\tau)^{1-\alpha}e^{-s\tau}]\nonumber\\
&+\alpha\tau^{-\alpha}\tau'(\lambda_se^{-b\tau}b'+\lambda_be^{-s\tau} +\eta_s\tau^{-1}),\label{eq:F''}
\end{align}
where the second equality is obtained by substituting the expression of $G'(\cdot).$ 

Solving~(\ref{eq:tau'}) for $\tau'$ gives 
\begin{eqnarray*}
	\tau' = -\frac{\eta_s +\lambda_be^{-s\tau}\tau+\lambda_sb'e^{-b\tau}\tau}{\lambda_b se^{-s\tau}+\lambda_s be^{-b\tau}}<0,
\end{eqnarray*}
which implies that the third term of (\ref{eq:F''}) is negative.
If we can show that the first two terms of (\ref{eq:F''}) are negative, i.e.,
\begin{eqnarray*}
	(\alpha-1 )G(b\tau)+(b\tau)^{1-\alpha}e^{-b\tau} \le 0,\\
	(\alpha-1 )G(s\tau)+(s\tau)^{1-\alpha}e^{-s\tau} \le 0,
\end{eqnarray*}
then we can conclude that $H''(s)< 0$, which means that there exists a unique maximizer of $H$.
To this end, define $$J(t) \triangleq (\alpha-1)G(t)+t^{1-\alpha}e^{-t},t\geq 0. $$
Because
\begin{align*}
	J'(t)=&(\alpha-1)G'(t)+(1-\alpha)t^{-\alpha}e^{-t}-t^{1-\alpha}e^{-t},\\
	=& -t^{1-\alpha}e^{-t}\leq 0,
\end{align*}
then for $t\geq 0$, $	J(t)\leq J(0)=0.$ Hence, we conclude that there exists a unique solution to the optimization problem and prove the asymptotic optimality of the proposed utility-threshold policy.

Taken together, the asymptotically optimal thresholds in the $n^{\rm th}$ system are 
\begin{equation}
\label{eq54}
v_{n,b}^* = v_{n,s}^* = \left(\frac{\kappa}{\tau(s_*)}\right)^{\alpha}m(n)
\end{equation}
by~(\ref{eq41a}), (\ref{eq50aa}) and~(\ref{eq50a}), where $s_*$ is the solution to~(\ref{eq53}) and $\tau(s_*)$ is the unique solution to \begin{equation}
\lambda _{b}e^{-s_*\tau} +\lambda _{s_*}e^{-b(s_*) \tau} =s_*\eta _{s_*}+\lambda _{b}.
\label{eq55}
\end{equation}
By~(\ref{eq72b}),~(\ref{eq:utilityrate_asymptotic}),~(\ref{eq46a}) and~(\ref{eq54}), the corresponding utility rate satisfies
\begin{equation}
\label{eq55a}
U_{n}^u(v_{n,b}^*,v_{n,s}^*) \sim n m(n)\kappa^{\alpha}[\lambda _{b}[s^*]^{\alpha }G(s^*\tau(s^*))+\lambda_s[b(s^*)]^{\alpha}G(b(s^*)\tau(s^*))]  ~~{\rm as} ~~n\rightarrow \infty.
\end{equation}



\subsection{Proof of Theorem~\ref{th:batch}}\label{ssec-proofs-batching}

Before proving Theorem~\ref{th:batch}, we introduce and prove the following lemma, which provides upper and lower bounds on the  limit of $E[{\cal M}(k)]$ as $k\to\infty$.


\begin{lemma}
\label{lemma:batch} The limit of $\frac{E[{\cal M}(k)]}{k^{\alpha+1}}$ exists as $k\rightarrow\infty$. Moreover, 
	\begin{eqnarray*}
		\left(1-\frac{3}{2}e^{-1/2} \right)^2\leq \lim_{k\to\infty}\frac{E[{\cal M}(k)]}{k^{\alpha+1}}\leq c\Gamma(1-\alpha).
	\end{eqnarray*}
\end{lemma}

\noindent {\bf Proof of Lemma~\ref{lemma:batch}}

First, note that $E{[\cal M(k)]}/k^{\alpha}$ is subadditive, which follows from the definition of ${\cal M(k)}$ and the fact that $k^\alpha$ is increasing. Fekete's Subadditive Lemma (Fekete 1923) implies that the limit of $E{[\cal M(k)]}/k^{\alpha+1}$ exists as $k\rightarrow\infty$.

Next, following Frenk {\sl et al.} (1987), we construct two sequences of i.i.d random variables, $U_{i,j}$ and $W_{i,j}$ for $i,j=1,\ldots,k$, with distribution $\tilde {F}(x)=\sqrt{F(x)}.$ Recalling that our matching utilities are denoted by $V_{i,j}$, we have that
\begin{eqnarray*}
	V_{i,j}&\overset{d}{=}& \max\{U_{i,j},W_{i,j}\},  ~~{\rm where}\\
	U_{i,j}&\overset{d}{=}& \tilde {F}^{-1}(X_{i,j}), ~~{\rm where} ~~X_{i,j}\sim \text{ Uniform}[0,1],\\
	W_{i,j}&\overset{d}{=}& \tilde {F}^{-1}(Y_{i,j}), ~~{\rm where} ~~Y_{i,j}\sim \text{ Uniform}[0,1].
\end{eqnarray*}

To obtain the lower bound, we consider a complete directed bipartite graph $G_k$ with vertices $B=\{b_1,...,b_k\},S=\{s_1,...,s_k\}$. There are directed edges $e(b_i,s_j)$ from vertex $b_i$ to $s_j$ and $e(s_j,b_i)$ from vertex $s_j$ to $b_i$. Each directed edge has a weight, with weight $U_{i,j}$ on $e(b_i,s_j)$ and $W_{i,j}$ on $e(s_j,b_i)$. Define $U_{i}^{(l)}$ to be the $l^{\rm th}$ largest weight among $U_{i,j},j=1,...,k$, and $W_{j}^{(l)}$ to be the $l^{\rm th}$ largest weight among $W_{i,j},i=1,...,k$.  Let $G_k(d)$ be the graph by removing edge $e(b_i,s_j)$ unless $U_{i,j}$ is one of the $d$ largest weights at $b_i$ and removing  edge $e(s_j,b_i)$ unless $W_{i,j}$ is one of the $d$ largest weights at $s_j$. Define $P(k,d)$ to be the probability that $G_k(d)$ contains a perfect matching. By Walkup (1980), we know that
\[1-P(k,2)\leq \frac{1}{5k},~~{\rm and} ~~ 1-P(k,d)\leq \frac{1}{122}\left(\frac{d}{k}\right)^{(d+1)(d-2)} ~~{\rm for} ~~d\geq 3.\]
Then
\begin{eqnarray}
	E[{\cal M}(k)] &=& P(k,2)E[{\cal M}(k)|G_k(2)\text{ contains a perfect matching}] \nonumber \\
	& & +(1-P(k,2))E[{\cal M}(k)|G_k(2)\text{ does not contain a perfect matching}]. \label{mp1}
\end{eqnarray}
The first expectation on the right side of~(\ref{mp1}) is lower bounded by
\begin{eqnarray*}
	E[{\cal M}(k)|G_k(2)\text{ contains a perfect matching}]&\geq& k E[\min \{U^{(2)},W^{(2)}\}], \nonumber \\ &=& kE[\min\{\tilde {F}^{-1}(X^{(2)}),\tilde {F}^{-1}(Y^{(2)})\}], \nonumber \\
	&=&kE[\tilde{F}^{-1}(\min\{X^{(2)},Y^{(2)}\})].
\end{eqnarray*}
Notice that
\begin{eqnarray}
	&&E[\tilde{F}^{-1}(\min\{X^{(2)},Y^{(2)}\})]\nonumber \\
	&=&E[\tilde{F}^{-1}(\min\{X^{(2)},Y^{(2)}\})I_{\left\{\min\{X^{(2)},Y^{(2)}\}\geq \sqrt{1-\frac{1}{k}}\right\}}]+E[\tilde{F}^{-1}(\min\{X^{(2)},Y^{(2)}\}) I_{\left\{\min\{X^{(2)},Y^{(2)}\}< \sqrt{1-\frac{1}{k}}\right\}}],\nonumber \\
	&\geq & P\left(\min\{X^{(2)},Y^{(2)}\}\geq  \sqrt{1-\frac{1}{k}}\right)F^{-1}\left(1-\frac{1}{k}\right).\label{mp3}\end{eqnarray}
The first term in~(\ref{mp3}) satisfies 
\begin{eqnarray}
	P\left(\min\{X^{(2)},Y^{(2)}\}\geq  \sqrt{1-\frac{1}{k}}\right)&=&P\left(X^{(2)}\geq  \sqrt{1-\frac{1}{k}},~Y^{(2)}\geq  \sqrt{1-\frac{1}{k}}\right),\nonumber \\
	&=&P\left(X^{(2)}\geq  \sqrt{1-\frac{1}{k}}\right)^2,\nonumber\\
	&=&\left(1- \left(1-\frac{1}{k}\right)^{k/2}-k \left(1-\frac{1}{k}\right)^{(k-1)/2}\left(1- \sqrt{1-\frac{1}{k}}\right)\right)^2,\nonumber\\
	&\rightarrow&\left(1-\frac{3}{2}e^{-1/2} \right)^2 ~\text{ as }k\rightarrow\infty.\label{mp4}
\end{eqnarray}
Equation~(\ref{mp4}) implies the lower bound 
\begin{eqnarray*}
	\lim\inf_{k\to\infty}\frac{{\cal M}(k)}{k^{\alpha+1}} &=& \lim\inf_{k\to\infty}\frac{{\cal M}(k)}{kF^{-1}(1-\frac{1}{k})}, \\
	&\geq& \left(1-\frac{3}{2}e^{-1/2} \right)^2.
\end{eqnarray*}

An upper bound for $E[{\cal M}(k)]$ is given by 
\begin{eqnarray*}
	E[{\cal M}(k)] &\leq& k\max\{U_1,\ldots,U_k\}, \\
	&=& ck\frac{\Gamma(k+1)\Gamma(1-\alpha)}{\Gamma(k+1-\alpha)}
\end{eqnarray*}
by equation~(3.7) in Malik (1966), where $c$ is a parameter of the Pareto distribution.
It follows that 
\begin{eqnarray*}
	\lim\sup_{k\to\infty} \frac{E[{\cal M}(k)]}{k^{\alpha+1}}\leq \lim_{k\to\infty} c\frac{\Gamma(k+1)\Gamma(1-\alpha)}{k^\alpha\Gamma(k+1-\alpha)}.
\end{eqnarray*}
The proof is completed by noting that 
\begin{eqnarray*}
	\lim_{k\to\infty}c\frac{\Gamma(k+1)\Gamma(1-\alpha)}{k^\alpha\Gamma(k+1-\alpha)}=c\Gamma(1-\alpha),
\end{eqnarray*}
which follows from equation~(1) of Tricomi and Erd$\acute {\rm e}$lyi (1951).

\noindent {\bf Proof of Theorem~\ref{th:batch}}

As before, we let $(B_n(t),S_n(t))$ be the number of buyers and sellers in the $n^{\rm th}$ system at time $t$. To analyze the long-run average performance of the batch-and-match policy, it suffices to study the queue length dynamics over a single batching cycle, which for ease of presentation we take to be the time interval $[0,\Delta)$, where by construction (i.e, because matches were made just prior to time 0) $\min\{B_n(0),S_n(0)\}=0$. For $t\in [0,\Delta)$, because there are no matches made in this time interval, we have 
\begin{eqnarray*}
	B_n(t)&=&B_n(0)+N_B^+(n\lambda t)-N_B^-(\eta \int_0^t B_n(r)dr),\\
	S_n(t)&=&S_n(0)+N_S^+(n\lambda t)-N_S^-(\eta \int_0^t S_n(r)dr),
\end{eqnarray*}
with its fluid model
\begin{eqnarray*}
	\bar{B}(t)&=&\bar{B}(0)+\lambda t -\eta \int_0^t \bar{B}(r)dr,\\
	\bar{S}(t)&=&\bar{S}(0)+\lambda t -\eta \int_0^t \bar{S}(r)dr.
\end{eqnarray*}
We know that $(\bar{B}_n(\cdot),\bar{S}_n(\cdot))$ converges to $(\bar{B}(\cdot),\bar{S}(\cdot))$ uniformly on compact sets in probability, where
\begin{eqnarray}
	\bar{B}(t) &=& \bar{B}(0)e^{-\eta t}+\frac{\lambda}{\eta} (1- e^{-\eta t}), \label{batch3}\\
	\bar{S}(t) &=& \bar{S}(0)e^{-\eta t}+\frac{\lambda}{\eta} (1- e^{-\eta t}). \label{batch4}
\end{eqnarray}
By the continuous mapping theorem, 
\begin{eqnarray*}
	\min(\bar{B}_n(t),\bar{S}_n(t))\rightarrow \min (\bar{B}(t),\bar{S}(t))\text{ in probability as} ~~n\to\infty.
\end{eqnarray*}
Because $(B_n(t),S_n(t))$ is upper bounded by the queue lengths with only arrivals and no abandonments, the dominated convergence theorem implies that 
\[
	\lim_{n\to\infty} E[\min\{\bar{B}_n(t),\bar{S}_n(t)\}] =E[\min\{\bar{B}(t),\bar{S}(t)\}].
\]
Because $\min\{B_n(0),S_n(0)\}=0$, it follows from~(\ref{batch3})-(\ref{batch4}) that just prior to matching, 
\begin{equation}
\label{batch6}
\min\{\bar{B}(\Delta),\bar{S}(\Delta)\}= \frac{\lambda}{\eta}(1-e^{-\eta \Delta}).
\end{equation}

By Lemma$~\ref{lemma:batch}$, there exists a constant $C$ such that 
\[
\lim_{n\to\infty} \frac{E[{\cal M}(k)]}{k^{\alpha+1}} = C \in \left[\left(1-\frac{3}{2}e^{-1/2} \right)^2, c\Gamma(1-\alpha)\right].
\]
Hence, $E[{\cal M}(k)]$ is regularly varying with index $\alpha$ by Assumption~\ref{as:Regularly_Varying}, and Lemma~\ref{lemma:Expectation_Max} implies that the utility rate $U_n^b(\Delta)$ of the batch-and-match policy with time window $\Delta$ satisfies
\begin{equation}
    \label{batch1}
\lim_{n\to\infty}\frac{U_n^b(\Delta)}{n^{\alpha+1}} = \frac{C\Bigl[\frac{\lambda}{\eta}(1-e^{-\eta \Delta})\Bigr]^{\alpha+1}}{\Delta}.
\end{equation}
Because the right side of~(\ref{batch1}) is concave, the first-order conditions corresponding to~(\ref{batch1}) imply that the asymptotically optimal time window $\Delta^*$ is the unique positive solution to~(\ref{batch2}) in the main text,
which is independent of the arrival rate $\lambda$ and decreasing in the abandonment rate $\eta$. Combining Lemma~\ref{lemma:batch}, equation~(\ref{batch2}) in the main text,  and~(\ref{batch1}) gives the upper bound in~(\ref{batch5}) in the main text, thereby concluding the proof of Theorem~\ref{th:batch}.

\noindent{\bf The Unbalanced Case}

For brevity's sake, we present the corresponding results for the unbalanced case and omit the details. The arrival rates are $(\lambda_b,\lambda_s)$, the abandonment rates are $(\eta_b,\eta_s)$ and $\rho_b=\lambda_b/\eta_b$, $\rho_s=\lambda_s/\eta_s$. The generalization of~(\ref{batch6}) is
\[
\min\{\bar{B}(t),\bar{S}(t)\}=\min\{\rho_b (1- e^{-\eta_b t}),\rho_s (1- e^{-\eta_s t})\},
\]
and the utility rate satisfies 
\begin{equation}
\label{batch7}
\lim_{n\to\infty}\frac{U_n^b(\Delta)}{n^{\alpha+1}} = \frac{C[\min\{\rho_b (1- e^{-\eta_b t}),\rho_s (1- e^{-\eta_s t})\}]^{\alpha+1}}{\Delta}.
\end{equation}
The asymptotically optimal time window $\Delta^*$ solves
\begin{eqnarray*}
	\max_{\Delta,\xi }&& \frac{\xi^{\alpha+1}}{\Delta}\\
	\text{s.t.} &&\xi\leq \frac{\lambda_b}{\eta_b}(1-e^{-\eta_b \Delta}),\\
	&& \xi\leq \frac{\lambda_s}{\eta_s}(1-e^{-\eta_s \Delta}).
\end{eqnarray*}
This optimization problem has a solution because it has a concave objective function and a convex feasible set.

\section{Examples} \label{sec-examples} 

In~\S\ref{ssec-examples-exp}-\ref{ssec-examples-uniform}, we consider one canonical matching utility distribution from
each of the three domains of attraction (Weibull, Gumbel and Frechet), represented, respectively, by U($a,b)$, exp($\nu$), and Pareto($c,\beta$). In each of these examples, we compute the utility rate of the upper bound in Lemma~\ref{lemma:bound}, the utility rate under the greedy policy from~(\ref{eq16}) in the main text, and the asymptotically optimal (or heuristic, in some cases) thresholds and corresponding utility rates for the population threshold policy and the utility threshold policy from Theorems~\ref{thm:Population_Based} and~\ref{thm:utility_based}, respectively. We continue to add the superscripts $+$, $g$, $p$ and $u$ to $U$ to denote the utility rate of the upper bound, the greedy policy, the population threshold policy and the utility threshold policy, respectively. We briefly consider matching utilities that come from a correlated Pareto distribution in~\S\ref{ssec-examples-correlated}. 

\subsection{Matching Utilities are Exponential}

\label{ssec-examples-exp} Let the matching utilities be exponential with
parameter $\nu$ and CDF $F(v)=1-e^{-\nu v}$ for $v\ge 0$, which falls under the $\alpha=0$ case in Theorem~\ref{thm:Population_Based}. The exponential is
in the domain of attraction of Type I, and so (see~(\ref{eq10})-(\ref{eq10a}) in~\S\ref{sec-evt}) $\mu=\gamma=0.5772\ldots$,
which is Euler's constant, $a_n=\frac{1}{\nu}$ and $b_n=\frac{\ln n}{\nu}$,
and hence $m(n)\sim (\gamma+\ln n)/\nu$. 

The utility rate of the upper bound is 
\begin{eqnarray}
U_n^+ & \sim & \frac{\lambda n}{\nu}\Biggl(\gamma+\ln\Biggl(\frac{n\lambda}{%
\eta}\Biggr)\Biggr) ~~{\rm by} ~{\rm Lemma}~\ref{lemma:bound}, \label{eq22a} \\
& \sim & \frac{\lambda}{\nu} n\ln n,  \label{eq50-mt}
\end{eqnarray}
and the utility rate of the greedy policy is 
\begin{eqnarray}
U_n^g & \sim & \frac{\lambda n}{\nu}\Biggl(\gamma+\ln\Biggl(\frac{\lambda}{%
\eta}\sqrt{\frac{2n}{\pi}}\Biggr)\Biggr) ~~{\rm by} ~~(\ref{eq16})~~{\rm in~the~main~text}, \label{eq23a} \\
& \sim & \frac{\lambda}{2\nu} n\ln n.  \label{eq51-mt}
\end{eqnarray}

As noted below Theorem~\ref{thm:Population_Based}, a range of population thresholds are asymptotically optimal in the $\alpha=0$ case. For concreteness, we consider $z_n^*=\frac{n}{\ln n}$, which has utility rate 
\begin{eqnarray}
U_n^p\Biggl(\frac{n}{\ln n}\Biggr) & \sim & \frac{\lambda n}{\nu}\Biggl(%
\gamma+\ln\Biggl(\frac{n}{\ln n}\Biggr)\Biggr) ~~{\rm by} ~(\ref{eq4b})~~{\rm in~the~main~text},  \label{eq51a} \\
& \sim & \frac{\lambda}{\nu} n\ln n.  \label{eq52-mt}
\end{eqnarray}
By (\ref{eq50-mt}),~(\ref{eq51-mt}) and~(\ref{eq52-mt}), the population threshold policy with
threshold $\frac{n}{\ln n}$ is asymptotically optimal and doubles the utility rate of the greedy
policy in the  limit.



Recall that Theorem~\ref{thm:utility_based} does not apply to the $\alpha=0$ case. Nonetheless, we apply the ideas in Theorems~\ref{thm:Population_Based} and~\ref{thm:utility_based} to derive a heuristic threshold level for the utility threshold policy. By considering the steady-state version of equation \eqref{eq:BnUtilityBased}  and differentiating, we obtain 
$$\frac{\eta B_n(\infty)}{\lambda n} = P( M(S_n(\infty))  \le v_n ) - \left(1  -  P( M(B_n(\infty))  \le v_n )\right),$$
which by symmetry yields
\begin{equation}\label{eq:Fluid_compute}
\frac{\eta B_n(\infty)}{\lambda n} = 2 P( M(B_n(\infty))  \le v_n ) - 1. 
\end{equation}
Now we heuristically assume that the utility threshold $v_n$ is such that it achieves a population level $B_n(\infty)$ that equals the optimal population threshold $z_n^*$, which for concreteness we again take to be $\frac{n}{\ln n}$. Substituting $\frac{n}{\ln n}$ for $B_n(\infty)$ in~(\ref{eq:Fluid_compute}) and noting that $P( M(n)  \le v_n ) = P( \nu M(n) - \ln n \le \nu v_n - \ln n ) \sim \exp(-\exp( -\nu v_n +\ln n ))$ by~(\ref{evt1}), we get
\begin{equation}
\frac{\eta}{\lambda\ln n} = 2 \exp\left[-\exp\left(-\nu v_n+\ln\left(\frac{n}{\ln n}\right)\right)\right]-1.
\label{eq64}
\end{equation}
Solving equation~(\ref{eq64}) gives the proposed threshold level, 
\begin{equation}
v_n^* = \frac{\ln n - \ln \ln n - \ln \ln \left(\frac{2\lambda \ln n}{\lambda \ln n+\eta}\right)}{\nu}.
\label{eq65}
\end{equation}
This heuristic approach does not generate a corresponding utility rate. 

\subsection{Matching Utilities are Pareto With Finite Mean}

\label{ssec-examples-pareto} Let the matching utilities have CDF 
$F(v)=1-(cv)^{-\beta}$, for $\beta>1, c>0$ and $cv\ge 1$, so that the mean matching
utility is finite and $\alpha=1/\beta$ in Theorem~\ref{thm:Population_Based}. The Pareto distribution is in the domain of attraction of the Frechet distribution, and hence (see~(\ref{eq10})-(\ref{eq10a}) in~\S\ref{sec-evt}) $b_n=0$, $a_n=(cn)^{1/\beta}$ and $\mu=\Gamma\Bigl(1-%
\frac{1}{\beta}\Bigr)$, where $\Gamma(n)$ is the gamma function. It follows
that 
\[
m(n)\sim (cn)^{1/\beta}\Gamma\Bigl(1-\frac{1}{\beta}\Bigr).
\]
By Lemma~\ref{lemma:bound} and~(\ref{eq16}) in the main text, 
\begin{equation}  \label{eq54a}
U_n^+ \sim \lambda\Bigl(\frac{c\lambda}{\eta}\Bigr)^{1/\beta}\Gamma\Bigl(1-%
\frac{1}{\beta}\Bigr)n^{1+1/\beta},
\end{equation}
and 
\begin{equation}  \label{eq55a-mt}
U_n^g \sim \lambda\Bigl(\frac{c\lambda}{\eta}\sqrt{\frac{2}{\pi}}\Bigr)%
^{1/\beta}\Gamma\Bigl(1-\frac{1}{\beta}\Bigr)n^{1+1/(2\beta)}.
\end{equation}
In contrast to~(\ref{eq50-mt}) and~(\ref{eq51-mt}) in the exponential case and to~(\ref{eq60}) and~(\ref{eq61}) in the uniform case, the upper bound and the greedy performance in~(\ref{eq54a})-(\ref{eq55a-mt}) have different growth rates in $n$. 

Because $\alpha=1/\beta$, part ii) of Theorem~\ref{thm:Population_Based} implies that 
\begin{equation}  \label{eq56}
z_n^*=\frac{\lambda}{\eta(1+\beta)}n.
\end{equation}
That is, the optimal threshold equals the mean size of either side of the market in the absence
of matching ($\lambda n/\eta$) times the factor $\frac{1}{1+\beta}$, which
is less than 1/2. Substituting~(\ref{eq56}) into~(\ref{eq5b}) in the main text gives the
utility rate 
\begin{equation}  \label{eq57}
U_n^p(z_n^*) \sim \lambda \Bigl(\frac{%
c\lambda }{\eta(1+\beta)}\Bigr)^{1/\beta}\left(\frac{\beta}{1+\beta}\right)\Gamma\Bigl(1-\frac{1}{\beta}\Bigr)%
n^{1+1/\beta},
\end{equation}
which has the same exponent of $n$ as the upper bound. Comparing the utility rate under the optimal population threshold policy to the upper
bound, we get 
\begin{equation}  \label{eq58}
\frac{ U_n^+}{U_n^p(z_n^*)} = \frac{1}{\left(\frac{1}{1+\beta}\right)^{1/\beta}
\left(\frac{\beta}{1+\beta}\right)},
\end{equation}
which converges to 4 as $\beta\to 1$, and
converges to 1 as $\beta\to\infty$.

Comparing the utility rate of the optimal population threshold policy to the utility
rate of the greedy policy yields 
\begin{equation}
\frac{U_n^p(z_n^*)}{U_n^g} = \left(\frac{\sqrt{\pi}}{\sqrt{2}(1+\beta)}\right)^{1/\beta}
\left(\frac{\beta}{1+\beta}\right)
n^{1/(2\beta)}, 
\label{eq58a}
\end{equation}
which converges to $\sqrt{\frac{\pi n}{32}}\approx 0.3133\sqrt{n}$ as $%
\beta\to 1$, and converges to $1$ as $\beta\to\infty$. 
For all finite values of $\beta$, the difference in performance between the two policies
becomes unbounded as $n\to\infty$.

The optimal utility threshold needs to be computed numerically using the results in Theorem~\ref{thm:utility_based}. To streamline the presentation, we consider the special case considered in the simulation experiments in~\S\ref{sec-simulation}, where $\lambda=\eta=1$, $c=1$ and $\beta=2$, and hence $\alpha=0.5$, $\kappa=1/\pi$ and $m(n)=\sqrt{\pi n}$. Using the fact that the incomplete gamma function $\gamma(0.5,x)=\sqrt{\pi}$erf$(\sqrt{x})$, by Theorem~\ref{thm:utility_based} we need to find the solution $z_*\in (0,\lambda/\eta)$ satisfying 
\begin{equation}
\label{eq58b}
\frac{z}{\sqrt{\pi \ln\left(\frac{2}{z+1}\right)}} = {\rm erf}\left(\sqrt{\ln\Bigl(\frac{2}{z+1}\Bigr)}\right).
\end{equation}
Given the solution $z_*$ to~(\ref{eq58b}), Theorem~\ref{thm:utility_based} and its proof imply that 
\begin{equation}
\label{eq58c}
v_n^* = \sqrt{\frac{nz_*}{\ln\left(\frac{2}{z_*+1}\right)}}
\end{equation}
and
\begin{equation}
\label{eq58d}
U^u_n(v_n^*) \sim \frac{2(nz_*)^{3/2}}{\sqrt{\ln\left(\frac{2}{z_*+1}\right)}}.
\end{equation}

\subsection{Matching Utilities are Uniform}

\label{ssec-examples-uniform} When the matching utilities are distributed as
U($a,b$) with $F(v)=\frac{v-a}{b-a}$ for $v\in [a,b]$, which is in the domain of attraction of the Weibull law, we have (see~(\ref{eq10})-(\ref{eq10a}) in~\S\ref{sec-evt}) $a_n=%
\frac{b-a}{n}$, $b_n=b$, $\mu=-\Gamma(2)=-1$, $m(n)\sim b-\frac{%
b-a}{n}$, and $\alpha=0$ in Theorem~\ref{thm:Population_Based}. By Lemma~\ref{lemma:bound} and~(\ref{eq16}) in the main text, 
\begin{eqnarray}
U_n^+ & \sim & n\lambda\Bigl(b-\frac{b-a}{\frac{\lambda}{\eta}n}\Bigr),
\label{eq33a} \\
& \sim & \lambda bn,  \label{eq60}
\end{eqnarray}
and 
\begin{eqnarray}
U_n^g & \sim & n\lambda\Bigl(b-\frac{b-a}{\frac{\lambda}{\eta}\sqrt{\frac{2n}{%
\pi}}}\Bigr),  \label{eq34a} \\
& \sim & \lambda bn.  \label{eq61}
\end{eqnarray}
By~(\ref{eq60}) and~(\ref{eq61}), the greedy policy is asymptotically optimal,
and so there is no need to consider a positive
threshold level for the population threshold policy. 

To heuristically analyze the utility threshold policy, we proceed as in~\S\ref{ssec-examples-exp}, where equation~(\ref{eq:Fluid_compute}) now becomes
$$
\frac{\eta B_n(\infty)}{\lambda n} = 2 \left(\frac{v_n-a}{b-a}\right)^{B_n(\infty)} - 1, 
$$
which can be rearranged as
$$
v_n = a +(b-a)\left(\frac{\eta B_n(\infty)}{2\lambda n}+\frac{1}{2}\right)^{1/B_n(\infty)}.
$$
Setting $B_n(\infty)=\frac{\lambda}{\eta}\sqrt{\frac{2n}{\pi}}$ from~(\ref{eq15}) in the main text, which is the expected number of available mates for an arriving agent under the greedy policy, leads to
\begin{equation}
\label{eq61b}
v_n^* = a + (b-a)\left(\frac{1}{\sqrt{2n\pi}}+\frac{1}{2}\right)^{\frac{\eta}{\lambda}\sqrt{\frac{\pi}{2n}}}.
\end{equation}

\subsection{Matching Utilities are Correlated}

\label{ssec-examples-correlated}




One of the advantages of our analysis, which essentially decouples the extremal behavior of the utilities and the dynamics of the agents in the fluid scale, is that we can enrich our model with complex dependencies in the utilities by directly importing results from extreme value theory for non-i.i.d. sequences of random variables. Assumptions 1, 2 and 3 can be shown to hold in substantial generality, well beyond the setting of i.i.d. utilities. The study and calibration of extremes under non-i.i.d. sequences is a well-developed topic in extreme value theory; e.g., see Leadbetter et al. (1983) and Smith and Weissman (1994). 

Our decoupling approach provides a great degree of flexibility for modelers -- informed by specific types of applications -- to incorporate correlated utilities. We consider two examples here. The first example involves utilities that can be decomposed by adding a common factor (or factors) and an idiosyncratic factor, thereby inducing a correlation effect.  Suppose that when a seller (buyer) arrives and finds $k$ buyers (sellers), the corresponding utilities $V_1, V_2, \ldots, V_k$ are given by 
\begin{equation}\label{eq:max_model}
V_i = \rho U_0 + \sqrt{1-\rho^2} U_i ~~{\rm for} ~ i=1,2, \ldots,k,
\end{equation} 
where $\rho\in [0,1)$, and $U_0,U_1,\ldots,U_k$ are i.i.d.\ with a Pareto($\frac{\sqrt{3}}{2},3$) distribution; i.e., for each $i$ we have $P(U_i \le u) = 1 - (\frac{2}{\sqrt{3}u})^3$ for $u \ge \frac{2}{\sqrt{3}}$, and $0$ otherwise. For example, the model in~(\ref{eq:max_model}) allows the common portion ($\rho U_0$) of the utility to quantify some characteristics of the seller, and the variable portion ($\sqrt{1-\rho^2}U_i$) to vary across the $k$ buyers based on their individual characteristics. The variance of each $V_i$ is $1$, independent of $\rho$. However, for $\rho>0$, Cov$(V_i,V_j) = \rho^2$ for each $i\neq j$, and the utilities have correlation $\rho^2$. 

Before computing the optimal thresholds and the corresponding utility rates for the population threshold policy and the utility threshold policy, we show that Assumptions~\ref{as:Regularly_Varying}-\ref{as:Max_Pareto} hold. We have $m(k) = \rho E[U_0] + \sqrt{1-\rho^2} E[\max_{i=1}^k U_i]  \sim \sqrt{1-\rho^2} \Gamma(1-\frac{1}{\beta}) k^{1/\beta} $, and  $\frac{M(k)}{m(k)} = \frac{ \rho U_0 + \sqrt{1-\rho^2} \max_{i=1}^k U_i}{m(k)} \sim   \frac{  \max_{i=1}^k U_i}{ E[\max_{i=1}^k U_i] } $ w.p.1, which implies that Assumptions \ref{as:Regularly_Varying} and \ref{as:Max_Pareto} are satisfied with $\alpha = 1/3$. 

By Theorem~\ref{thm:Population_Based}, the asymptotically optimal population threshold is $z_{n}^*=\frac{\lambda \alpha}{\eta (1+\alpha)}n$, independent of $\rho$, and the asymptotic utility rate is proportional to $\sqrt{1-\rho^2}$. By Theorem \ref{thm:utility_based}, the asymptotically optimal utility threshold is  $v_n^* = \sqrt{1-\rho^2} \Gamma(1-\frac{1}{\beta})  v_{\ast } n^{1/\beta} $ where $v_{\ast }$ is as given in the theorem statement, and the corresponding utility rate is increasing in $\sqrt{1-\rho^2}$. 

Hence, in both cases, the utility rate decreases as the correlation $\rho^2$ increases, as does the optimal utility threshold. This is to be expected because the increased utility from being patient is reduced in the presence of positive correlation. Less obvious is that the optimal population threshold is independent of $\rho^2$. The optimal population threshold is independent of $\rho^2$ because it trades off the higher utility from additional thickness (i.e., increased $z_n$), which depends on the variation in $(U_1,\ldots,U_k)$ in~(\ref{eq:max_model}) but not on the common term $U_0$ that dictates the correlation, and the higher abandonment rate, which is independent of $\rho^2$. 

Our second example involves cases in which the population size itself impacts the correlation among utilities (e.g., crowding network effects). Suppose the matching utilities for an arriving agent that observes $k$ potential matches are 
$$
V_i = 2^{-1/\beta} \max(c_k W,U_i) ~~{\rm for} ~i=1,2, \ldots,k,
$$
where $U_1,\ldots,U_k$ are i.i.d. with a Pareto($1,\beta$) distribution with $\beta>1$, $W$ is  independent of $U_i$ and is a Frechet($\beta$) distributed random variable that is drawn independently for each arriving agent, and $c_k = k^{1/\beta} \Gamma(1-\frac{1}{\beta})$ for each $k$. Assumptions \ref{as:Regularly_Varying} and \ref{as:Max_Pareto} are satisfied with $\alpha = 1/\beta$, and thus our results apply. 

The random variable $W$ could arise as a result of a separate mechanism in which crowding is incorporated. For instance, if $W$ itself was the result of another selection process that filtered customer arrivals (thus the assumption that $W$ is Frechet). However, even when $W$ has a distribution that is different than Frechet($\beta$), while Assumption \ref{as:Max_Pareto} may not hold as is, our proof techniques may be leveraged to compute optimal thresholds. We omit the details for the sake of brevity.

\section{Extreme Value Theory and Regularly Varying Functions}\label{sec-evt}

In this section, we collect some useful facts about extreme value theory and show that Assumption~\ref{as:Regularly_Varying} is satisfied by distributions that are subject
to the application of extreme value theory. Throughout this section we assume that the $V_{i}$s are i.i.d. random variables.

The central result in extreme value theory is that, for certain
distributions $F(v)$, the CDF of a properly normalized version of $M\left(
n\right) $ converges to a limiting CDF that is known as the generalized
Pareto distribution. More precisely, 
\begin{equation}
\label{evt1}
P\left( \frac{M\left( n\right) -b_{n}}{a_{n}}\leq x\right) \rightarrow
\Xi(ax+b;\xi)~~\mathrm{as}~~n\rightarrow \infty , 
\end{equation}
where 
\[
\Xi(x;\xi ):=\exp \left( -\left( 1+\xi x\right) ^{-1/\xi }\right)
,\quad \quad 1+\xi x>0 
\]%
and $\xi \in \mathcal{R}$. The case $\xi =0$ is interpreted as $%
G(x;0)=e^{-e^{-x}} $.

There are three domains of attraction: $\xi <0$ (Weibull), $\xi =0$ (the Gumbel), and $\xi >0$ (the Frechet). 
Distributions with bounded support (e.g., uniform, beta) typically\ belong to the Weibull domain of attraction. Distributions with finite moments of every order (often, but not always, with unbounded support) belong to the domain of attraction of the Gumbel distribution (e.g., exponential, gamma, normal, lognormal). Distributions with power-law-like decaying tails belong to the domain of attraction of the Frechet distribution (e.g., Pareto, Cauchy). The CDFs corresponding to the Weibull, Gumbel and Frechet domains of attraction are readily available by evaluating the corresponding values of $\xi$ in $G(x;\xi)$. 



Define $\bar{F}\left( x\right) =1-F\left( x\right) $ and let $w(F)=\sup \{x:F(x)<1\}$ be the upper endpoint of the support of $F$. If convergence to a generalized Pareto distribution with parameter $\xi$ holds (i.e. if extreme value theory applies), then the corresponding constants can be computed as follows:
\begin{equation}
{\rm Weibull} (\xi<0):~a_{n}=w(F)-\bar{F}^{-1}(n^{-1}),~b_{n}=w(F),~a=-1/\xi,~\mathrm{and}~b=1/\xi, \label{eq9}
\end{equation}%
\begin{equation}
{\rm Gumbel} (\xi=0):~a_{n}=\frac{1}{\bar{F}\left( b_{n}\right) }\int_{b_{n}}^{w(F)}\bar{F}%
\left( t\right) dt,~b_{n}=\bar{F}^{-1}(n^{-1}),~a=0,~\mathrm{and}~b=0, \label{eq9b}   
\end{equation}
\begin{equation}
{\rm Frechet} (\xi>0):~a_{n}=\bar{F}^{-1}(n^{-1}),~b_{n}=0,~a=1/\xi,~\mathrm{and}~b=-1/\xi,  \label{eq8}
\end{equation}%
where $\bar{F}^{-1}(\cdot )$ is the inverse of $\bar{F}\left( \cdot \right) $. 

We now switch our attention to $E[M\left( n\right) ]$. Theorem~2.1 in Pickands
(1968) shows that the first moment converges as long as $E[V^{1+\delta
}]<\infty $ for some $\delta >0$, which is satisfied for all the concrete
examples explored in this paper. This result implies that 
\begin{equation}
E[M\left( n\right) ]\sim b_{n}+a_{n}\mu~, \label{eq10}
\end{equation}%
where 
\begin{equation}
\mu=\int_{-\infty }^{\infty }x~d\Xi(ax+b;\xi)~ \mathrm{for}~ \xi \in (-\infty,\infty)
\label{eq10a}
\end{equation}%
is the mean of the distribution $\Xi(\cdot;\xi)$ computed according to the limiting value of the domain of attraction. 
The mean
in~(\ref{eq10a}) is $\mu=\gamma =0.5772\ldots$ if $\xi=0$, which is Euler's
constant; $\mu=\Gamma (1-\xi)$, where $\Gamma
(a)$ is the gamma function evaluated at $a>0$, if $\xi \in (0,1)$, which holds in the Frechet case if $E(V^{1+\delta})<\infty$; and $\mu=-\Gamma(1-\xi)$ if $\xi<0$.

We conclude this section by stating the relationship between extreme value theory and slowly varying functions, thereby showing that the extreme value distributions satisfy Assumption~\ref{as:Regularly_Varying}.

\textbf{Fact~1.} If $F$ belongs to the domain of attraction of the Gumbel distribution (i.e. $\xi=0$) then $b_n$ is slowly varying at infinity.

Fact~1 follows from Proposition 0.10 of Resnick (1987), combined
with the first exercise on p. 35 of Resnick (1987). 

\textbf{Fact~2.} (Resnick 1987, p. 52) If $F$ belongs to the domain of the Gumbel law then $a_{n}=o\left( b_{n}\right) $ as $%
n\rightarrow \infty $.

\textbf{Fact~3.} (Resnick 1987, p. 54) $F$ belongs to the domain of
attraction of the Frechet distribution (i.e. $\xi>0$) if and only if $a_{n}=0$ and $b_{n}=\bar{F}^{-1}\left( 1/n\right) $ as $n\rightarrow \infty $ and $F\left( \cdot
\right) $ is regularly varying with index $-\xi$.

\textbf{Fact~4.} (Resnick 1987, p. 59) $F$ belongs to the domain of
attraction of a Weibull distribution if and only if $w\left( F\right) <\infty $ and $%
\bar{F}\left( w\left( F\right) -x^{-1}\right) $ is regularly varying with
index $\xi $.

Facts~1 and~2 imply that $m\left( \cdot \right) $ is slowly
varying (i.e. regularly varying with index 0) for distributions in the Gumbel domain of attraction.
Fact~3 implies that $m\left( \cdot \right) $ is regularly varying with
index $\xi \in \left( 0,1\right) $ for distributions in the domain of attraction of the Frechet law. Finally, Fact~4 implies that $m\left(
\cdot \right) $ is slowly varying for distributions in the domain of attraction of the Weibull law.

\clearpage

\begin{figure}[hbt]
\centering
  \includegraphics[width=0.8\textwidth]{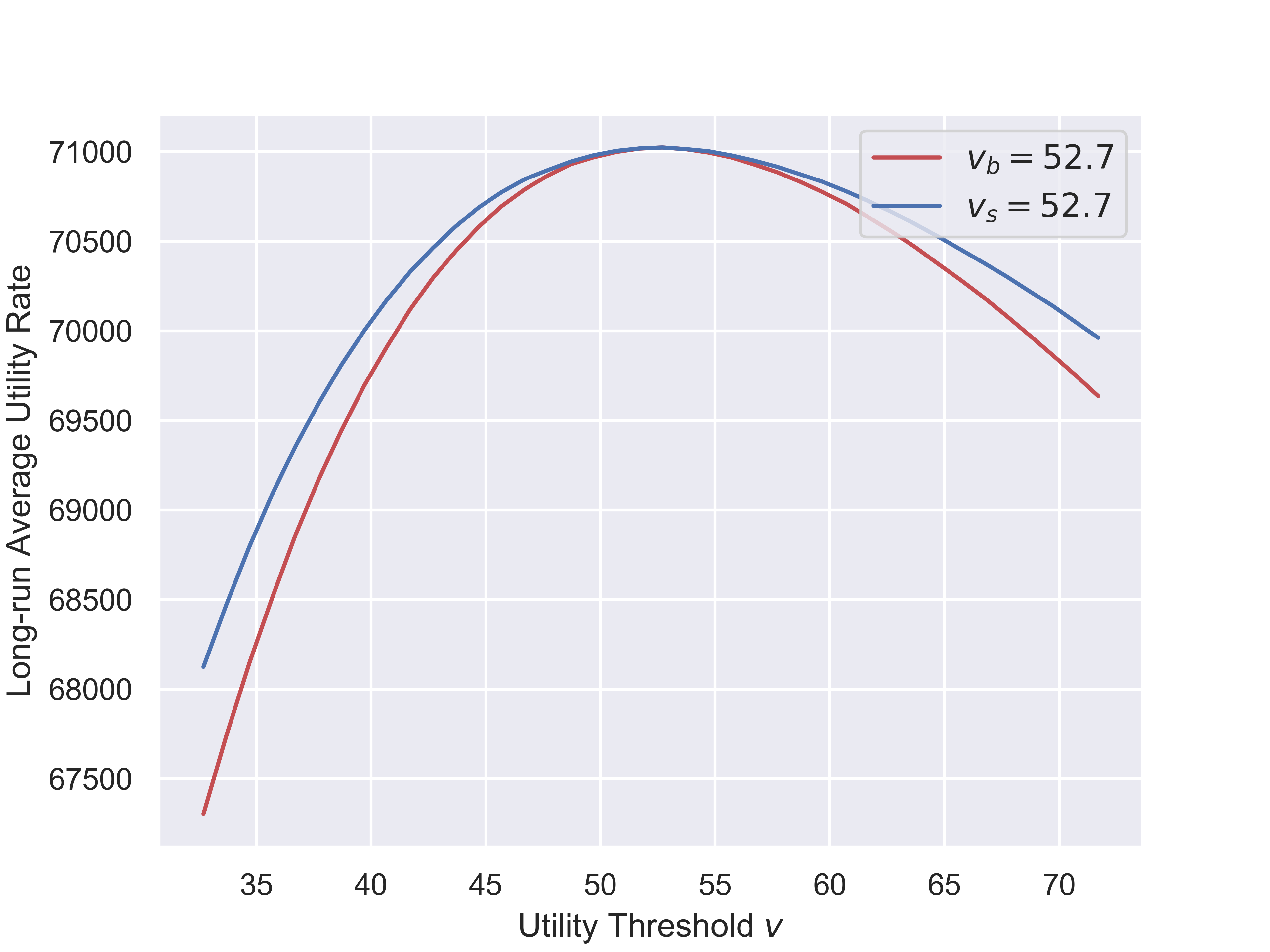}\label{Fig:sensitivity}
\end{figure}
\noindent{\bf Figure~1:} In the unbalanced numerical example in~\S\ref{sec-unbalanced}, the simulated utility rate vs. the utility threshold $v_b$ with $v_s=52.7$, and vs. $v_s$ with $v_b=52.7.$






\end{document}